%% file: sigma18-084.tex
\pdfoutput=1
\RequirePackage{ifpdf}
\ifpdf 
\documentclass[pdftex]{sigma}
\else
\documentclass{sigma}
\fi

\usepackage{enumitem}
\usepackage[mathscr]{eucal}
\usepackage{tensor}
\usepackage{wrapfig} 
\usepackage[all]{xy}
\usepackage{float} 
\usepackage{longtable}

\usepackage{mathtools} 
\usepackage{accents}
\usepackage[normalem]{ulem}
\usepackage{cancel}
\usepackage{upgreek}

\input xy
\xyoption{all}

\numberwithin{equation}{section}
\numberwithin{figure}{section}

\newtheorem{Theorem}{Theorem}[section]
\newtheorem{Corollary}[Theorem]{Corollary}
\newtheorem{Lemma}[Theorem]{Lemma}
\newtheorem{Proposition}[Theorem]{Proposition}
\newtheorem*{assumption}{Assumption}

\theoremstyle{definition}
\newtheorem{Definition}[Theorem]{Definition}
\newtheorem{Example}[Theorem]{Example}
\newtheorem{exm}[Theorem]{Example}

\newtheorem{rmk}[Theorem]{Remark}

\newcommand{\C}{\mathbb{C}}

\newcommand{\N}{\mathbb{N}}

\newcommand{\R}{\mathbb{R}}
\newcommand{\T}{\mathbb{T}} 
\newcommand{\Z}{\mathbb{Z}}

\def\epsilon{\varepsilon}

\newcommand{\ga}{\gamma}

\newcommand{\ep}{\epsilon}

\newcommand{\et}{\eta}

\newcommand{\si}{\sigma}

\newcommand{\om}{\omega}

\newcommand{\Ga}{\Gamma}

\newcommand{\etti}{{\widetilde{\et}}}

\newcommand{\zti}{{\widetilde{z}}}

\newcommand{\sm}{(M,\omega)} 
\newcommand{\cotan}{\mathrm{T}^*} 

\newcommand{\toric}{(M,\omega,\mu)} 
\newcommand{\toricphi}{(M,\omega,\Phi)} 

\newcommand{\is}{(M,\omega,\Phi)} 
\newcommand{\vat}{(M,\omega,\Phi=(J,H))} 

\DeclarePairedDelimiter\abs{\lvert}{\rvert}

\newcommand{\vungoc}{V\~u Ng\d{o}c}

\begin{document}
\allowdisplaybreaks

\newcommand{\arXivNumber}{1706.09935}

\renewcommand{\PaperNumber}{084}

\FirstPageHeading

\ShortArticleName{Faithful Semitoric Systems}

\ArticleName{Faithful Semitoric Systems}

\Author{Sonja HOHLOCH~$^{\dag^1}$, Silvia SABATINI~$^{\dag^2}$, Daniele SEPE~$^{\dag^3}$ and Margaret SYMINGTON~$^{\dag^4}$}

\AuthorNameForHeading{S.~Hohloch, S.~Sabatini, D.~Sepe and M.~Symington}

\Address{$^{\dag^1}$~Department of Mathematics - Computer Science, University of Antwerpen, \\
\hphantom{$^{\dag^1}$}~Campus Middelheim, Building G, M.G.211, Middelheimlaan 1, 2020 Antwerpen, Belgium}
\EmailDD{\href{mailto:sonja.hohloch@uantwerpen.be}{sonja.hohloch@uantwerpen.be}}
\URLaddressDD{\url{https://www.uantwerpen.be/en/staff/sonja-hohloch/}}

\Address{$^{\dag^2}$~Mathematisches Institut, Universit\"at zu K\"oln, Weyertal 86-90, D-50931 K\"oln, Germany}
\EmailDD{\href{mailto:sabatini@math.uni-koeln.de}{sabatini@math.uni-koeln.de}}
\URLaddressDD{\url{http://silvias79.wixsite.com/silvia-sabatini-math}}

\Address{$^{\dag^3}$~Universidade Federal Fluminense, Instituto de Matem\'atica, \\
\hphantom{$^{\dag^3}$}~Departamento de Matem\'atica Aplicada, Rua Professor Marcos Waldemar de Freitas Reis,\\
\hphantom{$^{\dag^3}$}~s/n, Bloco H, Campus do Gragoat\'a, CEP 24210-201, Niter\'oi, RJ, Brazil}
\EmailDD{\href{mailto:danielesepe@id.uff.br}{danielesepe@id.uff.br}}
\URLaddressDD{\url{https://sites.google.com/site/danielesepe/}}

\Address{$^{\dag^4}$~Department of Mathematics, Mercer University,\\
\hphantom{$^{\dag^4}$}~1501 Mercer University Drive, Macon, GA 31207, USA}
\EmailDD{\href{mailto:symington_mf@mercer.edu}{symington\_mf@mercer.edu}}
\URLaddressDD{\url{http://faculty.mercer.edu/symington_mf/}}

\ArticleDates{Received July 07, 2017, in final form July 30, 2018; Published online August 16, 2018}

\Abstract{This paper consists of two parts. The first provides a review of the basic properties of integrable and almost-toric systems, with a particular emphasis on the integral affine structure associated to an integrable system. The second part introduces faithful semitoric systems, a~generalization of semitoric systems (introduced by V\~u Ng\d{o}c and classified by Pelayo and V\~u Ng\d{o}c) that provides the language to develop surgeries on almost-toric systems in dimension~4. We prove that faithful semitoric systems are natural building blocks of almost-toric systems. Moreover, we show that they enjoy many of the properties that their (proper) semitoric counterparts do.}

\Keywords{completely integrable Hamiltonian systems; almost toric systems; semitoric systems; integral affine geometry; focus-focus singularities}

\Classification{37J35; 37J05; 53D20; 70H06}

\tableofcontents

\section{Introduction}\label{sec:introduction-1}
A driving problem in Hamiltonian mechanics and symplectic geometry is to classify integrable systems up to a suitable notion of equivalence. An {\em integrable system} is a triple $\is$, where $\sm$ is a $2n$-dimensional symplectic manifold and $\Phi\colon \sm \to \R^n$ is a smooth map whose components are in Poisson involution and functionally independent almost everywhere on~$M$. This paper introduces faithful semitoric systems, a category of integrable systems on $4$-dimensional symplectic manifolds that generalizes toric and semitoric systems and lays the foundation for studying almost-toric systems. A key feature of faithful semitoric systems is that they behave well under a process of taking appropriate subsystems, a fact that facilitates development of precise language to define, for faithful semitoric systems, integrable surgeries in the sense of Zung~\cite{zung_ii}.

In general, classification of integrable systems becomes a tractable problem only under assumptions that restrict the topology of fibers of the system. Intuitively, the greatest challenge comes from non-compactness of the group action arising from the flow of the Hamiltonian vector fields of the components of the moment map. Accordingly, full classifications were first established for toric systems in which the $\R^n$-action descends to a $\mathbb{T}^n$-action. Building upon the foundational results of Atiyah~\cite{atiyah} and Guillemin and Sternberg~\cite{guillemin-sternberg}, Delzant~\cite{delzant} classified toric systems on closed manifolds up to isomorphism. More recently, Karshon and Lerman~\cite{karshon-lerman} have extended Delzant's classification to non-compact toric manifolds, relying upon the local normal forms of Guillemin and Sternberg~\cite{GS2} and Marle~\cite{marle}.

Once one allows non-compactness of the group acting on the total space, complexity of both the fibers and of the total space can be reasonably controlled by restricting the singularities of the moment map. Symington \cite{symington} and \vungoc~\cite{vu-ngoc} have proposed a notion of almost-toric systems on 4-dimensional symplectic manifolds that includes toric systems but also allows for so-called focus-focus fibers, which can be thought of as the Lagrangian analog of the nodal fibers that arise in Lefschetz fibrations. The diffeomorphism types of closed manifolds that support an almost-toric system has been determined (cf.\ Leung and Symington \cite{leung_symington}), and recently, almost-toric systems have proved to be of independent interest in symplectic topology (cf.\ Vianna~\cite{vianna,vianna_inf,vianna_dp}).

While the classification problem of almost-toric systems has not been settled, even in the compact case, an important subclass of almost-toric systems has been completely understood: Pelayo and \vungoc~\cite{pelayo-vu-ngoc-inventiones,pelayo-vu-ngoc-constr} have classified semitoric systems, which were initially introduced by \vungoc\ in~\cite{vu-ngoc}. An integrable system $\vat$ is semitoric if it is almost-toric and if $J$ is a proper moment map of an effective Hamiltonian $S^1$-action. Semitoric systems, whose total spaces may be non-compact, share many fundamental properties with closed symplectic toric manifolds, like connectedness of the fibers of the moment map, but their classification is significantly more involved as the presence of focus-focus fibers introduces more data (see Pelayo and \vungoc~\cite{pelayo-vu-ngoc-inventiones}). While semitoric systems appear naturally both in symplectic topology and in (quantum) Hamiltonian mechanics (cf.\ Eliashberg, Polterovich, Le Floch, Pelayo and \vungoc~\cite{el_pol,lpv,pvn-spin}), the properness condition on $J$ excludes some familiar almost-toric integrable systems, such as the spherical pendulum (cf.\ Duistermaat~\cite{duistermaat}). For this reason, Pelayo, Ratiu and \vungoc~\cite{pvr,pvr_carto} introduce a family of almost-toric systems that share some of the main properties of semitoric systems, like connectedness of the fibers of the moment map, while allowing enough freedom to include examples such as the spherical pendulum. These systems are called proper semitoric. In such an almost-toric system~$\vat$ the moment map~$\Phi$ is proper and~$J$ is the moment map of an effective Hamiltonian $S^1$-action that satisfies some constraints on the sets of singular points and values (see Definition~1.3 in Pelayo, Ratiu and \vungoc~\cite{pvr_carto}).

Faithful semitoric systems, defined and studied in this paper, can be viewed as an extension of proper semitoric systems to the non-compact setting but were defined with different purposes in mind. For instance, faithful semitoric systems provide the appropriate setting to accommodate
non-compact systems that are convenient local models and building blocks for almost-toric systems. The essential difference between faithful semitoric systems and proper semitoric systems is that the moment maps of the former are merely required to be proper onto their image while the moment map of a proper semitoric system must be proper.

The definition of a faithful semitoric system is crafted so that appropriately chosen subsystems are again faithful semitoric. Specifically, given an open, connected subset $U$ of the moment map image of a faithful semitoric system, if the intersection of $U$ with any vertical line is
either empty or connected then restricting the moment map to the preimage of $U$ yields a faithful semitoric system. In fact, such a subsystem of a (proper) semitoric system is faithful semitoric. In a forthcoming paper, \cite{HSSS-surgeries}, the process of taking such subsystems is an essential ingredient in the definition of surgeries in the category of faithful semitoric systems. We plan to apply those surgeries to the determination of which Hamiltonian $S^1$-spaces underlie compact semitoric systems (\cite{HSSS-lifting}, forthcoming), thus completing the work started in Hohloch, Sabatini and Sepe~\cite{HSS}. The language of faithful semitoric systems also allows one to have a more conceptual understanding of the local-to-global arguments in the classification of semitoric systems (cf.\ Pelayo and \vungoc~\cite{pelayo-vu-ngoc-inventiones,pelayo-vu-ngoc-constr}); this is also going to be explored in a separate paper.

The main results of this paper are as follows:
\begin{enumerate}[label=(\Alph*), ref = (\Alph*), leftmargin=*]\itemsep=0pt
\item \label{item:30} The fibers of a faithful semitoric system are connected (see Theorem \ref{thm:connected}).
\item \label{item:1} A connected component of a fiber of an almost-toric system admits an open neighborhood that is isomorphic to a faithful semitoric system (see Proposition~\ref{prop:at-vat}).
\item \label{item:2} Using terminology analogous to that for (proper) semitoric systems (cf.\ Pelayo, Ratiu and \vungoc~\cite{pvr_carto}), faithful semitoric systems possess {\em cartographic homeomorphisms} (see Theorem~\ref{prop:rh}). These are homeomorphisms of the moment map image onto subsets of $\R^2$ that encode the induced {\em $\Z$-affine structures} (see Section~\ref{sec:z-affine-structure}) and generalize inclusion of the moment map image for toric systems into the ambient Euclidean space. In particular, the monodromy introduced by focus-focus fibers is encoded via {\em vertical cuts}. (It is worth mentioning that, for semitoric systems, the images of cartographic homeomorphisms are {\em convex}, possibly noncompact, polygons, as shown in \vungoc~\cite{vu-ngoc}).
\item \label{item:5} The space of all cartographic homeomorphisms of a given faithful semitoric system is described (see Theorem~\ref{thm:different_carto}), generalizing the analogous result for a semitoric system (cf.\ \vungoc~\cite{vu-ngoc}). This description can be used to construct an invariant of the isomorphism class of a faithful semitoric system analogous to the {\em semitoric polygon} of Pelayo and \vungoc~\cite[Definition~4.5]{pelayo-vu-ngoc-inventiones} (see Lemmas~\ref{lemma:im_invariant} and~\ref{lemma:im_inv_full}).
\item \label{item:6} Given a faithful semitoric system $\vat$ and a~cartographic homeomorphism $f\colon \Phi(M) \to \R^2$, the composition $f\circ\Phi$ may lack the smoothness required of a moment map. We provide a method for smoothing $f\circ\Phi$ to obtain a faithful semitoric system isomorphic to $\vat$ whose moment map image equals $f (\Phi(M) )$ on the complement of arbitrarily small neighborhoods of the cuts used to define~$f$ (see Theorem~\ref{thm:rect-embedding}).
\end{enumerate}

Result \ref{item:30} is a simple consequence of the work by Pelayo, Ratiu and \vungoc\ \cite{pvr,pvr_carto} and is included to highlight that connectedness of the fibers follows under the slightly weaker assumption that the moment map be proper onto its image.

Result \ref{item:1} establishes faithful semitoric systems as building blocks for almost-toric systems. While it is probably known to experts in the area, we could not find a complete, self-contained proof of this fact and decided to include it, along with proofs of basic topological facts leading up to it (see Section~\ref{sec:at_generalities}).

Results \ref{item:2} and \ref{item:5} are not surprising in light of the work in Pelayo, Ratiu and \vungoc\ \cite{pvr_carto,vu-ngoc}. However, because of the difference in purpose between those papers and this one, even the fundamental notion of isomorphism of systems differs. We use a notion of isomorphism that applies more generally to any integrable system (Definition~\ref{defn:cihs}), but we also prove that the presence of focus-focus points imposes restrictions on isomorphisms between faithful semitoric systems that causes the notions to align in the (proper)
semitoric context (see Proposition~\ref{prop:im_equ}). Furthermore, we provide explicit proofs of Results~\ref{item:2} and~\ref{item:5} for the following reasons:

\begin{itemize}[leftmargin=*]\itemsep=0pt
\item Our alternative proof of the existence of cartographic homeomorphisms in the case in which the defining cuts disconnect the moment map image allows us to avoid the `homotopy argument' of Pelayo, Ratiu and \vungoc\ \cite[Step~5 of the proof of Theorem~B]{pvr_carto}.
\item The description of the set of cartographic homeomorphisms of a~faithful semitoric system is analogous to that of a semitoric system. However, the potential for infinitely many focus-focus points in a faithful semitoric system gives a richer behavior, as can be seen by comparing Section~\ref{sec:set-cart-home} with \vungoc\ \cite[Section~4]{vu-ngoc}.
\end{itemize}

Result \ref{item:6} ensures the existence of $\boldsymbol{\eta}$-cartographic systems in the isomorphism class of a~faithful semitoric system -- with respect to any of the notions of isomorphism of Definition~\ref{defn:im_faithful_st}. Result~\ref{item:6} is useful for applications because
one of the cartographic homeomorphisms of an $\boldsymbol{\eta}$-cartographic system is the identity on the complement of a small neighborhood of the vertical cuts, and hence the moment map of the $\boldsymbol{\eta}$-cartographic system can be thought of as a~\smash{`smoothing'} of that cartographic homeomorphism. As one application, $\boldsymbol{\eta}$-cartographic systems play an important role in defining surgeries of faithful semitoric systems (\cite{HSSS-surgeries}, forthcoming). Also, Result~\ref{item:6} allows one to make precise the notion that the image of a~cartographic homeomorphism is a limit of moment map images (see Proposition~\ref{prop:limit}).

{\bf Structure of the paper.}
This paper is split in two parts. Part \ref{sec:integrable_systems}, consisting of Sections~\ref{sec:integr-syst} and~\ref{sec:at}, defines and explains notions that we use throughout the paper, while also establishing notation.

While Section \ref{sec:integr-syst} should serve as a self-contained primer to guide readers unfamiliar with the subject through this paper (and the forthcoming~\cite{HSSS-surgeries}), it may be of interest to experts in the field as well, for a few ideas which do not appear in many other places. For instance, we introduce the notion of a {\em faithful} integrable system, which is one whose moment map image is homeomorphic, as a subset of~$\R^n$, to the leaf space of the system (see Section~\ref{sec:faithful-moment-maps}). Moreover, in Section~\ref{sec:faithful}, we elaborate on the notion of {\em cartographic homeomorphisms} that is introduced in Pelayo, Ratiu and \vungoc~\cite{pvr_carto} for proper semitoric systems and establish some general properties for these objects. Note that the section should not be taken as an exhaustive reference for either topological or symplectic aspects of integrable systems.

In Section \ref{sec:at}, almost-toric systems are defined and their basic properties are explored. In particular, the neighborhood of a connected component of a fiber is described (see Section~\ref{sec:at_generalities}) and, in preparation for the next section on faithful semitoric systems, we describe properties of systems that are both faithful and almost-toric (Section~\ref{sec:faithful_at}).

Part \ref{part:faithful-semitoric}, consisting of Sections~\ref{sec:defin-basic-prop} and~\ref{sec:dscc-home-boldsymb}, contains the definition of faithful semitoric systems as well as all the main results described above. While Part~\ref{part:faithful-semitoric} uses notions and ideas that appear in Part~\ref{sec:integrable_systems}, it is sufficiently self-contained that readers who are familiar with the basic properties of integrable and almost-toric systems can skip Part \ref{sec:integrable_systems} and refer to it as they read Part~\ref{part:faithful-semitoric}. Section~\ref{sec:defin-basic-prop} contains main results~\ref{item:30} and~\ref{item:1} as Theorem~\ref{thm:connected} and Proposition~\ref{prop:at-vat}, characterizations of the moment map image of faithful semitoric systems (Corollary~\ref{cor:contractible_mom_map}), and a useful criterion to determine which saturated subsystems of faithful semitoric systems are
faithful semitoric (Proposition~\ref{prop:sub_mostly_vat}). The relation between faithful semitoric systems and the (proper) semitoric systems of \vungoc~\cite{vu-ngoc} and Pelayo, Ratiu and \vungoc~\cite{pvr_carto} is studied in Section~\ref{sec:dscex-relat-other}, where examples are provided of proper semitoric systems that are not semitoric in the sense of \vungoc~\cite{vu-ngoc}, and of faithful semitoric systems that are not proper semitoric (see Example~\ref{exm:relation}). Section~\ref{sec:dscp-noti-isom} introduces various notions of isomorphism for faithful semitoric systems, and relates the notion of isomorphism in Pelayo, Ratiu and \vungoc~\cite[Definition~1.5]{pvr_carto} with that of Definition~\ref{defn:cihs} (see Remark~\ref{rmk:toric_im} and Proposition~\ref{prop:im_equ}). Section~\ref{sec:geom-impl-s1} explores
the consequences of the presence of the $S^1$-action.

Section~\ref{sec:dscc-home-boldsymb} contains all results concerning faithful semitoric systems and cartographic homeomorphisms. The existence of
cartographic homeomorphisms (Theorem~\ref{prop:rh}) is proved in Section~\ref{sec:cuts-cart-home}, which also establishes some useful topological properties of the complements of the cuts needed to define cartographic homeomorphisms. Section~\ref{sec:set-cart-home} describes the set of cartographic homeomorphisms associated to a given faithful semitoric system, paying particular attention to the subtleties that arise from allowing infinitely many focus-focus points (see Theorem~\ref{thm:different_carto}). By understanding the set of cartographic homeomorphisms associated to a given faithful semitoric system, we classify faithful semitoric systems with no focus-focus points up to any of the notions of isomorphisms of Definition~\ref{defn:im_faithful_st} (see Lemma~\ref{lemma:class_no_ff}), and construct an invariant of faithful semitoric systems with at least one focus-focus point up to isomorphisms of integrable systems (see Lemma~\ref{lemma:im_inv_full}). The latter can be viewed as a first step towards achieving a classification of faithful semitoric systems as integrable systems. Finally, Section~\ref{sec:choos-an-appr} proves that, in some sense, cartographic homeomorphisms can be made smooth everywhere by modifying them on arbitrarily small neighborhoods of the defining cuts. This is the content of Theorem~\ref{thm:rect-embedding}, which can be used to establish the existence of $\boldsymbol{\eta}$-cartographic faithful semitoric systems in any given isomorphism class (Theorem~\ref{thm:eta-carto}).

\subsection*{Notation and conventions}

\subsubsection*{Topological conventions}

\begin{itemize} [leftmargin=*]\itemsep=0pt
\item A subset of a topological space is endowed with the subspace topology unless otherwise stated.
\item A {\em pair} of topological spaces $(Y,Z)$ consists of a topological space $Y$ together with a~subset $Z \subset Y$ endowed with the relative topology. A~topological embedding of pairs of topological spaces $(Y_1,Z_1)$, $(Y_2,Z_2)$ is a topological embedding $\chi\colon Y_1 \to Y_2$ that restricts to a topological embedding of $Z_1$ into $Z_2$. A~homeomorphism between pairs of topological spaces is a~topological embedding onto the target.
\item A map $f\colon X \to Y$ between topological spaces $X$ and $Y$ is {\em proper} if for any compact subset $K \subset Y$, the preimage $f^{-1}(K) \subset X$ is compact; it is {\em proper onto its image} if the map $f\colon X \to f(X)$ is proper.
\end{itemize}

\subsubsection*{Smoothness conventions}
Throughout the paper, we use the following standard conventions on smoothness.
\begin{itemize}[leftmargin=*]\itemsep=0pt
\item Following Joyce \cite{joyce}, consider the subspace $[0,+\infty[\,^n\subset \R^n$ and let $M$ be a topological space. An {\em $n$-dimensional smooth atlas with corners} on $M$ is a set $\mathcal{A} := \left\{\left(U_i,\chi_i\right)\right\}$, where
 \begin{itemize}[leftmargin=*]\itemsep=0pt
\item the set $\left\{U_i\right\}$ is an open cover of $M$;
\item for each $i$, there is an open set $V_i\subset [0,+\infty[\,^n$ such that the map $\chi_i \colon U_i \to V_i$ is a homeomorphism; and
\item for all $i$, $j$ with $U_i \cap U_j \neq \varnothing$, the map $\chi_j \circ \chi_i^{-1} \colon \chi_i(U_i\cap U_j) \to \chi_j(U_i\cap U_j)$ is a diffeomorphism.
\end{itemize}
 A {\em smooth manifold with corners} of dimension $n$ is a Hausdorff, second countable topological space together with an $n$-dimensional smooth atlas with corners. A~{\em smooth structure with corners} on a topological space is an equivalence class of smooth atlases with corners, where two atlases are deemed equivalent if their union is again a smooth atlas with corners of a given dimension.

 \item A smooth atlas with corners is {\em $\Z$-affine} (or {\em integral affine}) if the transition maps $\chi_j \circ \chi_i^{-1} \colon$ $\chi_i(U_i\cap U_j) \to \chi_j(U_i\cap U_j)$ of the atlas are of the form
 \begin{gather*}
 x \mapsto Ax + b,
 \end{gather*}
for some $(A,b) \in \mathrm{AGL}(n;\Z)=\mathrm{GL}(n;\Z) \ltimes \R^n$, where $n$ is the dimension of the atlas. A~{\em $\Z$-affine manifold with corners} of dimension $n$ is a Hausdorff, second countable topological space together with an $n$-dimensional $\Z$-affine atlas with corners. A {\em $\Z$-affine structure with corners} on a topological space is an equivalence class of $\Z$-affine atlases with corners, where two atlases are deemed equivalent if their union is again a $\Z$-affine atlas with corners of a given dimension.

{\em Smooth atlases} and {\em $\Z$-affine atlases}, without corners, (and the corresponding manifold structures) are defined as above, with the stipulation that the images of the coordinate charts are subsets of~$\R^n$. Also, note that, because the transition maps of are smooth, a $\Z$-affine atlas (with or without corners) is also a smooth atlas, and hence defines a unique smooth structure.

 \item Let $A \subset \R^n$ be a subset. A map $f\colon A \to \R^m$ is said to be {\em smooth} if for all $x \in A$ there exists an open neighborhood $U_x \subset \R^n$ of $x$ and a smooth map $f_x \colon U_x
 \to \R^m$ that is a local extension of~$f$.
\item A map $f \colon A \subset \R^n \to \R^m$ is a {\em smooth embedding} if it is a diffeomorphism onto its image.
\item Manifolds are assumed to be without boundary or corners unless otherwise stated.
\end{itemize}

\subsubsection*{Boundary conventions}\label{sec:bdy_conditions}
Two types of boundaries of subsets $X\subset\R^n$ are dealt with in this paper whenever $X$ is a smooth manifold with corners embedded in~$\R^n$. The topological boundary, the closure of $X\subset\R^n$ minus its interior, is denoted $\mathrm{Bdy}(X)$. Meanwhile, its boundary as a manifold with corners, $X\cap \mathrm{Bdy}(X)$, is denoted $\partial X$. For instance if $X=\{(x,y)\,|\, \abs{x}<1 \ {\rm and}\ \abs{y}\le 1\}$, then
\begin{gather*}
 \mathrm{Bdy}(X) = \{ (x,y)\,|\, \abs{x}=1, \abs{y}\le 1 \ {\rm or}\ \abs{y}=1, \abs{x}\le 1\}
\end{gather*}
and
\begin{gather*}
 \partial X = \{(x,y)\,|\, \abs{y}=1, \abs{x}< 1\}.
\end{gather*}

\subsubsection*{Group conventions}
Throughout the paper, the identification $S^1 \cong \R/2\pi\Z$
is used tacitly.

\part{Primer on integrable and almost-toric systems}\label{sec:integrable_systems}
This part introduces the basic notions regarding integrable and
almost-toric systems
that are used throughout the paper. Section \ref{sec:integr-syst}
 introduces integrable systems, their subsystems and some of their
 invariants up to isomorphisms. While most of these notions are standard and
 appear in more comprehensive texts on the topology and geometry of
 integrable systems, such as \cite{bol+fom,vu_ngoc_book}, the notions of faithfulness of
 a system and the definition of a cartographic homeomorphism,
 introduced in Section \ref{sec:integr-syst}, seem to be new. Section
 \ref{sec:at} defines almost-toric systems in dimension 4 and proves
 some of their fundamental properties. Almost-toric systems
 generalize toric systems by allowing the presence of
 so-called {\em focus-focus leaves} (see Section~\ref{sec:at_generalities}).

\section{Completely integrable Hamiltonian systems}\label{sec:integr-syst}
This section presents the category of integrable systems and defines a coarse topological invariant: the {\em leaf space} of an integrable system (see Definition \ref{defn:leaf_space}). Systems whose moment map images can be identified with their leaf spaces play an important role in this paper and are studied in Section~\ref{sec:faithful-moment-maps}; we call such systems {\em faithful}. In Sections~\ref{sec:z-affine-structure}--\ref{sec:toric-leaf}, we endow large subsets of the leaf space of an integrable system with an {\em $\Z$-affine structure}. First, following Duistermaat~\cite{duistermaat}, we show how the part of the leaf space corresponding to regular leaves inherits such a structure in Section~\ref{sec:z-affine-structure}. Second, we identify a class of systems that are isomorphic to systems equipped with Hamiltonian torus actions of maximal dimension: these are called {\em weakly toric}, are related to symplectic toric manifolds and are studied in Section~\ref{sec:toric-systems}. (This notion is a mild generalization of systems of {\em toric type} introduced by \vungoc\ in \cite[Definition~2.1]{vu-ngoc}, see Remark~\ref{rmk:toric_type}.) Third, Section~\ref{sec:toric-leaf} extends the $\Z$-affine structure on the regular part of the leaf space to include singular leaves that admit a neighborhood supporting a Hamiltonian torus action of maximal dimension. In Section \ref{sec:developing-maps}, we recall a~fundamental property that $\Z$-affine structures enjoy, namely that they can be {\em developed}. Finally, following Pelayo, Ratiu and \vungoc\ \cite{pvr_carto}, we introduce the notion of {\em cartographic homeomorphisms}, which, intuitively, can be thought of a way to encode the above $\Z$-affine structure in a way that is compatible with singular orbits of the system. Throughout this section there is no restriction on the dimension of the phase space
of an integrable system.

\subsection{Definition, subsystems and leaf spaces}\label{sec:dscd-subsyst-leaf}
Let $\sm$ be a symplectic manifold. Given a smooth function $H \in C^{\infty}(M)$, its {\em Hamiltonian vector field} $X_H \in \mathfrak{X}(M)$ is the unique vector field defined implicitly by the equation $\omega(X_H,\cdot) = dH$. The symplectic form $\omega$ induces a {\em Poisson structure} on $M$, i.e., a Lie bracket $\{\cdot,\cdot\}$ on $C^{\infty}(M)$ that satisfies the Leibniz identity, defined by $ \{H_1,H_2 \}:= \omega (X_{H_1},X_{H_2} )$, for $H_1,H_2\in C^{\infty}(M)$. These notions allow to introduce the category of integrable systems.

\begin{Definition}\label{defn:cihs}
For any $n \geq 1$, the {\em category of completely integrable Hamiltonian systems with $n$ degrees of freedom}, denoted by $\mathcal{IS}(n)$, has objects and morphisms as follows:
 \begin{itemize}[leftmargin=*]\itemsep=0pt
 \item {\bf Objects}: {\em completely integrable Hamiltonian systems} $(M, \om, \Phi)$ where $(M, \om)$ is a~$2n$-dimen\-sio\-nal symplectic manifold and
 \begin{gather*}
 \Phi:=(H_1,\ldots,H_n) \colon \ \sm \to \R^n
 \end{gather*}
 a smooth map satisfying
 \begin{itemize}[leftmargin=*]\itemsep=0pt
 \item $\{H_i,H_j\} =0$ for all $i,j =1,\ldots,n$, where $\{\cdot,\cdot\}$ is the Poisson bracket induced by $\omega$;
 \item $\Phi$ is a submersion on an open dense subset, i.e., there exists an open, dense subset $V \subset M$ such that, for all $p \in V$, $d_pH_1 ,\ldots, d_pH_n$ are linearly independent.
 \end{itemize}
Sometimes, for brevity, $\Phi$ is referred to as a (completely) integrable (Hamiltonian) system. Its component $H_i$ is called the $i^{\rm th}$ {\em integral $($of motion$)$}.
 \item {\bf Morphisms}: {\em isomorphisms of integrable systems} $(\Psi,\psi)$, where, for $i=1,2$, $(M_i,\omega_i,\Phi_i)$ is a~completely integrable Hamiltonian system, $\Psi \colon (M_1,\omega_1) \to (M_2,\omega_2)$ is a symplectomorphism, $\psi \colon \Phi_1(M_1) \to \Phi_2(M_2)$ is a diffeomorphism, and the following diagram commutes:
 \begin{gather*}
 \xymatrix{ (M_1,\omega_1) \ar[r]^-{\Psi} \ar[d]_-{\Phi_1} &
 (M_2,\omega_2) \ar[d]^-{\Phi_2} \\
		 \Phi_1(M_1) \ar[r]_-{\psi} & \Phi_2(M_2).}
 \end{gather*}
 \end{itemize}
\end{Definition}

Given $\Phi=(H_1,\ldots,H_n) \colon \sm \to \R^n$, for $1 \leq i \leq n$, the Hamiltonian vector field associated to $H_i$ is denoted by $X_i$. If the flows of the vector fields $X_1,\ldots,X_n$ are complete, then there is a Hamiltonian $\R^n$-action on $\sm$, one of whose moment maps is precisely $\Phi$.

Throughout this paper, unless otherwise stated,
\begin{center}
 {\it integrable systems have compact fibers},
\end{center}
so the above completeness assumption is satisfied and $\Phi$ is a moment map for a Hamiltonian $\R^n$-action (upon identifying the Lie algebra of $\R^n$ with $\R^n$). For this reason, it is referred to as the {\em moment map} of the system.

Many of the integrable systems considered in this paper arise from restricting a given system to a subset.

\begin{Definition}\label{defn:subsystem}
A {\em subsystem} of an integrable system $\is$ is an integrable system \linebreak $ (V,\omega|_{V}, \Phi|_{V} )$ where $V$ is an open subset of~$M$. If $V=\Phi^{-1}(U)$ for some open subset~$U$ of~$\Phi(M)$, the subsystem $(V,\omega|_{V}, \Phi|_{V})$ is also referred to as the {\em subsystem of~$\is$ relative to~$U$}.
 \end{Definition}
\begin{rmk}\label{rmk:subsystems} Subsystems of integrable systems with compact fibers need not have compact fibers. If $\is$ is an integrable system with compact fibers, for any point $p \in M$, the subsystem $ (M\smallsetminus \{p \}, \omega_{M\smallsetminus \{p\}}, \Phi|_{M\smallsetminus \{p\}} )$ does not have compact fibers. If $p$ is chosen to be regular in the sense of Definition~\ref{defn:singular_points_fibers}, then the above example shows that subsystems of integrable systems supporting a Hamiltonian $\R^n$-action need not support a Hamiltonian $\R^n$-action.
\end{rmk}

The most accessible feature of an integrable system is the image of~$\Phi$. However, it is the leaf space of the system that reliably reflects some of the topological structure.

\begin{Definition}\label{defn:leaf_space}
 Given an integrable system $\is$,
 \begin{itemize}[leftmargin=*]\itemsep=0pt
 \item a {\em leaf} is a connected component of a fiber of $\Phi$;
 \item its {\em leaf space} is the topological space $\mathcal{L} : = M/{\sim}$, where $p \sim p'$ if $p$ and $p'$ belong to the same leaf, endowed with the quotient topology;
 \item a subsystem $ (V,\omega|_{V}, \Phi|_{V} )$ of $\is$ is {\em saturated} (with respect to the quotient map \smash{$q \colon \!M {\to} \mathcal{L}$}) if any leaf of $\is$ that intersects $V$ is contained in $V$.
 \end{itemize}
\end{Definition}

\begin{figure}[h] \centering
 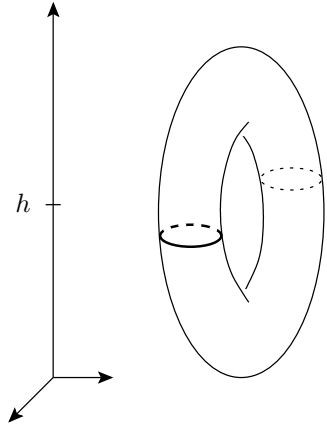
\caption{Level set of $h$ consists of two circles, but any sufficiently small open neighborhood of the heavily marked leaf contains a subsystem that is saturated with respect to the quotient map $q \colon \mathbb{T}^2 \to
 \mathcal{L}$ but that fails to be the preimage of an open set in the image of the height function.} \label{LABEL1}
\end{figure}

\begin{rmk}\label{rmk:saturated_vs_relative} Given an integrable system $\is$, for any open subset $U \subset \Phi(M)$, the subsystem of $\is$ relative to $U$ is saturated in the sense of Definition~\ref{defn:leaf_space}. However, a~saturated subsystem of $\is$ need not be a~subsystem relative to an open subset $U \subset \Phi(M)$, as illustrated by the following example. Indeed, consider the standard height function on~$\mathbb{T}^2$ embedded submanifold in~$\R^3$ as in Fig.~\ref{LABEL1}. Endowing $\mathbb{T}^2$ with an area form, the height function gives rise to an integrable system. Observe that any sufficiently small open neighborhood of the leaf marked heavily in Fig.~\ref{LABEL1} contains an open neighborhood $V$ that is saturated with respect to the quotient map $q \colon \mathbb{T}^2 \to \mathcal{L}$ but that fails to be the preimage of an open set in the image of the height function.
\end{rmk}

Given an integrable system $\is$ with leaf space $\mathcal{L}$, the moment map $\Phi$ factors through~$\mathcal{L}$, inducing a continuous map $\pi\colon \mathcal{L} \to \Phi(M)$ that makes the following diagram commute:
\begin{gather*}
 \xymatrix{& \sm \ar[dl]_-{q} \ar[dr]^-{\Phi} & \\
 \mathcal{L} \ar[rr]^-{\pi} & & B: =\Phi(M) \subset \R^n.}
\end{gather*}

Leaf spaces descend to isomorphism classes of integrable systems
and behave well with respect to saturated subsystems:

 \begin{Lemma}\label{lemma:leaf_space_invariant}\quad
 \mbox{}
 \begin{enumerate}[label={\rm (\arabic*)}, ref=(\arabic*),leftmargin=*]\itemsep=0pt
 \item \label{item:11} Isomorphic integrable systems possess homeomorphic leaf spaces.
 \item \label{item:12} The leaf space of a saturated subsystem of $\is$ embeds naturally in the leaf space of~$\is$.
 \end{enumerate}
 \end{Lemma}
 \begin{proof} To prove \ref{item:11}, suppose that integrable systems $(M_1,\omega_1,\Phi_1)$ and $(M_2,\omega_2,\Phi_2)$ are isomorphic via the pair $(\Psi,\psi)$. For $i=1,2$, denote the leaf space of
 $(M_i,\omega_i,\Phi_i)$ by $\mathcal{L}_i$ and let $q_i\colon M_i \to \mathcal{L}_i$ denote the quotient map. Let $L_1 \subset M_1$ be a leaf of $\Phi_1$ in the sense of Definition~\ref{defn:leaf_space}. Because
 $\Phi_1(L_1)$ is constant and $\Phi_2 \circ \Psi = \psi \circ \Phi_1$, the continuity of $\Psi$ ensures that $\Psi(L_1) \subset M_2$ is contained in a leaf of $\Phi_2$. Therefore, there exists a map $\zeta\colon \mathcal{L}_1 \to \mathcal{L}_2$ such that $\zeta \circ q_1 = q_2 \circ \Psi$. (Observe that $\pi_2 \circ \zeta = \psi \circ \pi_1$, where, for $i=1,2$, $\pi_i \colon \mathcal{L}_i \to \Phi_i (M_i )$ is the induced
 continuous map, since $ (\Psi,\psi )$ is an isomorphism of integrable systems.) Continuity of $q_2 \circ \Psi$ and the equa\-li\-ty $\zeta \circ q_1 = q_2 \circ \Psi$ imply that $\zeta$ is continuous, since $q_1$ is a quotient map. Applying the above argument with $\big(\Psi^{-1},\psi^{-1}\big)$ yields the existence of a~continuous inverse $\zeta^{-1} \colon \mathcal{L}_2 \to \mathcal{L}_1$ of $\zeta$, thus showing that $\zeta$ is a homeomorphism as desired.

To prove \ref{item:12}, let $(V,\omega|_V,\Phi|_V)$ be a saturated subsystem of $\is$, and denote the leaf spaces of $(V,\omega|_V,\Phi|_V)$ and of $\is$ by $\mathcal{L}_V$ and $\mathcal{L}$ respectively. The above argument shows that there is a continuous map $\iota\colon \mathcal{L}_V \to \mathcal{L}$. (This fact does not need the subsystem to be saturated.) The fact that the subsystem is saturated implies that $\iota$ is injective and that its image, which is equal to $q(V)$, is an open subset of $\mathcal{L}$. To see that $\iota$ is a~topological embedding, observe that $V = q^{-1}(q(V))$ as $V$ is saturated with respect to $q$. Therefore, the fact that the inclusion $V \hookrightarrow M$ is a topological embedding implies that $\iota\colon \mathcal{L}_V \to \mathcal{L}$ is a~topological embedding as desired.
 \end{proof}

\begin{rmk}\label{rmk:no_embed}The condition on the subsystem being saturated in part~\ref{item:12} of Lemma~\ref{lemma:leaf_space_invariant} is necessary. For instance, consider the integrable system $H\colon \big(\R^2,\mathrm{d} x \wedge \mathrm{d}y\big) \to \R$ given by \smash{$H(x,y) = x$}. The leaf space of the subsystem $\big(\R^2 \smallsetminus \{(0,0)\}, \mathrm{d} x \wedge \mathrm{d} y, H|_{\R^2 \smallsetminus\{(0,0)\}}\big)$ is not Hausdorff, while that of $\big(\R^2,\mathrm{d} x \wedge \mathrm{d}y, H\big)$ is.
\end{rmk}

The following result illustrates the fact that some topological properties of the phase space of an integrable system are reflected in the topology of the associated leaf space.
\begin{Lemma}\label{lemma:leaf_conn} The set of connected components of the leaf space of an integrable system and the set of connected components of its total space are in bijection.
\end{Lemma}
\begin{proof} Let $\is$ be an integrable system, let $\mathcal{L}$ denote its leaf space and denote the quotient map by $q \colon M \to \mathcal{L}$. Since~$q$ is continuous and surjective, the induced map from the set of connected components of $M$ to that of connected components of~$\mathcal{L}$ is surjective. To establish that the induced map is injective, we show that the image under $q$ of a~disconnected subset is disconnected. First note that, by virtue of being a~manifold, $M$ is locally connected and hence each connected component is open. Consider a disconnected subset $M_1 \sqcup M_2 \subset M$ and suppose its image $q(M_1 \sqcup M_2)$ is connected. Consider $p \in q(M_1) \cap q(M_2)$. Then $M_1\cap q^{-1}(p)$ and $M_2 \cap q^{-1}(p)$ provide a disconnection of $q^{-1}(p)$ because $M_1$ and $M_2$ are open and disjoint. But this is impossible because the fibers
 of $q$ are by definition connected. The injectivity of the induced map on connected components follows immediately because the preimage under $q$ of a connected component of $L$ must be connected. Thus, the map that~$q$ induces on connected components is bijective as desired.
\end{proof}

Topological properties of the leaves of an integrable system determine some topological properties of the associated leaf space, as illustrated by the following result.

 \begin{Lemma}\label{lemma:haus}
 Given an integrable system $\is$, if{\samepage
 \begin{enumerate}[label={\rm (\alph*)}, ref=(\alph*), leftmargin=*]\itemsep=0pt
 \item \label{item:13} its leaves are compact and locally connected, or
 \item \label{item:14} its fibers are connected,
 \end{enumerate}
 then the associated leaf space $\mathcal{L}$ is Hausdorff.}
 \end{Lemma}

 \begin{proof} Let $q\colon M \to \mathcal{L}$ denote the quotient map and let $\pi \colon \mathcal{L} \to \Phi(M)$ be the induced continuous map. Suppose that~\ref{item:13} holds and let $L_1,L_2 \subset M$ be distinct leaves of~$\is$. Since, for $i=1,2$, $L_i$ is compact, it is closed; moreover, since $M$ is a normal topological space (being a~smooth manifold), there exist disjoint open sets $V_1, V_2 \subset M$ such that, for $i=1,2$, $L_i \subset V_i$. Observe that the map $\Phi \colon M \to \R^n$ satisfy the hypotheses of \cite[Corollary~3.4]{mrcun}. Therefore, for $i=1,2$, the open set $V_i$ contains an open neighborhood $\hat{V}_i$ of $L_i$ that is saturated with respect to $q$. In particular, the subsets $q\big(\hat{V}_1\big), q\big(\hat{V}_2\big) \subset \mathcal{L}$ are open, disjoint and, for $i=1,2$, $[L_i] \in q\big(\hat{V}_i\big)$, where $[L_i]$ denotes the equivalence class of any point on $L_i$. This shows that $\mathcal{L}$ is Hausdorff. On the other hand, suppose that~\ref{item:14} holds. Then the map $\pi \colon \mathcal{L} \to \Phi(M)$ is injective. Continuity of $\pi$ and Hausdorffness of $\Phi(M)$ imply that $\mathcal{L}$ is also Hausdorff, as desired.
 \end{proof}

 \begin{rmk}\label{rmk:zung_haus} Lemma \ref{lemma:haus} should be compared with \cite[Proposition~3.3]{zung-symplectic} whose hypotheses are that the leaves of the integrable system be compact and that all singular points be
 non-degenerate (see Section~\ref{sec:non-degen-sing} and references therein). Using the linearization theorems for non-degenerate singular orbits (cf.~\cite{duf_mol,eliasson-thesis,miranda-zung}), it can be shown that non-degeneracy of all singular points implies local connectedness of leaves. However, as observed by Zung \cite[p.~187]{zung-symplectic}, there are integrable systems with degenerate orbits whose leaf spaces are Hausdorff. To the best of our knowledge, all such examples satisfy either condition \ref{item:13} or condition~\ref{item:14} of Lemma~\ref{lemma:haus}.
 \end{rmk}

A natural way to enhance the topological data encoded in the leaf space of an integrable system is to identify the singular leaves of the system.

\begin{Definition}\label{defn:singular_points_fibers} Let $\is$ be an integrable system with $n$ degrees of freedom with leaf space $\mathcal{L}$.
 \begin{itemize}[leftmargin=*]\itemsep=0pt
 \item A point $p \in M$ is {\em singular} if $\operatorname{rk} D_p \Phi < n$. Otherwise it is {\em regular}.
 \item A leaf of $\Phi$ is {\em singular} if it contains a singular point. Otherwise it is {\em regular}.
 \item The subset $\mathcal{L}_{\mathrm{sing}} \subset \mathcal{L}$ consisting of the image of singular leaves of $\Phi$ is said to be the {\em singular part} of $\mathcal{L}$, while its complement $\mathcal{L}_{\mathrm{reg}}$ is said to be the {\em regular part}.
 \end{itemize}
\end{Definition}

\begin{Definition}\label{defn:pair_leaf_regular} The pair $(\mathcal{L}, \mathcal{L}_{\mathrm{reg}})$ associated to an integrable system $\is$ is called the pair of {\em leaf and regular leaf spaces} of the system.
\end{Definition}

The above association descends to isomorphism classes of systems and behaves well with respect to saturated subsystems (see Lemma~\ref{lemma:leaf_space_invariant}):

\begin{Lemma}\label{lemma:pair_leaf_spaces}\quad
 \begin{enumerate}[label={\rm (\arabic*)}, ref= (\arabic*), leftmargin=*]
 \item \label{item:9} Isomorphic integrable systems possess homeomorphic pairs of leaf and regular leaf spaces.
 \item \label{item:10} The pair of leaf and regular leaf spaces of saturated subsystems of $\is$ naturally embed in the pair of leaf and regular leaf spaces of $\is$.
 \end{enumerate}
\end{Lemma}
\begin{proof} Let $(M_1,\omega_1,\Phi_1)$ and $(M_2,\omega_2,\Phi_2)$ be integrable systems that are isomorphic via the pair $(\Psi,\psi)$. Part~\ref{item:9} is an immediate consequence of the fact that $\Psi$ and $\psi$ are diffeomorphisms. The proof of~\ref{item:10} is analogous to the proof of part~\ref{item:12} of Lemma~\ref{lemma:leaf_space_invariant} and is left to the reader.
\end{proof}

\subsection{Faithful integrable systems}\label{sec:faithful-moment-maps}

In light of Lemma \ref{lemma:leaf_space_invariant}, it is helpful to distinguish those cases in which a moment map image can be identified, as a~topological space, with the leaf space. (To the best of our knowledge, the following notion does not appear elsewhere in the literature.)

\begin{Definition}\label{defn:faithful}
An integrable system $\is$ with leaf space $\mathcal{L}$ is said to be {\em faithful} if the induced map $\pi\colon \mathcal{L}\to B=\Phi(M)$
is a homeomorphism. Here, $B\subset\R^n$ is equipped with the subset topology.
\end{Definition}

Faithfulness implies connectivity of the fibers, but it is stronger than the latter (see Example~\ref{rmk:not_weaker} for an integrable system with connected fibers that fails to be faithful). Faithfulness also guarantees that the property of a subsystem being saturated coincides with there existing an open subset of the moment map image relative to which it is a~subsystem.

\begin{Lemma}\label{lemma:sat=rel} A subsystem of a faithful integrable system $\is$ is a subsystem relative to an open subset $U \subset \Phi(M)$ if and only if it is saturated.
 \end{Lemma}
\begin{proof} In light of Remark \ref{rmk:saturated_vs_relative}, it suffices to show that if $(V,\omega|_V, \Phi|_V)$ is a saturated subsystem then there is an open set of $\Phi(M)$ with respect to which $(V,\omega|_V, \Phi|_V)$ is a relative subsystem. We show that $U= \Phi(V)$ is the desired subset. Because $V$ is open and saturated with respect to the quotient map $q\colon M \to \mathcal{L}$, the subset $q(V) \subset \mathcal{L}$ is open. Since the system is faithful, $U = \pi(q(V))$, where $\pi \colon \mathcal{L} \to \Phi(M)$ is the induced map, is open in~$\Phi(M)$. Observe that
 \begin{gather*} 
 \Phi^{-1}(U) = q^{-1}\big(\pi^{-1}(U)\big) = q^{-1}(q(V)) = V,
 \end{gather*}
where the second equality follows from the fact that $\pi$ is injective, while the last equality follows from the fact that $V$ is saturated with respect to $q$.
 \end{proof}

The following result, which follows from Lemmas \ref{lemma:leaf_space_invariant} and~\ref{lemma:leaf_conn}, details how faithfulness behaves with respect to taking isomorphism classes and subsystems, and is stated below without proof.

\begin{Corollary}\label{cor:faithfulness_preserved} If an integrable system $\is$ is faithful, then so is every integrable system isomorphic to it, and every saturated subsystem. In particular,
\begin{itemize}[leftmargin=*]\itemsep=0pt
 \item faithful integrable systems form a full subcategory of integrable systems, and
 \item if $U \subset \Phi(M)$ is open, the set of connected components of $U$ is in bijective correspondence with the set of connected components of the total space of the subsystem relative to~$U$.
 \end{itemize}
\end{Corollary}

If the fibers of an integrable system are compact, necessary and sufficient conditions for faithfulness can be phrased without reference to the leaf space.

\begin{Lemma}\label{lemma:faithful_general} An integrable system $\is$ with compact fibers is faithful if and only if $\Phi$ has connected fibers and is proper onto its image.
\end{Lemma}

\begin{proof} Let $\mathcal{L}$ and $B$ be the leaf space and moment map image, respectively, of~$\is$. Let $q\colon M\to\mathcal{L}$ be the quotient map and let $\pi\colon \mathcal{L}\to B$ be the induced map. Suppose first that $\Phi$ has connected fibers and is proper onto its image. By definition of $\mathcal{L}$, the continuous map $\pi$ is a bijection because the fibers of~$\Phi$ are connected. It remains to show that $\pi$ is a closed map. In general, a continuous proper map to a metrizable space is closed (cf.\ Palais~\cite{palais-proper}). To see that $\pi$ is proper, consider an arbitrary compact set $K\subset B$. The preimage $\Phi^{-1}(K)$ is compact because~$\Phi$ is proper onto its image. Furthermore, $\pi^{-1}(K)=q\big(\Phi^{-1}(K)\big)$ because $q$ is surjective, so $\pi^{-1}(K)$ is compact because $q$ is continuous and $\Phi^{-1}(K)$ is compact. Therefore $\pi$ is proper. Then, since $B\subset \R^2$ is metrizable, the map $\pi$ is also closed. Consequently, $\pi$ is a~homeomorphism.

Conversely, suppose that $\pi$ is a homeomorphism; in particular, it is a bijection, which implies that $\Phi$ has connected fibers. (Thus, by Lemma~\ref{lemma:haus}, $\mathcal{L}$ is Hausdorff.) It remains to prove that~$\Phi$ is proper onto its image; since~$\pi$ is a homeomorphism, it suffices to check that $q$ is proper. To this end, note that, since $\pi$ is a homeomorphism, $\mathcal{L}$ is second countable and metrizable; moreover, $M$ is locally compact, the fibers of $q$ are Hausdorff, and~$q$ has compact and connected fibers. Therefore the result of Mr\v{c}un \cite[Theorem~3.3]{mrcun} can be applied: any open neighborhood of any given fiber of $q$ (= leaf of the system) contains an open neighborhood that is the union of compact connected components of fibers of~$q$. Connectedness of the fibers of~$q$ implies that this neighborhood is $q$-saturated. Therefore, any open neighborhood of a fiber of $q$ contains a~$q$-saturated neighborhood of the fiber. Arguing as in del Hoyo \cite[Proposition~2.1.3]{del_hoyo}, it follows that~$q$ is `sequentially proper' at any point, i.e., if $\{x_n\} \subset M$ is a sequence such that $\{q(x_n)\} \subset \mathcal{L}$ converges, then $\{x_n\}$ has a converging subsequence. Since for both~$M$ and~$\mathcal{L}$ compactness is equivalent to sequential compactness, the above `sequential properness' implies properness of~$q$ as desired.
\end{proof}

Following the ideas of Symington \cite{symington}, the moment map images of faithful almost-toric systems (see Definition~\ref{defn:almost-toric}) can be used to infer topological information regarding their total spaces. (More generally, if the singular orbits of an integrable system are understood, we would expect the above philosophy to extend to more general families of faithful systems.) Examples~\ref{rmk:not_weaker} and~\ref{exm:no_proper} illustrate two ways in which the above principle breaks down in the absence of faithfulness~-- first if the fibers of the moment map are not connected, and second if the moment map is not proper onto its image.

\begin{figure}[h] \centering
 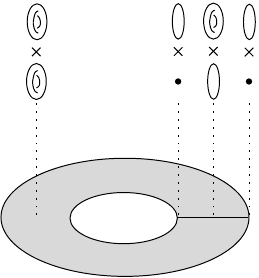
 \caption{Moment map image of the integrable system $(M,\omega,g_2\circ\Phi)$ of Example~\ref{rmk:not_weaker} showing the fibers over the horizontal segment $[1,e] \times \{0\}$. (For visual reasons, we have replaced the symbol for the union with $\times$.)} \label{LABEL2}
\end{figure}

\begin{exm}\label{rmk:not_weaker} Let $A\subset\R^2$ be the closed annulus centered at the origin with inner and outer radii~$1$ and $e$, i.e.,
\begin{gather*}
 A = \big\{(a,b) \in \R^2 \,|\, 1 \leq a^2 + b^2 \leq e^2\big\}.
\end{gather*}
Following Symington \cite{symington}, if $A$ is the moment map image of a faithful integrable system with elliptic singularities (see Definition~\ref{defn:cihs_elliptic}), then the total space of the system must be~$S^2 \times \mathbb{T}^2$. Such a system can be constructed as follows. Consider the toric system (see Definition~\ref{defn:toric_sys}) whose underlying symplectic toric manifold is
\begin{gather*}
 \big(S^2 \times S^1 \times \R, \, \operatorname{pr}^*_1 \omega_{S^2} + \operatorname{pr}^*_2\omega_{S^1 \times \R}, \, \tilde{\mu}:= \operatorname{pr}^*_1\mu + \operatorname{pr}^*_2 \zeta\big),
\end{gather*}
where $S^2\subset \R^3$ is the standard unit sphere in Euclidean space, $\omega_{S^2}$ is the standard area form on~$S^2$, $\mu \colon S^2 \to \R$ is the height function, $\omega_{S^1 \times \R} = d\theta \wedge d \zeta$, where $\theta$ is a mod~1 coordinate on~$S^1$ and~$\zeta$ is the standard coordinate on~$\R$, $\zeta \colon S^1 \times \R \to \R$ is projection onto the second component, and $\operatorname{pr}_1 \colon S^2
\times S^1 \times \R \to S^2$ and $\operatorname{pr}_2 \colon S^2 \times S^1 \times \R \to S^1 \times \R$ are the standard projections. The group $\Z$ acts freely and properly on $ S^2 \times S^1 \times \R$ by translations in the last component; the quotient by this action is diffeomorphic to $S^2 \times \mathbb{T}^2$ (this identification is henceforth used tacitly), and inherits a symplectic form $\omega_{S^2 \times \mathbb{T}^2}$. Define $g_k\colon [-1,+1] \times \R \to A\subset\R^2$ by
\begin{gather*}
 g_k(x,y)=\big(e^{\frac{1}{2}(x+1)}\cos (k\pi(y + 1) ), e^{\frac{1}{2}(x+1)}\sin(k\pi(y + 1))\big),
\end{gather*}
where $k \in \N$. The map $g_1 \circ \tilde{\mu}$ descends to a smooth map $\bar{\Phi}$ and the triple $\big(S^2 \times \mathbb{T}^2, \omega_{S^2 \times \mathbb{T}^2}, \bar{\Phi}\big)$ is a faithful system with elliptic singularities (see Definition~\ref{defn:cihs_elliptic}) whose moment map image is~$A$.

To illustrate how the failure of faithfulness can disrupt the relationship between the topology of the total space and that of the moment map image of an integrable system, we construct the following infinite family of integrable systems on $S^2 \times S^2$ whose moment map image equals $A$. Let $\is$ be the toric system that underlies the compact symplectic toric manifold
\begin{gather*}
 \big(S^2 \times S^2, \operatorname{pr}^*_1\omega_{S^2} + \operatorname{pr}^*_2\omega_{S^2}, \hat{\mu}:=\operatorname{pr}^*_1\mu + \operatorname{pr}^*_2\mu\big),
\end{gather*}
where $S^2$, $\omega_{S^2}$ and $\mu\colon S^2 \to \R$ are as above, and, for $i=1,2$, $\operatorname{pr}_i \colon S^2 \times S^2 \to S^2$ is projection onto the $i$th component. The moment map image $\hat{\mu}\big(S^2 \times S^2\big)$ is the square $R:=[-1,1] \times [-1, 1] \subset \R^2$. Fix $k \in \N$. The moment map of the integrable system $(M,\omega, g_k\circ\Phi)$ is proper (onto its image) as the total space is compact, and its image equals $A$ by construction. However, $(M,\omega, g_k\circ\Phi)$ fails to be faithful, as the fibers of $g_k\circ\Phi$ are not necessarily connected, for the fiber $(g_k \circ \Phi)^{-1}(a,b)$ is given by:
\begin{itemize}[leftmargin=*]\itemsep=0pt
 \item[$\bullet$] $k$ disjoint copies of $S^1$ if $a^2 + b^2 \in \big\{1,e^2\big\}$ and $(a,b) \notin \{(1,0),(e,0)\}$,
\item[$\bullet$] $2$ points and $k-1$ disjoint copies of $S^1$ if $(a,b) \in \{(1,0),(e,0)\}$,
 \item[$\bullet$] $2$ disjoint copies of $S^1$ and $k-1$ disjoint copies of $\mathbb{T}^2$ if $(a,b) \in \{(t,0) \,|\, 1 < t <e \}$, and
 \item[$\bullet$] $k$ disjoint copies of $\mathbb{T}^2$ otherwise.
\end{itemize}
In particular, in spite of being an integrable system with elliptic singularities, $(M,\omega, g_k\circ\Phi)$ has singular fibers in the preimage of the interior of the moment map image $A$ (see Fig.~\ref{LABEL2}), unlike any faithful system with elliptic singularities (see Remark~\ref{rmk:loc_tor_faithful}). Moreover, the topology of the total space of the above family of integrable systems with elliptic singularities differs significantly from that of the family of the faithful integrable system constructed above.
\end{exm}

\begin{exm}\label{exm:no_proper} With notation as in Example \ref{rmk:not_weaker}, let $(M',\omega',\Phi')$ be the subsystem of $\is$ relative to the subset $[-1,1] \times [-1, 1[$. Observe that $M'$ is diffeomorphic to $S^2 \times \R^2$. By construction, the integrable system $(M',\omega',g_1\circ\Phi')$ has connected fibers and its moment map image equals $A$. However, its moment map is not proper onto its image, as the preimage of $A$ is not compact. Unlike any faithful system with compact moment map image, the total space of the above system is not compact.
 \end{exm}

\subsection{$\Z$-affine structure on the regular part of the leaf space}\label{sec:z-affine-structure}
When the fibers of an integrable system are compact, the regular leaf space inherits a {\em geometric} structure. This is a consequence of
the Liouville--Arnol'd theorem, which provides a local normal form for a neighborhood of a regular leaf (cf.\ Cushman and Bates \cite[Appendix~D]{cush_bat}, Duister\-maat~\cite{duistermaat}, Guillemin and Sternberg \cite[Chapter~44]{gui_ster_st}, Sepe and \vungoc\ \cite[Theorem~3.36]{sepe_vu_ngoc} for various versions of a proof). Let $\Omega$ be the canonical symplectic form on $\cotan \mathbb{T}^n \cong\mathbb{T}^n\times\R^n$ for which the projection $\operatorname{pr}_2 \colon (\cotan\mathbb{T}^n,\Omega)\to \R^n$ defines an integrable system.

\begin{Theorem}[Liouville--Arnol'd]\label{thm:la} Let $\is$ be an integrable system with $n$ degrees of freedom and let~$F $ be a regular, compact leaf. Then there exist open neighborhoods $V\subset (M,\omega)$ of $F$ and $W\subset (\cotan \mathbb{T}^n,\Omega)$ of $\mathbb{T}^n \times \{0\}$, the latter saturated w.r.t.\ $\operatorname{pr}_2$, such that the subsystems of $\is$ and of $\big(\cotan \mathbb{T}^n,\Omega,\operatorname{pr}_2\big)$, relative to~$V$ and~$W$, respectively, are isomorphic via a pair $(\Psi,\psi)$ where $\psi(\Phi(F)) = 0$.
\end{Theorem}

\begin{rmk}\label{rmk:choice_g} Identify the Lie algebra of $\mathbb{T}^n$ with $\R^n$. Using the notation of Theorem~\ref{thm:la}, the composition $\psi \circ \Phi|_V$ is the moment map of a free, effective Hamiltonian $\mathbb{T}^n$-action, i.e., the Hamiltonian vector fields of its components have $2\pi$-periodic flows. Moreover, if $\Phi = (H_1,\ldots,H_n)$ with $H_1$ being the moment map of an effective Hamiltonian $S^1$-action, then the diffeomorphism~$\psi$ can be taken to be of the form
 \begin{gather}
 \psi(x_1,\ldots,x_n) = \big(\psi^{(1)}, \dots, \psi^{(n)}\big)(x_1, \dots, x_n) \nonumber\\
\hphantom{\psi(x_1,\ldots,x_n)}{} = \big(x_1 + a, \psi^{(2)}(x_1,\ldots,x_n),\ldots, \psi^{(n)}(x_1,\ldots,x_n)\big), \label{eq:1}
 \end{gather}
where $a \in \R$ is a constant.
\end{rmk}

\begin{Corollary}\label{cor:regular_leaf_space} Let $\is$ be an integrable system with compact fibers. The regular leaf space $\mathcal{L}_{\mathrm{reg}} \subset \mathcal{L}$ is open and inherits a structure of smooth manifold uniquely defined by requiring that the restriction of the quotient map $q|_{q^{-1}(\mathcal{L}_{\mathrm{reg}})}$ be a submersion onto $\mathcal{L}_{\mathrm{reg}}$. In particular, the restriction $\pi|_{\mathcal{L}_{\mathrm{reg}}}$ is smooth.
\end{Corollary}

\begin{proof} Moerdijk and Mr\v{c}un prove, in~\cite[Section 2.4]{moerdijk_mrcun}, that the leaf space of a submersion whose fibers are compact can be endowed with the structure of a smooth manifold uniquely defined by demanding that the quotient map be a submersion. The version of the Liouville--Arnol'd theorem given in Theorem~\ref{thm:la} implies that, for any point $p\in \mathcal{L}_{\mathrm{reg}}$, the corresponding regular leaf $F_p$ has an open neighborhood that is saturated by regular leaves. Therefore, $\Phi|_{q^{-1}(\mathcal{L}_{\mathrm{reg}})}$ is a~submersion whose leaf space is naturally isomorphic to $\mathcal{L}_{\mathrm{reg}}$, an open subset of $\mathcal{L}$. Since the fibers of $\Phi$ are compact by hypothesis, the result of Moerdijk and Mr\v{c}un implies the desired result.
\end{proof}

A symplectomorphism $\varphi$ of $\big(\T^n\times\R^n,\Omega\big)$ that preserves the fibers of $\operatorname{pr}_2$ must have the form $\varphi=\big(\varphi^{(1)},\varphi^{(2)}\big)$ where $\varphi^{(1)}(t,x)=\big(A^{-1}\big)^T t+f(x)$ and $\varphi^{(2)}(t,x)=Ax+c$ for some $A\in{\rm GL}(n,\Z)$, some $c\in\R^n$, and a smooth function $f\colon \R^n\to\R^n$ such that the matrix $A^{-1}\frac{\partial f}{\partial x}$ is symmetric (cf.\ Symington \cite[Lemma~2.5]{symington}.) This implies that maps of the form $\psi \circ \pi$, where $\psi$ is as in Theorem~\ref{thm:la}, can be used to define a~$\Z$-affine atlas on~$\mathcal{L}_{\mathrm{reg}}$.

\begin{Definition}\label{defn:iam} For any $n \geq 1$, the category of $n$-dimensional {\em $\Z$-affine manifolds}, denoted by $\mathcal{A}{\rm ff}_{\Z}(n)$, has objects and morphisms as follows:
 \begin{itemize}[leftmargin=*]\itemsep=0pt
 \item {\bf Objects}: {\em $\Z$-affine manifolds} (with corners), as defined in the Notation and conventions section.
 \item {\bf Morphisms}: {\em $\Z$-affine maps}, i.e., maps $f\colon (N_1,\mathcal{A}_1)\to (N_2,\mathcal{A}_2)$ that are given in local $\Z$-affine coordinates, by restrictions of elements of $\mathrm{AGL}(n;\Z)=\mathrm{GL}(n;\Z) \ltimes \R^n$.
 \end{itemize}
\end{Definition}

If $(N_2,\mathcal{A}_2)$ is a $\Z$-affine manifold and $f\colon N_1 \to N_2$ is locally a homeomorphism then there exists a unique (up to isomorphism) $\Z$-affine structure $\mathcal{A}_1$ on $N_1$ that makes $f$ into a $\Z$-affine morphism. The structure~$\mathcal{A}_1$ is henceforth referred to as being {\em induced} by~$f$.

\begin{exm}\label{exm:z-affine} For any $n\ge 1$, denote by $\mathcal{A}_0$ both the $\Z$-affine structure on~$\R^n$ and the $\Z$-affine structure with corners on $[0,\infty[\,^n$ obtained by declaring the standard coordinates $x_1,\ldots,x_n$ to be $\Z$-affine. Then the {\em standard $\Z$-affine structure} on an open subset of~$\R^n$ is the $\Z$-affine structure induced by inclusion of the subset in~$\big(\R^n,\mathcal{A}_0\big)$. Likewise, an open subset of the subspace $[0,\infty[\,^n$ also inherits the standard $\Z$-affine structure from inclusion in $([0,+\infty[\,^n,\mathcal{A}_0)$.
\end{exm}

For $\Z$-affine manifolds, it makes sense to consider (the sheaf of) {\em $\Z$-affine functions}, i.e., (locally defined) smooth functions that, in local $\Z$-affine coordinates $(x_1,\ldots,x_n)$, are given by
\begin{gather*}
 \sum\limits_{i=1}^n k_i x_i + c,
\end{gather*}
where $k_i \in \Z$ and $c\in\R$. The local normal form provided by the Liouville--Arnol'd theorem (Theorem~\ref{thm:la}) implies that the regular leaf space of an integrable system with compact fibers can be characterized, as a $\Z$-affine manifold, by the sheaf of functions that generate $2\pi$-periodic flows tangent to the fibers of the quotient map.

\begin{Corollary}\label{cor:reg_iam} Let $\is$ be an integrable system with compact fibers. Then the subset $\mathcal{L}_{\mathrm{reg}} \subset \mathcal{L}$ inherits a $\Z$-affine structure $\mathcal{A}_{\mathrm{reg}}$, uniquely defined by the property that locally defined $\Z$-affine functions from $(\mathcal{L}_{\mathrm{reg}},\mathcal{A}_{\mathcal{L}_\mathrm{reg}})$ to $(\R,\mathcal{A}_0)$ correspond, by taking the pull-back along the restriction to $q^{-1}(\mathcal{L}_{\mathrm{reg}})$ of the quotient map $q \colon M \to \mathcal{L}$, to functions on $q^{-1} (\mathcal{L}_{\mathrm{reg}} ) \subset M$ whose Hamiltonian vector fields are tangent to the fibers of $q$ and have $2\pi$-periodic flows.
\end{Corollary}

\begin{Corollary}\label{cor:functor_is_aff} For each $n \geq 1$, there is a functor $\mathcal{IS}(n) \to \mathcal{A}{\rm ff}_{\Z}(n)$ that, on objects, is precisely the map $\is \mapsto (\mathcal{L}_{\mathrm{reg}},\mathcal{A}_{\mathrm{reg}})$ given by Corollary~{\rm \ref{cor:reg_iam}}.
\end{Corollary}

\begin{proof} The above functor is completely determined by the following property, which can be checked directly: an isomorphism of integrable systems induces a $\Z$-affine (iso)morphism between the associated $\Z$-affine manifolds.
\end{proof}

Furthermore, the correspondence between integrable systems and $\Z$-affine manifolds given by Corollary \ref{cor:functor_is_aff} behaves well under restriction to saturated subsystems.

\begin{Corollary}\label{cor:is_ias_subs} Given an integrable system $\is$, the natural inclusion of the leaf space of a saturated subsystem into the leaf space of $\is$ corresponds to a $\Z$-affine embedding of the regular leaf space of the former into the regular leaf space of the latter.
\end{Corollary}

Finally, the above discussion allows further refinement of the set of invariants that can be associated to an integrable system. Given $\is$, associate the pair $(\mathcal{L},(\mathcal{L}_{\mathrm{reg}},\mathcal{A}_{\mathrm{reg}}))$ to it, where $(\mathcal{L},\mathcal{L}_{\mathrm{reg}})$ is the pair of leaf and regular leaf spaces of $\is$ and $\mathcal{A}_{\mathrm{reg}}$ is the $\Z$-affine structure given by Corollary~\ref{cor:reg_iam}. This association descends to isomorphism classes of systems. In this case, isomorphisms of pairs are homeomorphisms of the underlying topological pairs that restrict to $\Z$-affine isomorphisms on the $\Z$-affine subspace.

\begin{rmk}\label{rmk:z-affine_faithful} For faithful integrable systems $\is$, the leaf space $\mathcal{L}$ can be identified topologically with the moment map image $B = \Phi(M)$. Under this correspondence, $\mathcal{L}_{\mathrm{reg}}$ is identified with the subset of regular values $B_{\mathrm{reg}} \subset B$. By Corollary~\ref{cor:reg_iam}, $B_{\mathrm{reg}}$ inherits a $\Z$-affine structure denoted by $\mathcal{A}_{\mathrm{reg}}$ which, in general, is {\em not isomorphic} to the standard one as a~subset of~$\R^n$.
\end{rmk}

\subsection{Developing maps}\label{sec:developing-maps}

Given an $n$-dimensional $\Z$-affine manifold (with corners) $(N,\mathcal{A})$, let $\tilde{N}$ denote its universal cover. The universal covering map $\mathsf{q} \colon \tilde{N} \to N$ induces a $\Z$-affine structure $\tilde{\mathcal{A}}$ on $\tilde{N}$, making $\mathsf{q} \colon (\tilde{N},\tilde{\mathcal{A}}) \to (N,\mathcal{A})$ into a $\Z$-affine morphism (cf.\ Goldman and Hirsch~\cite{gh} for details). Identify~$\tilde{N}$ with the space of paths starting at $\mathsf{x}_0$, up to homotopy relative to endpoints. Fix a~basepoint $\mathsf{x}_0 \in N$, a point $\tilde{\mathsf{x}}_0 \in \tilde{N}$ with $\mathsf{q}(\tilde{\mathsf{x}}_0) = \mathsf{x}_0$, and a $\Z$-affine coordinate chart $\phi_0 \colon U_0 \to \R^n$ defined near $\mathsf{x}_0$. By shrinking $U_0$ if necessary, it may be assumed that there exists an open neighborhood~$\tilde{U}_0$ of~$\tilde{\mathsf{x}}_0$ such that $\mathsf{q}|_{\tilde{U}_0} \colon \tilde{U}_0 \to U_0$ is a~diffeomorphism. Then $\phi \circ \mathsf{q}|_{\tilde{U}_0}\colon \tilde{U}_0 \to \R^n$ defines a $\Z$-affine chart for $(\tilde{N},\tilde{\mathcal{A}})$. This map can be extended to a smooth map $\mathrm{dev} \colon \tilde{N} \to \R^n$ (cf.\ Goldman and Hirsch~\cite{gh} and references therein for a proof and further details). Intuitively, developing maps are constructed by means of `analytic continuation' of a $\Z$-affine chart). Such a map is called the {\em developing map} of $(N,\mathcal{A})$ (relative to the choices $(\mathsf{x}_0, \tilde{\mathsf{x}}_0, \phi_0)$). Moreover, there is a~representation $\mathfrak{a} \colon \pi_1(N;\mathsf{x}_0) \to \mathrm{AGL}(n;\Z)$, called {\em the affine holonomy of $(N,\mathcal{A})$}, which is intertwined with $\mathrm{dev}$ as follows: for all $\gamma \in \pi_1(N;\mathsf{x}_0)$, the following diagram commutes
 \begin{gather*}
 \xymatrix{(\tilde{N},\tilde{\mathcal{A}}) \ar[r]^-{\mathrm{dev}} \ar[d]_-{\cdot
 \gamma} & (\R^n,\mathcal{A}_0) \\
 (\tilde{N},\tilde{\mathcal{A}}) \ar[r]_-{\mathrm{dev}} & (\R^n,\mathcal{A}_0)
 \ar[u]_-{\mathfrak{a}(\gamma)},}
 \end{gather*}
where $\cdot \gamma$ denotes the $\Z$-affine isomorphism of $(\tilde{N},\tilde{\mathcal{A}})$ induced by the natural action of $\pi_1(N;\mathsf{x}_0)$ on~$\tilde{N}$. Note that, using the fundamental groupoid of $N$, the information of a~developing map and of the affine holonomy can be packaged and conveyed independent of choices (cf.\ Crainic, Fernandes and Mart\'inez-Torres~\cite{cfmt}).

\begin{rmk}\label{rmk:diff_choices_dev}\quad
\begin{itemize}[leftmargin=*]\itemsep=0pt
\item If $\mathrm{dev}, \mathrm{dev}' \colon \tilde{N} \to \R^n$ are developing maps for $(N,\mathcal{A})$ constructed using different choices then there exists a unique $h \in \mathrm{AGL}(n;\Z)$ such that $\mathrm{dev}' = h \circ \mathrm{dev}$ (cf.\ Goldman and Hirsch~\cite{gh}).
 \item If $(N,\mathcal{A})$ is a $\Z$-affine manifold with corners, then the image of any codimension-$k$ face of~$\tilde{N}$ ($0 < k \leq n$) under a developing map is the intersection of~$k$ linear hyperplanes of $\R^n$ whose normals can be chosen to span a~{\em unimodular} sublattice of~$\Z^n$, i.e., this span is a direct summand of~$\Z^n$. This can be seen as follows. By definition of $\Z$-affine chart on a~manifold with corners, the above statement holds locally. The way in which developing maps are constructed (cf.\ Goldman and Hirsch~\cite{gh}) implies that the local statement is sufficient to obtain the result globally.
 \item In general, developing maps need not be covering maps and their images can be rather complicated (cf.\ Sullivan and Thurston~\cite{sull_thur} for pathological examples).
 \end{itemize}
\end{rmk}

\subsection{(Visible) Toric systems}\label{sec:toric-systems}

This section describes the connection between integrable toric actions and integrable systems. Following Karshon and Lerman~\cite{karshon-lerman}, we begin by introducing the category of integrable toric actions on symplectic manifolds.

\begin{Definition}\label{defn:stm}
 For any $n \geq 1$, the category of {\em symplectic toric manifolds} of dimension $2n$, denoted by $\mathcal{TM}(2n)$, has objects and morphisms as follows:
 \begin{itemize}[leftmargin=*] \itemsep=0pt
 \item {\bf Objects}: {\em symplectic toric manifolds}, i.e., $2n$-dimensional symplectic manifolds~$\sm$ endowed with an effective Hamiltonian $\mathbb{T}^n$-action with moment map $\mu \colon \sm \to \mathfrak{t}^*$, where~$\mathfrak{t}^*$ denotes the dual of the Lie algebra of~$\mathbb{T}^n$. (Observe that~$M$ is not required to be compact.) A symplectic toric manifold is henceforth denoted by the triple $(M,\omega,\mu)$ and, for brevity, referred to as a {\em toric manifold}.
 \item {\bf Morphisms}: {\em isomorphisms of symplectic toric manifolds}, i.e., given $(M_i,\omega_i,\mu_i)$ for $i=1,2$, a pair $(\Psi,\xi)$, where $\Psi\colon (M_1,\omega_1) \to (M_2,\omega_2)$ is a symplectomorphism, and $\xi \in \mathfrak{t}^*$, making the following diagram commute
 \begin{gather*}
 \xymatrix{(M_1,\omega_1) \ar[d]_-{\mu_1} \ar[r]^-{\Psi} &
 (M_2,\omega_2) \ar[d]^-{\mu_2} \\
 \mathfrak{t}^* \ar[r]_-{+ \xi} & \mathfrak{t}^*,}
 \end{gather*}
 where $+ \xi \colon \mathfrak{t}^* \to \mathfrak{t}^*$ denotes translation by $\xi$. Equivalently, $\Psi$ is $\mathbb{T}^n$-equivariant.
 \end{itemize}
\end{Definition}

Henceforth, for each $n \geq 1$, fix an isomorphism $\mathfrak{t}^* \cong \R^n$ so that the standard lattice in $\mathfrak{t}^*$ (dual to $\ker \big(\exp \colon \mathfrak{t} \to \mathbb{T}^n\big)$) is mapped to $\Z^n$.

Note that the components of the moment map $\mu \colon \sm \to \R^n$ Poisson commute and $\mu$ is a submersion almost everywhere, due to the Marle-Guillemin-Sternberg local normal form for Hamiltonian actions of compact Lie groups (cf.\ Guillemin, Sternberg and Marle \cite{GS2,marle}). Therefore, taking $\Phi = \mu$, call $\is$ the integrable system {\em underlying}~$\toric$. Because isomorphisms of toric manifolds induce isomorphisms of underlying integrable systems, for each $n\ge 1$, the function from $\mathcal{TM}(2n)$ to $\mathcal{IS}(n)$ that maps a toric manifold to its underlying integrable system defines a `forgetful functor' $\mathcal{F} \colon \mathcal{TM}(2n) \to \mathcal{IS}(n)$.

It is useful to identify the integrable systems that underlie symplectic toric manifolds and systems that are isomorphic to such.

\begin{Definition}\label{defn:toric_sys} An integrable system $\is$ is {\em toric} if there exists a toric manifold $(M,\omega,\mu)$ such that $\is = \mathcal{F}(M,\omega,\mu)$. An integrable system is {\em weakly toric} if it is isomorphic to a~toric one.
\end{Definition}

\begin{rmk}\label{rmk:toric_type} \vungoc\ introduced the notion of integrable systems of {\em toric type} in \cite[De\-fi\-nition~2.1]{vu-ngoc}. Such systems are necessarily weakly toric in the sense of Definition~\ref{defn:toric_sys}; however there are two important differences, namely:
 \begin{itemize}[leftmargin=*]\itemsep=0pt
 \item the moment map of an integrable system of toric type in the sense of \cite[Definition~2.1]{vu-ngoc} is necessarily proper, and
 \item an integrable system is of toric type in the sense of \cite[Definition~2.1]{vu-ngoc} if it is isomorphic as an integrable system to a~toric system via an isomorphism of the form $(\mathrm{id}, \psi)$.
 \end{itemize}
Thus weakly toric systems as in Definition \ref{defn:toric_sys} can be viewed as a slight generalization of systems of toric type.
 \end{rmk}

\begin{Example}\label{ex:reg_leaf_toric_nbhd} Let $\is$ be an integrable system with compact fibers and $q\colon M\to \mathcal{L}$ the quotient map to its leaf space. A $\Z$-affine coordinate chart on $\mathcal{L}_{\mathrm{reg}}\subset\mathcal{L}$ yields, pre-composing with the quotient map, the moment map of a locally defined effective (but not unique) Hamiltonian $\mathbb{T}^n$-action (see Corollary~\ref{cor:reg_iam}), and hence a toric system on the preimage by $q$ of the domain of the coordinate chart.
\end{Example}

\begin{rmk}\label{rmk:uniqueness} The image of $\mathcal{F}$ is a subcategory of~$\mathcal{IS}(n)$, but it is not {\em full}, i.e., there are morphisms between objects in the image of~$\mathcal{F}$ that are not the image of morphisms under~$\mathcal{F}$. For instance, if $h \in \mathrm{AGL}(n;\Z)$ is an element different from the identity, the toric manifolds $\toricphi$ and $(M,\omega,h\circ \Phi)$ are not isomorphic, but the underlying integrable systems are. In fact, the integral affine group $\mathrm{AGL}(n;\Z)$ completely captures the failure of the functor to be full: Suppose that $(M_1,\omega_1,\Phi_1)$ and $(M_2,\omega_2,\Phi_2)$ are isomorphic toric systems with isomorphism denoted by~$(\Psi, \psi)$. Then $\psi \colon \Phi_1(M_1) \to \Phi_2(M_2)$ is given, on each connected component of~$\Phi_1(M_1)$, by the restriction of an element in $\mathrm{AGL}(n;\Z)$. This follows from the fact that $\Psi \circ \Phi_1$ is the moment map of an effective Hamiltonian $\mathbb{T}^n$-action whose components Poisson commute with the components of~$\Phi$.
\end{rmk}

The classification by Delzant \cite{delzant} of {\em compact} toric manifolds $\toric$ is well-known. In particular, the manifold~$M$, the symplectic form $\omega$ and the moment map $\mu$ up to isomorphism are determined by the image of the moment map,~$\mu(M)$. Two properties of a~compact toric manifold $\toric$ that are important for the classification are that~$\mu$ has connected fibers and that~$\mu(M)$ is a convex polytope. These properties follow from work of Atiyah~\cite{atiyah} and Guillemin, Sternberg~\cite{guillemin-sternberg} and hold for the more general family of effective Hamiltonian $\mathbb{T}^k$-actions on compact $2n$-dimensional symplectic manifolds for any $0 < k \leq n$.

If the underlying symplectic manifold is {\em not} compact, neither of the above properties need hold (cf.\ Karshon and Lerman~\cite{karshon-lerman}). Nevertheless, Karshon and Lerman~\cite{karshon-lerman} achieve a~classification of these objects. The starting point is the fact that orbits are tori that have a neighborhood that can be put in normal form (cf.\ Guillemin and Sternberg~\cite{GS2}, Marle~\cite{marle}). This is the content of the following result, stated below without proof.

 \begin{Theorem}\label{thm:gsm} Let $\toric$ be a toric manifold. For each orbit $\mathcal{O}$ of dimension $k$, there exist
 \begin{itemize}[leftmargin=*]\itemsep=0pt
\item open neighborhoods $V \subset M$ of $\mathcal{O}$ and $W\subset\cotan \T^k \times \R^{2(n-k)}\cong\T^k\times\R^{2n -k}$ of $\T^k \times \{\mathbf{0}\}$;
\item a symplectomorphism $\Psi \colon (V,\omega) \to (W,\omega_{\mathrm{can}}\oplus\omega_0)$ sending $\mathcal{O}$ to $\T^k \times \{\mathbf{0}\}$;
\item an element $A \in \mathrm{GL}(n-k;\Z)$;
\item a translation $-\mu(\mathcal{O}) \colon \mu(U) \to \R^n$
 \end{itemize}
making the following diagram commute
 \begin{gather*}
 \xymatrix{ (V,\omega) \ar[r]^-{\Psi} \ar[d]_-{\mu} &
 (W,\omega_{\mathrm{can}}\oplus\omega_0) \ar[d]^-{\operatorname{pr}_2\oplus A\circ{\bf q}} \\
 \mu(U) \ar[r]_-{-\mu(\mathcal{O})} & \R^k\times\R^{n-k},}
 \end{gather*}
 where
 \begin{itemize}[leftmargin=*]\itemsep=0pt
 \item $\omega_{\mathrm{can}}$ is the canonical symplectic form on $T^*\T^k\cong \T^k\times \R^k$;
 \item $\operatorname{pr}_2\colon \T^k\times \R^k\to \R^k$ is projection onto the second factor;
 \item $\omega_0 = \sum\limits_{i=1}^{n-k} dx_i\wedge dy_i$ with respect to standard symplectic coordinates on $\R^{2(n-k)}$;
 \item ${\mathbf q}=\big(\frac{x_1^2+y_1^2}{2},\ldots, \frac{x_{n-k}^2+y_{n-k}^2}{2} \big)$.
 \end{itemize}
 \end{Theorem}

Observe that in the statement of Theorem~\ref{thm:gsm} the open set $W$ is saturated by the fibers of $\operatorname{pr}_2\oplus A \circ {\bf q}$.

A simple, yet important, consequence of Theorem \ref{thm:gsm} is the following fact, stated below without proof.

\begin{Corollary}\label{cor:orbits_leaves} Any $\mathbb{T}^n$-orbit in a toric manifold is a connected component of a fiber of the moment map. In particular, for a toric system, orbits and leaves coincide.
 \end{Corollary}

Corollary \ref{cor:orbits_leaves} sets toric systems apart from general integrable systems.

Another important ingredient in the classification of toric manifolds by Karshon and Lerman~\cite{karshon-lerman} is the following notion.

 \begin{Definition}\label{defn:orbital_mom_map} Let $\toric$ be a toric manifold. The {\em orbital moment map} of~$\toric$ is the unique continuous map $\chi\colon M/\mathbb{T}^n \to \R^n$ such that $\mu = \chi \circ q$, where $q \colon M \to M/\mathbb{T}^n$ is the quotient map.
 \end{Definition}

The orbital moment map is an essential invariant of toric manifolds; this is the content of the following result.

\begin{Proposition}\label{prop:toric_manifolds} Suppose that the toric manifolds $(M_1,\omega_1,\mu_1)$ and $(M_2,\omega_2,\mu_2)$ are isomorphic via $(\Psi,\xi)$. Then there exists a homeomorphism $\zeta \colon \mathcal{L}_1 \to \mathcal{L}_2$ between the orbit/leaf spaces of $(M_1,\omega_1,\mu_1)$ and $(M_2,\omega_2,\mu_2)$ respectively satisfying $\pi_2 \circ \zeta = + \xi \circ \pi_1$, where, for $i=1,2$, $\pi_i \colon \mathcal{L}_i \to \R^n$ is the orbital moment map of $(M_i,\omega_i,\mu_i)$.
\end{Proposition}

 \begin{proof} The argument is entirely analogous to the one in the proof of item \ref{item:11} of Lemma~\ref{lemma:leaf_space_invariant}.
 \end{proof}

In analogy with Corollary \ref{cor:reg_iam}, Theorem \ref{thm:gsm} implies the following result, stated below without proof.

\begin{Corollary}\label{cor:prop_toric} Let $\toric$ be a toric manifold with underlying toric system $\is$. Then
 \begin{enumerate}[label = {\rm \arabic*)}, ref= \arabic*), leftmargin=*]\itemsep=0pt
 \item \label{item:18} The orbit space $M/\mathbb{T}^n$ of $\toric$ is a $\Z$-affine manifold with corners uniquely characterized as in Corollary~{\rm \ref{cor:reg_iam}}.
 \item \label{item:19} $M/\mathbb{T}^n$ is canonically homeomorphic to the leaf space $\mathcal{L}$ of $\is$ so that the orbital moment map is identified with the continuous map $\pi \colon \mathcal{L} \to \R^n$. $($Throughout the rest of the statement, this identification is used tacitly.$)$
 \item \label{item:20} The image under $\pi$ of a codimension-$k$ face of~$\mathcal{L}$, where $0 < k \leq n$, is the intersection of~$k$ hyperplanes of $\R^n$ whose normals can be chosen to span a unimodular sublattice of~$\Z^n$, i.e., a sublattice $\Xi \subset \Z^n$ such that the quotient $\Z^n/\Xi$ does not have torsion $($see Remark {\rm \ref{rmk:diff_choices_dev})}.
 \item\label{item:21} With the above identification, the $\Z$-affine structure on $\mathcal{L}$, denoted by $\mathcal{A}_{\mathcal{L}}$, extends $\mathcal{A}_{\mathrm{reg}}$, i.e., the inclusion $(\mathcal{L}_{\mathrm{reg}},\mathcal{A}_{\mathrm{reg}}) \hookrightarrow (\mathcal{L},\mathcal{A}_{\mathcal{L}})$ is a~$\Z$-affine embedding.
 \item \label{item:22} The regular leaf space satisfies $\mathcal{L}_{\mathrm{reg}} = \mathcal{L} \smallsetminus \partial \mathcal{L}$.
 \end{enumerate}
 \end{Corollary}

The developing map $\mathrm{dev} \colon \tilde{\mathcal{L}} \to \R^n$ for $(\mathcal{L},\mathcal{A}_{\mathcal{L}})$ makes the following diagram commute
\begin{gather*}
 \xymatrix{ & \tilde{\mathcal{L}} \ar[dl]_-{\mathsf{q}}
 \ar[dr]^-{\mathrm{dev}} & \\
 \mathcal{L} \ar[rr]_-{\pi} & & \R^n,}
\end{gather*}
where $\mathsf{q} \colon \tilde{\mathcal{L}} \to \mathcal{L}$ is the universal covering map and, under the identification of part~\ref{item:19} of Corollary~\ref{cor:prop_toric}, $\pi$~is the orbital moment map of $\toric$ (cf.\ Karshon and Lerman \cite[Proposition~1.1, Remarks~1.4 and~1.5]{karshon-lerman}). In particular, if $\mathcal{L}$ is simply connected, this proves the following result.

\begin{Corollary}\label{cor:dev_map_sc} Let $\toric$ be a toric manifold with associated toric system~$\is$. If the leaf space $\mathcal{L}$ is simply connected, then the map $\pi \colon \mathcal{L} \to \R^n$ is a developing map for the induced $\Z$-affine structure.
 \end{Corollary}

The remaining ingredients in the classification up to isomorphism are topological invariants depending on $\mathrm{H}^2(\mathcal{L};\Z)$ (cf.\ Karshon and Lerman \cite[Theorem~1.3]{karshon-lerman} for a precise statement). Motivated by the classification theorem of Karshon and Lerman~\cite[Theorem~1.3]{karshon-lerman}, we introduce a class of toric systems whose corresponding toric manifolds are determined up to isomorphism by their moment map images.

\begin{Definition}\label{defn:delzant-system} A toric system $\is$ is said to be {\em visible} if it is faithful and $\Phi(M)$ is contractible.
\end{Definition}

The reason for the terminology in Definition~\ref{defn:delzant-system} is that, in order to reconstruct a visible toric system, it suffices to know (or `look at') its moment map image by Karshon and Lerman \cite[Theorem~1.3]{karshon-lerman}.

\begin{Lemma}\label{lemma:delzant_connected} The total space of a visible toric system is connected.
\end{Lemma}

\begin{proof} This follows from faithfulness of the moment map, connectedness of the moment map image, and the equivalence of the cardinality of the set of components of the leaf space and of the total space of an integrable system (see Lemma~\ref{lemma:leaf_conn}).
\end{proof}

The induced $\Z$-affine structure on the leaf space of a visible toric system is (isomorphic to) the standard one as a subset of~$\R^n$. More precisely, the following holds.

 \begin{Corollary}\label{rmk:z-aff_delz_sys} Let $\is$ be a visible toric system with $n$ degrees of freedom. Then its moment map image inherits an $\Z$-affine structure from~$\is$ that agrees with the standard one as a subset of~$\R^n$.
 \end{Corollary}
\begin{proof}Since $\is$ is visible, it is faithful in particular and, therefore, its moment map image $B$ can be canonically identified with its leaf space~$\mathcal{L}$ . The latter can in turn be canonically identified with the orbit space of the toric system~$\toric$ associated to~$\is$ (see part~\ref{item:19} of Corollary~\ref{cor:prop_toric}). By part~\ref{item:18} of Corollary~\ref{cor:prop_toric}, this implies that~$B$ inherits a $\Z$-affine structure. Unraveling the above identifications, Corollary~\ref{cor:dev_map_sc} implies that $\Z$-affine structure agrees with the standard one as a subset of~$\R^n$ as desired.
 \end{proof}

Using Karshon and Lerman \cite[Theorem 1.3]{karshon-lerman}, it is straightforward to obtain the classification of visible toric systems.

\begin{Proposition}\label{prop:ds} Two visible toric systems $(M_i, \omega_i, \Phi_i)$, $i =1,2$, are isomorphic if and only if there is an element $h\in\mathrm{AGL}(n;\Z)$ such that $\Phi_2(M_2)=h\circ\Phi_1(M_1)$. Furthermore, the two visible toric systems underlie isomorphic toric manifolds if and only if~$\Phi_1(M_1)$ and~$\Phi_2(M_2)$ agree up to translation.
\end{Proposition}

\begin{proof} Because $(M_i, \omega_i, \Phi_i)$, $i =1,2$, are toric systems, there exist toric manifolds $(M_i, \omega_i, \mu_i)$ such that $\mu_i=\Phi_i$, $i=1,2$. Under that equivalence, faithfulness of the moment maps $\Phi_i$ corresponds to the orbital moment maps of the toric manifolds being embeddings (this follows from part~\ref{item:19} of Corollary~\ref{cor:prop_toric}). Therefore, since $\mu_i(M_i)=\Phi_i(M_i)$, $i=1,2$, is contractible, Theorem~1.3 of Karshon and Lerman~\cite{karshon-lerman} implies the toric manifolds $(M_i,\omega_i,\mu_i)$ are determined up to isomorphism by their moment map images. By the definition of $\mathcal{TM}(2n)$, the category of symplectic toric manifolds (Definition~\ref{defn:stm}), $(M_i,\omega_i,\mu_i)$, $i=1,2$ belong to the same isomorphism class if and only if $\mu_1(M_1)$ and $\mu_2(M_2)$ differ by a translation. The criterion for isomorphism of the toric systems then follows from Remark~\ref{rmk:uniqueness} and connectedness of the total space (Lemma~\ref{lemma:delzant_connected}).
\end{proof}

Unlike {\em compact} toric manifolds (and their associated toric systems), the moment map image of a visible toric system need not be convex; however, Theorem~\ref{thm:gsm} together with the defining properties of visible toric systems, imply the following result.

\begin{Corollary}\label{cor:delz_sys_loc_convex} The moment map image of a visible toric system $\is$ is locally convex, i.e., for all $c \in \Phi(M)$, there exists an open neighborhood $U$ of $c$ in $\Phi(M)$ that is convex as a~subset of Euclidean space.
\end{Corollary}

\subsection{Weakly toric leaf spaces}\label{sec:toric-leaf}
Weakly toric systems provide examples of integrable systems whose leaf spaces are naturally endowed with the structure of a~$\Z$-affine manifold with corners. This is not a phenomenon to be expected in general. However it is natural to ask, what is the largest subset of the leaf space of an integrable system that does inherit the structure of a $\Z$-affine manifold with corners?

That question motivates the following notions of a {\em weakly toric leaf} and the {\em weakly toric leaf space} of an integrable system.

\begin{Definition}\label{defn:lt} Given an integrable system $\is$, a leaf $L \subset M$ is said to be {\em weakly toric} if there exists a~connected open neighborhood $V\subset M$ such that the subsystem $(V,\omega|_V,\Phi|_V)$ is weakly toric.
\end{Definition}

\begin{rmk}\label{rmk:lt} A weakly toric leaf $L$ is a smoothly embedded submanifold that must be compact even though the system $\is$ may have non-compact fibers. Also, the open neighborhood~$V$ of Definition~\ref{defn:lt} is saturated with respect to the quotient map $q \colon M \to \mathcal{L}$ and all leaves contained in~$V$ are weakly toric. Finally, in analogy with Remark~\ref{rmk:choice_g}, if $(\Psi,\psi)$ denotes the isomorphism between $(V,\omega|_V,\Phi|_V)$ and a toric system, and the first component of $\Phi|_V$ is the moment map of an effective Hamiltonian $S^1$-action, then $\psi$ can be taken to be of the form of equation \eqref{eq:1}.
\end{rmk}

\begin{Definition}\label{defn:lt_leaf} Given an integrable system $\is$ with leaf space $\mathcal{L}$, the subset $\mathcal{L}_{\mathrm{wt}} \subset \mathcal{L}$ corresponding to weakly toric leaves is called the {\em weakly toric leaf space} associated to $\is$.
 \end{Definition}

A priori, if the fibers of an integrable system are not necessarily compact, the associated weakly toric leaf space may be empty. In contrast, when the fibers are required to be compact, the Liouville--Arnol'd theorem (Theorem~\ref{thm:la}) implies the following result:

\begin{Corollary}\label{lemma:lt_contains_reg} If $\is$ has compact fibers, then $\mathcal{L}_{\mathrm{reg}} \subset \mathcal{L}_{\mathrm{wt}}$. In particular, $\mathcal{L}_{\mathrm{wt}} \subset \mathcal{L}$ is dense.
 \end{Corollary}

In fact, if an integrable system has compact fibers, its locally weakly toric leaf space inherits a $\Z$-affine structure.

\begin{Proposition}\label{cor:lt_iam} The weakly toric leaf space $\mathcal{L}_{\mathrm{wt}}$ of an integrable system with compact fibers inherits a structure of $\Z$-affine manifold with corners denoted by $\mathcal{A}_{\mathrm{wt}}$. This structure is uniquely defined by the property that locally defined $\Z$-affine functions from $(\mathcal{L}_{\mathrm{wt}},\mathcal{A}_{\mathcal{L}_\mathrm{wt}})$ correspond to functions on $q^{-1}(\mathcal{L}_{\mathrm{wt}}) \subset M$ whose Hamiltonian vector fields are tangent to the fibers of~$q$ and have $2\pi$-periodic flows.
\end{Proposition}

\begin{proof}Fix an integrable system $\is$ with $n$ degrees of freedom with compact fibers. To show that the weakly toric leaf space $\mathcal{L}_{\mathrm{wt}}$ is Hausdorff, suppose that $L_1$, $L_2$ are distinct weakly toric leaves. The aim is to show that, for $i=1,2$, there exists an open neighborhood~$V_i$ of~$L_i$, saturated with respect to $q$ and containing solely weakly toric leaves, with $V_1 \cap V_2 = \varnothing$. First, observe that $L_1,L_2 \subset M$ are closed; thus there exist $V_1', V_2' \subset M$ open, disjoint subsets with $L_i \subset V_i'$ for $i=1,2$. Since, for $i=1,2$, $L_i$ is weakly toric, it follows that $L_i \subset \Phi^{-1}(\Phi(L_i))$ is open. Therefore, using the result of Mr\v{c}un \cite[Theorem~3.3]{mrcun}, it may be assumed, without loss of generality, that $V_i'$ is saturated with respect to $\Phi$ and, therefore, with respect to $q$. Since, for $i=1,2$, $L_i$ is weakly toric, there exists an open neighborhood $V_i \subset V_i'$ of $L_i$ saturated with respect to $q$ containing solely weakly toric leaves; $V_1$ and $V_2$ are the desired separating open subsets. The space $\mathcal{L}_{\mathrm{wt}}$ is second countable because it is the image of a~second countable space under an open map. Indeed $q^{-1}(\mathcal{L}_{\mathrm{wt}}) \subset M$ is second countable as it is an open subset of a~smooth manifold, and the topological quotient map $q|_{ q^{-1}(\mathcal{L}_{\mathrm{wt}})} \colon q^{-1}(\mathcal{L}_{\mathrm{wt}}) \to \mathcal{L}_{\mathrm{wt}}$ is open because the quotient maps in the local models for toric manifolds are open. Finally, since $\mathcal{L}_{\mathrm{wt}}$ is locally homeomorphic to a subset of Euclidean space, it is locally compact, thereby implying~$\mathcal{L}_{\mathrm{wt}}$ is paracompact (by virtue of being a locally compact, Hausdorff, second countable space).

Next we define an open cover of $\mathcal{L}_{\mathrm{wt}}$ and coordinate charts whose codomain is $[0,+ \infty[\,^n$. Let $L \subset M$ be a~weakly toric leaf; by definition, there exists an open neighborhood $V \subset M$ of $L$ whose corresponding subsystem is weakly toric. By restricting $V$ if necessary, it may be assumed that the corresponding subsystem is isomorphic to a~local model for toric manifolds. Since $q|_{ q^{-1}(\mathcal{L}_{\mathrm{wt}})}$ is open, $q(V)$ is an open subset of $\mathcal{L}_{\mathrm{wt}}$; moreover, the above isomorphism implies that there is a map $\chi\colon q(V) \to [0,+\infty[^n$ that is locally a homeomorphism. Since $L \subset M$ is arbitrary, the above reasoning defines an open cover of $\mathcal{L}_{\mathrm{wt}}$, denoted by~$\{U_i\}$, and, for each $i$, a map $\chi_i \colon U_i \to [0,+\infty[\,^n$. In fact, $\mathcal{A}_{\mathrm{wt}}:=\{(U_i,\chi_i)\}$ is an $n$-dimensional $\Z$-affine atlas with corners. To see this, fix $i$, $j$ with $U_i \cap U_j \neq \varnothing$. Then, unraveling the above definitions, we obtain an isomorphism $(\Psi_{ij},\psi_{ij})$ of saturated subsystems of local models of toric manifolds with $\psi_{ij} = \chi_j \circ \chi_i^{-1}$; observing that $\psi_{ij}$ is necessarily the restriction of an element in $\mathrm{AGL}(n;\Z)$ shows the desired result. Finally, the defining property of $\mathcal{A}_{\mathrm{wt}}$ follows directly from its definition and the fact it holds for the (weakly toric) leaf spaces of toric manifolds.
\end{proof}

\begin{rmk}\label{sec:locally-weakly-toric} Note that, by Proposition \ref{cor:lt_iam}, the weakly toric leaf space of an integrable system is, tautologically, the largest subset of the leaf space of an integrable system that inherits the structure of a~$\Z$-affine manifold with corners from the integrable system. Moreover, if the system has compact fibers, the inclusion $(\mathcal{L}_{\mathrm{reg}},\mathcal{A}_{\mathrm{reg}}) \hookrightarrow (\mathcal{L}_{\mathrm{wt}},\mathcal{A}_{\mathrm{wt}})$ is a~$\Z$-affine morphism.
\end{rmk}

As expected, the $\Z$-affine manifold with corners $(\mathcal{L}_{\mathrm{wt}}, \mathcal{A}_{\mathrm{wt}})$ associated to~$\is$ is an invariant of the isomorphism class of~$\is$ and behaves well with respect to restriction to saturated subsystems, i.e., statements analogous to Corollaries~\ref{cor:functor_is_aff} and~\ref{cor:is_ias_subs} hold for this (possibly larger) $\Z$-affine manifold with corners. This allows one to associate to an integrable system~$\is$ the pair $(\mathcal{L},(\mathcal{L}_{\mathrm{wt}},\mathcal{A}_{\mathrm{wt}}))$, the latter being an invariant of the isomorphism class of $\is$.

\begin{rmk}\label{rmk:loc_tor_faithful} If $\is$ is faithful and has compact fibers, the weakly toric leaf space corresponds to an open, dense subset denoted by $B_{\mathrm{wt}} \subset B=\Phi(M)$ and its boundary as a~manifold with corners (corresponding to singular weakly toric leaves) satisfies $\partial B_{\mathrm{wt}} \subset B \cap \mathrm{Bdy}(B)$, where the inclusion may be strict.
\end{rmk}

\subsection{Cartographic maps}\label{sec:faithful}
Let $\is$ be a faithful integrable system with compact fibers. If the inclusion $B = \Phi(M)$ $\hookrightarrow \R^n$ is an $\Z$-affine embedding when restricted to $B_{\mathrm{wt}} \cong \mathcal{L}_{\mathrm{wt}}$ (as is the case for faithful toric systems), then the $\Z$-affine geometry of~$B_{\mathrm{wt}}$ is `captured' by the moment map image, i.e., it is the standard one as a subset of~$\R^n$. In this case, the subsystem relative to $B_{\mathrm{wt}}$ is toric; expanding the terminology introduced by Pelayo, Ratiu and \vungoc\ in \cite[Definition~4.2]{pvr_carto}, we introduce the following notion.

\begin{Definition}\label{defn:carto_mom_map} Let $\is$ be a faithful integrable system with compact fibers. The moment map is said to be {\em cartographic} if the inclusion $B = \Phi(M) \hookrightarrow \R^n$ is an $\Z$-affine embedding when restricted to $B_{\mathrm{wt}} \cong \mathcal{L}_{\mathrm{wt}}$.
 \end{Definition}

When the moment map of a faithful integrable system $\is$ with compact fibers is not cartographic, one could ask whether there is an isomorphic system whose moment map {\em is} cartographic. Existence of such a system is tantamount to finding a smooth embedding of the moment map image~$B$ into~$\R^n$ whose restriction to the weakly toric leaf space $B_{\mathrm{wt}}$ is a~$\Z$-affine embedding into $\big(\R^n,\mathcal{A}_0\big)$. A~necessary condition is that the affine holonomy of $(B_{\mathrm{wt}},\mathcal{A}_{\mathrm{wt}})$ be trivial; however, that condition is not sufficient, as the following example illustrates.

\begin{figure}[h] \centering
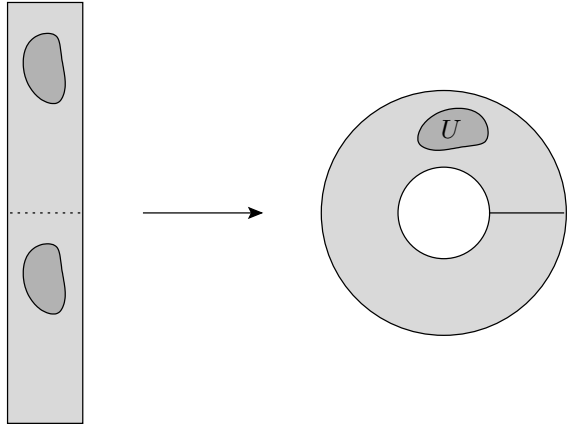
\caption{The preimage of the set $U$ under the map $g_2$ has two connected components.} \label{LABEL3}
\end{figure}

\begin{exm} For notation and further details, the reader is referred to Example~\ref{rmk:not_weaker} and Fig.~\ref{LABEL2}. Let $R$ be the open square $]{-}1,1[ \times\ ]{-}1, 1[ \subset \R^2$. Consider the unique toric manifold $\big(\mathbb{T}^2\times R,\omega,\mu\big)$ defined by the orbital moment map (see Definition~\ref{defn:orbital_mom_map}) $g_2\colon R\to\R^2$ given by
\begin{gather*}
 g_2(x,y)=\big(e^{\frac{1}{2}x + \frac{1}{2}}\cos (2\pi (y + 1 ) ),e^{\frac{1}{2}x + \frac{1}{2}}\sin (2\pi (y + 1 ) )\big).
\end{gather*}
Observe that $\omega$ is not the standard symplectic structure induced by inclusion of $ \T^2\times R$ in $\big( \T^2\times \R^2,\omega_{\rm can}\big)$. Moreover, if $\Phi\colon \mathbb{T}^2\times R\to R\subset\R^2$ denotes the projection onto the second factor, then $\mu = g_2 \circ \Phi$. In fact, $(\T^2\times R,\omega,\Phi)$ defines a faithful integrable system. Note that for any open simply connected subset $U\subset \{(x_1,x_2)\,|\, 1< x_1^2+x_2^2<e^2\}$, the subsystem of $\big(\mathbb{T}^2\times R,\omega,\mu\big)$ relative to~$U$ is the union of at least two disjoint subsystems, each of which is isomorphic to a subsystem of the toric system $\big(\T^2\times R,\omega,\Phi\big)$ (see Fig.~\ref{LABEL3}). However, the integrable system $\big(\mathbb{T}^2\times R,\omega,\Phi\big)$ as a whole is not isomorphic to any system with a~cartographic moment map as~$g_2$, which is not injective, is the unique map (up to composition on the right with an $\Z$-affine diffeomorphism of $(\R^2,\mathcal{A}_0)$) such that $g_2 \circ \Phi$ generates an effective Hamiltonian $\mathbb{T}^2$-action on $\big(\mathbb{T}^2 \times \R, \omega\big)$.
\end{exm}

The faithful integrable systems with compact fibers considered in Section \ref{sec:at} and Part~\ref{part:faithful-semitoric} allow for singular fibers (focus-focus fibers, see Section~\ref{sec:almost-toric-orbits}) that induce non-trivial affine holonomy on the weakly toric leaf space (see Theorem~\ref{thm:holonomy}). Thus for such an integrable system there is not necessarily a system in its isomorphism class that has a cartographic moment map. However, following the insight of Symington~\cite{symington} and \vungoc~\cite{vu-ngoc}, for such systems it is reasonable to ask whether there is a homeomorphism of the moment map image that, when restricted to an {\em open, dense} subset of the weakly toric leaf space, is a $\Z$-affine embedding. That question motivates introducing the following notion.

\begin{Definition}\label{defn:rh} Let $\is$ be a faithful integrable system with $n$ degrees of freedom and whose fibers are compact, and set $B= \Phi(M)$. A {\em cartographic pair} $(f,S)$ for $\is$ consists of a topological embedding $f \colon B \to f(B) \subset \R^n$, called {\em cartographic homeomorphism}, and an open, dense subset $S \subset B_{\mathrm{wt}}$, with the property that $f|_{S}\colon (S,\mathcal{A}_{\mathrm{wt}}|_S)\to(\R^n,\mathcal A_0)$ is a $\Z$-affine smooth embedding. If such a pair exists, $\is$ is said to {\em admit a~cartographic homeomorphism} whose image is said to be {\em cartographic}.
\end{Definition}

\begin{exm}\label{exm:toric} Any faithful toric system admits a natural choice of cartographic pair. Namely, if~$\is$ is one such system, then $(\mathrm{id},B)$ is a cartographic pair.
 \end{exm}

\begin{rmk}\label{rmk:smooth_locus} If $(f,S)$ is a cartographic pair for $\is$, there is no guarantee that $S$ is {\em maximal}, i.e., that it is not a strict subset of another open, dense subset of~$B_{\mathrm{wt}}$ on which~$f$ restricts to a $\Z$-affine embedding.
\end{rmk}

The following lemma provides a simple but useful way to adjust a given cartographic ho\-meo\-morphism.

\begin{Lemma}\label{lem:infinitely_many_carto} If $(f,S)$ is a cartographic pair for $\is$, then, for any $h \in \mathrm{AGL}(n;\Z)$, $(h \circ f,S)$ is also a cartographic pair.
\end{Lemma}

Intuitively, cartographic homeomorphisms should be thought of as continuous extensions of restrictions of developing maps to suitable domains; for instance, if~$S$ is simply connected, $S$~can be identified with a~dense subset of a fundamental domain in the universal cover of~$B_{\mathrm{wt}}$.

If $S=B$ holds in Definition~\ref{defn:rh}, then $f\circ\Phi$ is a~cartographic moment map. If not, a~cartographic homeomorphism at least provides a~dense subset of the total space on which the system is isomorphic to a toric, and possibly visible toric system. More precisely:

\begin{Lemma}\label{lemma:carto_toric} Let $\is$ be a faithful integrable system with compact fibers. If $(f,S)$ is a~cartographic pair, then the integrable system $\big(\Phi^{-1}(S),\omega|_{\Phi^{-1}(S)},(f \circ \Phi)|_{\Phi^{-1}(S)}\big)$ is toric. If, in addition, $S$ is contractible, then the above system is visible toric.
\end{Lemma}

Cartographic homeomorphisms restrict appropriately when taking saturated subsystems.

\begin{Lemma}\label{cor:carto_subsys} Suppose that $\is$ is a faithful integrable system with compact fibers with cartographic pair~$(f,S)$. Let $U \subset B = \Phi(M)$ be an open subset. Then $(f|_{U},S \cap U)$ is a~cartographic pair for the subsystem relative to~$U$.
\end{Lemma}

Combining Lemmas~\ref{lemma:carto_toric} and~\ref{cor:carto_subsys}, we obtain the following simple description of cartographic homeomorphisms when restricted to open subsets of the moment map image whose corresponding subsystems are visible toric.

\begin{Corollary}\label{cor:carto_Delzant} Let $\is$ be a faithful integrable system with compact fibers with cartographic pair~$(f,S)$. Suppose that $U \subset S$ is an open subset with the property that the subsystem relative to~$U$ is visible toric. Then~$f|_U$ is the restriction of an
 element $h_U \in \mathrm{AGL}(n;\Z)$ and $\big(h_U^{-1} \circ f,S\big)$ is a cartographic pair for~$\is$ with $\big(h_U^{-1} \circ f\big)|_U= \mathrm{id}|_U$.
\end{Corollary}

\begin{proof} By assumption, $\big(\Phi^{-1}(U),\omega|_{\Phi^{-1}(U)},\Phi|_{\Phi^{-1}(U)}\big)$ is a visible toric system and, by Lem\-ma~\ref{lemma:carto_toric}, $\big(\Phi^{-1}(U),\omega|_{\Phi^{-1}(U)},(f \circ \Phi)|_{\Phi^{-1}(U)}\big)$ is also visible toric. In fact, the pair $(\mathrm{id},f|_U)$ defines an isomorphism between these two systems. Because $U$ is connected (by virtue of being contractible), there exists an element $h_U \in \mathrm{AGL}(n;\Z)$ such that $f|_U =h_U$ (see Remark~\ref{rmk:uniqueness}). Lemma~\ref{lem:infinitely_many_carto} gives that $\big(h_U^{-1} \circ f,S\big)$ is a cartographic pair for~$\is$ that, by construction, satisfies $\big(h_U^{-1}\circ f\big)|_U = \mathrm{id}|_U$.
\end{proof}

Finally, it is important to notice that the property of admitting a cartographic homeomorphism is independent of the choice of representative of the isomorphism class of a faithful integrable system.

\begin{Corollary}\label{cor:pull_back_iso_carto} Let $(M_1,\omega_1,\Phi_1)$ and $(M_2,\omega_2,\Phi_2)$ be faithful integrable systems with compact fibers isomorphic via $(\Psi,\psi)$ and let $(f_1,S_1)$ be a cartographic pair for $(M_1,\omega_1,\Phi_1)$. Then $(f_2:=f_1 \circ \psi^{-1}, S_2: = \psi(S_1))$ is a cartographic pair for $(M_2,\omega_2,\Phi_2)$.
\end{Corollary}

In fact, cartographic images of isomorphic systems are homeomorphic via maps that extend $\Z$-affine isomorphisms. More precisely, the following holds.

\begin{Corollary}\label{cor:carto_iso} Consider, for $i=1,2$, a faithful integrable system $(M_i,\omega_i,\Phi_i)$ with $n$ degrees of freedom and compact fibers with a~cartographic pair~$(f_i,S_i)$. Assume that there exists an isomorphism $(\Psi,\psi)$ between $(M_1,\omega_1,\Phi_1)$ and $(M_2,\omega_2,\Phi_2)$. Then the cartographic images $f_1(B_1)$ and $f_2(B_2)$ are homeomorphic by a map that, when restricted to each connected component of $f_1\big(S_1\cap \psi^{-1}(S_2)\big)$, is the restriction of an element of~$\mathrm{AGL}(n;\Z)$.
\end{Corollary}

\begin{proof}Fix an isomorphism $(\Psi,\psi)\colon (M_1,\omega_1,\Phi_1)\to(M_2,\omega_2,\Phi_2)$. The map $g:= f_2 \circ \psi \circ f_1^{-1} \colon f_1(B_1) \to f_2(B_2)$ is a homeomorphism as it is the composition of homeomorphisms. In fact, we claim that it is the unique extension of a $\Z$-affine isomorphism $f_1\big(S_1 \cap \psi^{-1}(S_2)\big) \to f_2(\psi(S_1)\cap S_2)$. To see this, begin by observing that $S_1 \cap \psi^{-1}(S_2)$ and $ \psi(S_1) \cap S_2$ are open and dense in $B_1$ and $B_2$ respectively. This implies that for each $i=1,2$, $f_i\big(S_1 \cap \psi^{-1}(S_2)\big) $ is open and dense in $f_i(B_i)$ as $f_i$ is a~homeomorphism. Therefore, $g$ is determined uniquely by its restriction to $f_1\big(S_1 \cap \psi^{-1}(S_2)\big) $, which maps homeomorphically onto $f_2(\psi(S_1)\cap S_2)$. This restriction is a $\Z$-affine isomorphism, by the definition of cartographic homeomorphisms and Corollary~\ref{cor:functor_is_aff}. Since the $\Z$-affine structures on $f_1\big(S_1 \cap \psi^{-1}(S_2)\big) $ and on $f_2 (\psi(S_1)\cap S_2)$ are isomorphic to the ones induced by inclusion into $\big(\R^n,\mathcal{A}_0\big)$, the restriction of the above $\Z$-affine isomorphism to each connected component of $f_1\big(S_1 \cap \psi^{-1}(S_2)\big) $ is the restriction of an element of $\mathrm{AGL}(n;\Z)$.
\end{proof}

\begin{rmk}[relation to the ideas in Pelayo, Ratiu and \vungoc~\cite{pvr_carto}] Let $\is$ be a faithful system with compact fibers that admits a~cartographic pair. Corollary~\ref{cor:pull_back_iso_carto} implies that the set of all cartographic images of~$\is$ is an invariant of the isomorphism class of $\is$ that, in some sense, encodes the $\Z$-affine structure on its weakly toric leaf space~$B_{\mathrm{wt}}$. Following Pelayo, Ratiu and \vungoc~\cite[Definition~3.4]{pvr_carto}, the above association can be called a~{\em cartographic invariant} of~$\is$. (Observe that Pelayo, Ratiu and \vungoc~\cite{pvr_carto} work with a~stricter notion of isomorphism than the one given in Definition~\ref{defn:cihs}, see Pelayo, Ratiu and \vungoc~\cite[Definition~1.5]{pvr_carto} and Section~\ref{sec:dscp-noti-isom} below).
 \end{rmk}

One of the aims of this paper is to establish the existence of a certain kind of cartographic pairs for faithful semitoric systems (see Definition~\ref{defn:vat} and Theorem~\ref{prop:rh}), and to describe the collection of all such cartographic pairs (see Theorem~\ref{thm:different_carto}).

\section{Almost-toric systems}\label{sec:at}
Motivated by Symington \cite{symington} and \vungoc\ \cite{vu-ngoc}, this section introduces and studies the fundamental properties of {\em almost-toric systems}, a category of integrable systems generalizing that of weakly toric systems (Definition \ref{defn:toric_sys}) in dimension~4 by allowing for the presence of {\em focus-focus} fibers, which are the Lagrangian analog of nodal fibers in Lefschetz fibrations. To define this category formally, we first recall the notion of almost-toric singular orbits, which are a special class of non-degenerate singular
orbits, focusing on the 4-dimensional case. This is achieved in Section~\ref{sec:singular-orbits}. Section~\ref{sec:at_generalities} defines
almost-toric systems on 4-dimensional manifolds and establishes fundamental properties of leaves and their neighborhoods. Seeing as the faithful semitoric systems of Part~\ref{part:faithful-semitoric} are both faithful and almost-toric, Section~\ref{sec:faithful_at} collects results about systems that satisfy both properties.

\subsection{Singular orbits}\label{sec:singular-orbits}
\subsubsection{Non-degenerate singular orbits in arbitrary dimension}\label{sec:non-degen-sing}
The singular orbits considered in Section \ref{sec:almost-toric-orbits} are {\em non-degenerate} singular orbits, a condition that is briefly recalled below and should be thought of as a~`symplectic' Morse--Bott condition. (For more details regarding non-degenerate orbits, the reader is referred to Bolsinov and Fomenko \cite[Section~1.8]{bol+fom}, \vungoc\ \cite[Section~3.3]{vu_ngoc_book}, and references therein.) Throughout this subsection, let~$\is$ be an integrable system so that $\Phi$ is the moment map of an effective Hamiltonian $\R^n$-action; for any $\mathbf{t} \in \R^n$, denote by $\phi^{\mathbf{t}}\colon \sm \to \sm$ the symplectomorphism induced by acting via $\mathbf{t}$. Moreover, for any $p \in M$,
 denote by $\mathcal{O}_p$ the $\R^n$-orbit through $p$. If~$p$ is singular, then every point in $\mathcal{O}_p$ is singular; thus the
 notion of {\em singular orbit} is well-defined. Next, we introduce the following useful notion.

\begin{Definition}\label{defn:rank} Given an integrable system $\is$, the {\em rank} of a point $p \in M$ is given by $\operatorname{rk} D_p \Phi$.
 \end{Definition}

\begin{rmk}\label{rmk:rank_action} With the above notation, if $p \in M$ is a point of rank $0 \leq k \leq n$, the existence of a Hamiltonian $\R^n$-action implies that the orbit $\mathcal{O}_p$ is a $k$-dimensional immersed, isotropic submanifold of $\sm$ that is diffeomorphic to $\R^{k -c(p)} \times \mathbb{T}^{c(p)}$, where $0 \leq c(p) \leq k$ is called the {\em degree of closedness} of $\mathcal{O}_p$ in Zung \cite[Definition~3.4]{zung-symplectic}. In particular, the rank and the degree of closedness of an orbit are well-defined notions.
 \end{rmk}

Following Bolsinov and Fomenko \cite[Section~1.8.3]{bol+fom}, fix a singular orbit $\mathcal{O} \subset M$ of rank $0 \leq k < n$ and let $p
\in \mathcal{O}$; since, for all $\mathbf{t} \in \R^n$, $\phi^{\mathbf{t}}$ is a symplectomorphism sending $\mathcal{O}$ to itself, it follows that, for all $\mathbf{t} \in \R^n$, $D_p \phi^{\mathbf{t}}$ is a~symplectomorphism of $\big((T_p \mathcal{O}_p)^{\omega}/T_p \mathcal{O}_p,\Omega\big)$, where $(T_p \mathcal{O}_p)^{\omega}$ is the symplectic orthogonal of~$T_p \mathcal{O}_p$ and $\Omega$ is the symplectic form induced by performing linear reduction. Thus we obtain a Lie algebra homomorphism $\R^n \to \operatorname{Sp}\big((T_p\mathcal{O}_p)^{\omega}/T_p \mathcal{O}_p,\Omega\big)$. In fact, this homomorphism only depends on the orbit and not on the choice of point; this is because the action is by an abelian Lie group. Choosing local Darboux coordinates, it is possible to identify $\operatorname{Sp}\big((T_p \mathcal{O}_p)^{\omega}/T_p \mathcal{O}_p, \Omega\big)$ with $\operatorname{Sp}(2(n-k);\R)$; therefore, by taking derivative at the identity, we obtain a Lie algebra homomorphism $\R^n \to \mathfrak{sp}(2(n-k);\R)$ whose image is denoted by~$\mathfrak{h}_{\mathcal{O}}$.

\begin{Definition}\label{defn:non-dege} A singular orbit $\mathcal{O}$ of rank $0 \leq k < n$ is said to be {\em non-degenerate} if $\mathfrak{h}_{\mathcal{O}} \subset \mathfrak{sp}(2(n-k);\R)$ is a Cartan subalgebra.
 \end{Definition}

\begin{rmk}\label{rmk:checking} Since $\mathfrak{sp}(2n;\R)$ is semisimple, its Cartan subalgebras are maximal Abelian and self-normalizing. A~criterion to check that a fixed point in an integrable system with $n$-degrees of freedom is non-degenerate is as follows (cf.\ Bolsinov and Fomenko \cite[Definitions~1.24 and~1.25]{bol+fom}) for details. Let~$p$ be a singular point of rank $0$ in $\is$, where $\Phi=(H_1,\ldots,H_n)$. Then, for all $i=1,\ldots,n$, the Hamiltonian vector field of $H_i$ vanishes at~$p$; thus it makes sense to consider its {\em linearization} at $p$ denoted by $X_{i}^{\mathrm{Lin}}(p) \in \mathfrak{sp}(2n;\R)$. The point $p$ is non-degenerate if $X_{1}^{\mathrm{Lin}}(p),\ldots, X_{n}^{\mathrm{Lin}}(p)$ are linearly independent and if there exists a linear combination $\lambda_1 X_{1}^{\mathrm{Lin}}(p) + \cdots + \lambda_{n}X_{n}^{\mathrm{Lin}}(p)$ with $2n$ distinct, non-zero eigenvalues.
 \end{rmk}

Cartan subalgebras of $\mathfrak{sp}(2(n-k);\R)$ have been classified up to conjugacy in Williamson~\cite{williamson} using the standard isomorphism between $\mathfrak{sp}(2(n-k);\R)$ and $\mathrm{Sym}(2(n-k);\R)$, where the latter is the Lie algebra of symmetric bilinear forms on the linear symplectic vector space $\R^{2(n-k)}$ with Lie bracket given by the Schouten bracket, and the isomorphism sends a quadratic polynomial to its Hamiltonian vector field. The classification of Williamson~\cite{williamson} is recalled below without proof.

\begin{Theorem}\label{thm:will}
 Fix a positive integer $n$ and let $\mathfrak{h} \subset \mathrm{Sym}\left(2n;\R\right)$ be a
 Cartan subalgebra. Then there exist canonical coordinates
 $x_i,y_i$ of the linear symplectic vector space $\R^{2n}$, a triple
 $\left(k_{\mathrm{e}},k_{\mathrm{h}},k_{\mathrm{ff}}\right) \in
 \Z^3_{\geq 0}$ with $k_{\mathrm{e}}+k_{\mathrm{h}}+2
 k_{\mathrm{ff}} = n$, and a basis $H_1,\ldots,H_n$ of
 $\mathfrak{h}$ such that
 \begin{gather*}
 H_i =
 \begin{cases}
 \dfrac{x^2_i + y^2_i}{2} & \text{if } i=1,\ldots,k_{\mathrm{e}}, \\
 x_iy_i & \text{if } i = k_{\mathrm{e}}+1,\ldots, k_{\mathrm{e}} + k_{\mathrm{h}},
 \end{cases}
 \end{gather*}
and, if $i=k_{\mathrm{e}}+k_{\mathrm{h}} + 1, k_{\mathrm{e}}+k_{\mathrm{h}} + 3, \ldots, k_{\mathrm{e}}+k_{\mathrm{h}} + 2j-1, \ldots, k_{\mathrm{e}}+k_{\mathrm{h}} + 2 k_{\mathrm{ff}} - 1$, then
\begin{gather*}
 H_i = x_iy_{i+1}-x_{i+1}y_i, \qquad H_{i+1} = x_iy_i + x_{i+1}y_{i+1}.
\end{gather*}
Moreover, the triple $(k_{\mathrm{e}},k_{\mathrm{h}},k_{\mathrm{ff}})$ determines $\mathfrak{h}$ up to conjugacy.
 \end{Theorem}

\begin{Definition}\label{defn:will_triple} Given a Cartan subalgebra $\mathfrak{h} \subset \mathfrak{sp}(2n;\R)$, the triple $(k_{\mathrm{e}},k_{\mathrm{h}},k_{\mathrm{ff}}) \in \Z^3_{\geq 0}$ with $k_{\mathrm{e}}+k_{\mathrm{h}}+2 k_{\mathrm{ff}} = n$ classifying it up to conjugacy is called the {\em Williamson triple} of $\mathfrak{h}$, where $k_{\mathrm{e}}$, $k_{\mathrm{h}}$, $k_{\mathrm{ff}}$ are referred to as the number of {\em elliptic}, {\em hyperbolic} and {\em focus-focus} components.
 \end{Definition}

Going back to non-degenerate singular orbits, adapting and following Zung \cite[Definition~3.4]{zung-symplectic}, we introduce the following terminology.

\begin{Definition}\label{defn:will_quadr} Given an integrable system $\is$ and a non-degenerate singular orbit $\mathcal{O} \subset M$, its {\em Williamson type} is the element $(k,c, k_{\mathrm{e}},k_{\mathrm{h}},k_{\mathrm{ff}}) \in \Z^5_{\geq 0}$, where~$k$ is the rank of~$\mathcal{O}$, $c$~is its degree of closedness, and $(k_{\mathrm{e}},k_{\mathrm{h}},k_{\mathrm{ff}}) \in \Z^3_{\geq 0}$ is the Williamson triple of~$\mathfrak{h}_{\mathcal{O}}$.
 \end{Definition}

Non-degenerate singular orbits can be {\em linearized} (cf., for instance, Dufour, Molino, Eliasson, Miranda and Zung \cite{duf_mol,eliasson-thesis,miranda-zung} amongst others). While the various linearization results are beyond the scope of this article, it is worthwhile observing that, in the absence of hyperbolic blocks, i.e., if $k_{\mathrm{h}} = 0$, and when the orbits are compact, i.e., if $k = c$, the linearization result is stronger: there exist canonical coordinates which also put the moment map in standard form. This is part of the motivation for introducing the singular orbits studied in Section~\ref{sec:almost-toric-orbits}. To conclude this subsection, we state the following characterization of non-degenerate, compact singular orbits of {\em purely elliptic type}, i.e., whose Williamson types are given by elements of the form $(k,k,n-k,0,0)$, relating them to singular weakly toric leaves (cf.\ Dufour, Molino and Eliasson \cite{duf_mol,eliasson-thesis}). Such orbits are henceforth referred to as {\em elliptic tori}.

\begin{Theorem}\label{thm:comp_ell} Let $\mathcal{O}$ be an elliptic torus in an integrable system $\is$. Then there exists a~$($connected$)$ open neighborhood $V \subset M$ of~$\mathcal{O}$ whose corresponding subsystem is weakly toric $($see Definition~{\rm \ref{defn:toric_sys})}. In particular, $\mathcal{O}$ is a singular weakly toric leaf.
 \end{Theorem}

Theorem \ref{thm:comp_ell} is the crucial ingredient in proving that elliptic tori can be linearized (cf.\ Dufour, Molino and Eliasson \cite{duf_mol,eliasson-thesis}).

Integrable systems with compact fibers all of whose singular orbits are elliptic tori play an important role in the classification of integrable systems with compact fibers.

\begin{Definition}\label{defn:cihs_elliptic} An integrable system with compact fibers is said to have {\em elliptic singularities} if all its singular orbits are of purely elliptic type, i.e., are elliptic tori.
 \end{Definition}

Theorem \ref{thm:comp_ell} readily yields the following result, showing that integrable systems with elliptic singularities generalize (weakly) toric systems.

\begin{Corollary}\label{cor:elliptic_singularities} An integrable system with compact fibers has elliptic singularities if and only if all its leaves are weakly toric.
 \end{Corollary}

Integrable systems with elliptic singularities have been classified by Boucetta and Molino (cf.~\cite{bm}), generalizing a~construction due to Duistermaat, Dazord and Delzant (cf.~\cite{dd, duistermaat}).

\subsubsection{Almost-toric orbits}\label{sec:almost-toric-orbits}
Motivated by the work of Symington and \vungoc\ \cite{symington,vu-ngoc}, we distinguish the following family of singular orbits.

\begin{Definition}\label{defn:almost_toric} A singular orbit $\mathcal{O}$ in an integrable system $\is$ is said to be {\em almost-toric} if it is compact and non-degenerate without hyperbolic blocks.
\end{Definition}

While Definition \ref{defn:almost_toric} makes sense for integrable systems on symplectic manifolds of any dimension (cf.\ Izosimov~\cite{izo_at} for results in this direction), we are concerned solely with almost-toric orbits in integrable systems with two degrees of freedom. Therefore, throughout the rest of this section, $\is$ is an integrable system on a 4-dimensional symplectic manifold unless otherwise stated. If $\mathcal{O}$ is an almost-toric orbit of $\is$, its Williamson type (see Definition~\ref{defn:will_quadr}) is constrained to be of one of three types, namely
 \begin{itemize}[leftmargin=*]\itemsep=0pt
 \item {\em elliptic-elliptic}, given by $(0,0,2,0,0)$,
 \item {\em elliptic-regular}, given by $(1,1,1,0,0)$,
 \item {\em focus-focus}, given by $(0,0,0,0,1)$.
 \end{itemize}
 The first two are elliptic tori of dimension 0 and 1 respectively. On the other hand, focus-focus points are completely characterized by the following local normal form (cf.\ Chaperon, Eliasson, \vungoc\ and Wacheux \cite{chaperon,eliasson-thesis,vu_ngoc_wacheux}):

\begin{Theorem}\label{thm:ff_normal} Let $p$ be a focus-focus point in the integrable system $\is$ with two degrees of freedom and consider $\big(\R^4,\omega_0\big)$, where $\omega_0 = dx_1\wedge dy_1 + dx_2\wedge dy_2$. There exist open neighborhoods $V \subset M$ and $W \subset \R^4$ of~$p$ and the origin respectively, such that $(V,\omega|_V,\Phi|_V)$ is isomorphic to $(W,\omega_0|_W,\mathbf{q}|_W)$ via a pair $(\Psi,\psi)$ with $\psi(p) = 0$, where $\mathbf{q} = (q_1,q_2)$, $q_1 = x_1y_2 - x_2y_1$ and $q_2 = x_1y_1 + x_2y_2$.
\end{Theorem}

\begin{rmk}\label{rmk:ff_basic} In Theorem \ref{thm:ff_normal}, observe that the flow of the Hamiltonian vector field~$X_{q_1}$ is periodic, and that~$W$ can be chosen to be saturated with respect to the effective Hamiltonian $S^1$-action whose moment map is given by~$q_1$.
 \end{rmk}

An immediate consequence of Theorem \ref{thm:ff_normal} is the following result, stated below without proof.

\begin{Corollary}\label{rmk:discrete} In an integrable system with two degrees of freedom, the set of focus-focus points is discrete.
\end{Corollary}

Elliptic tori and focus-focus points differ significantly. An immediate topological difference is that the former are leaves of the system while for the latter are not. A~crucial geometric difference lies in the fact that the former support a local effective Hamiltonian $\mathbb{T}^2$-action whose moment map has components that Poisson commute with the integrals of the system (we say that the action is {\em system-preserving}), while the latter possesses only a unique (up to sign) system-preserving effective Hamiltonian $S^1$-action
(cf.\ Zung \cite[Proposition~4]{zung-focus-focus}).

\begin{figure}[h] \centering
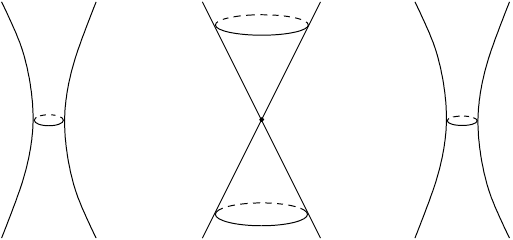
\caption{Fibers in an open neighborhood of a focus-focus point. Regular fibers are diffeomorphic to cylinders, while the singular fiber is given by the union of two Lagrangian planes transversally intersecting at the focus-focus point.} \label{LABEL4}
\end{figure}

\begin{rmk}\label{rmk:geom_prop} Let $p$ be a focus-focus point in an integrable system $\is$ and let~$V$ be an open neighborhood of~$p$ that can be put in local normal form. Then
 \begin{itemize}[leftmargin=*]\itemsep=0pt
 \item the restriction $\Phi|_V$ is open;
 \item $p$ is the only singular point of $\Phi|_V$;
 \item the fibers of $\Phi|_V$ are connected;
 \item a fiber of $\Phi|_V$ is either diffeomorphic to a cylinder (if it does not contain $p$) or given by the union of two Lagrangian planes intersecting transversally at~$p$. In particular, if the latter fiber is denoted by $L$, then $(V \cap L) \smallsetminus \{p\}$ consists of two connected components, each diffeomorphic to a cylinder (see Fig.~\ref{LABEL4});
 \item there exist smooth sections $\sigma_1$, $\sigma_2$ defined near $\Phi(p)$ whose image lies inside~$V$ with the property that $\sigma_1(\Phi(p)), \sigma_2(\Phi(p))$ lie in distinct connected components of $(V \cap L) \smallsetminus \{p\}$, where~$L$ denotes the leaf through~$p$.
 \end{itemize}
\end{rmk}

\subsection{Definition and fundamental properties}\label{sec:at_generalities}

Following Symington \cite{symington} and \vungoc\ \cite{vu-ngoc}, we introduce a~category of integrable systems of two degrees of freedom that generalize weakly toric systems on 4-dimensional manifolds while retaining significant similarities.

\begin{Definition}\label{defn:almost-toric} An integrable system $\is$ on a 4-dimensional symplectic manifold is {\em almost-toric} if $\Phi$ is proper onto its image and all of its singular orbits are almost-toric.
\end{Definition}

Throughout this paper, we consider almost-toric systems only on 4-dimensional symplectic manifolds. For this reason and for the sake of brevity, henceforth we drop the dependence on dimension when referring to an almost-toric system. The following result can be used to describe almost-toric systems in an equivalent fashion; seeing as it follows directly from Zung~\cite[Proposition~3.5]{zung-symplectic}, its proof is omitted.

\begin{Lemma}\label{lemma:compact} Let $\is$ be an integrable system on a 4-dimensional symplectic manifold with compact fibers, all of whose singular orbits are non-degenerate without hyperbolic blocks. Then the singular orbits are compact and, in particular, almost-toric.
\end{Lemma}

The above notion of almost-toric system differs slightly from that in \vungoc\ \cite{vu-ngoc}, for Definition~\ref{defn:almost-toric} only requires that $\Phi$ be proper {\em onto its image}, as opposed to being proper.

\begin{rmk}\label{rmk:at_full_cat} Almost-toric systems form a full subcategory of~$\mathcal{IS}(2)$, henceforth referred to as the {\em
 category of almost-toric systems} and denoted by~$\mathcal{AT}$.
\end{rmk}

\begin{rmk}\label{rmk:sat_subsyst_at} A saturated subsystem of an almost-toric system is almost-toric.
\end{rmk}

The restriction on the types of singular orbits in an almost-toric system $\is$ implies the singular leaves of $\Phi$ are either elliptic tori or contain at least one focus-focus singular orbit. That dichotomy arises because the local normal form for elliptic singular orbits implies such orbits (which are elliptic tori) make up whole connected components of a fiber of~$\Phi$. Henceforth, leaves that contain focus-focus orbits are referred to as {\em focus-focus leaves}. Denote by $\mathcal{L}_{\mathrm{ff}}$ the set of points in the leaf space $\mathcal L$ corresponding to leaves containing focus-focus singular orbits. Then $\mathcal{L}_{\mathrm{sing}} = \mathcal{L}_{\mathrm{e}} \cup \mathcal{L}_{\mathrm{ff}}$, where $\mathcal{L}_{\mathrm{e}}$ is the elliptic part of $\mathcal{L}$. Elements of $\mathcal{L}_{\mathrm{ff}}$ are called {\em focus-focus values $($in the leaf space$)$}.

\begin{Definition}\label{defn:multiplicity} Let $c \in \mathcal{L}_{\mathrm{ff}}$ be a focus-focus value in the leaf space of an almost-toric system $\is$. The {\em multiplicity} of $c$, denoted $r_c \geq 1$, is the number of focus-focus singular orbits in the corresponding leaf of~$\Phi$.
\end{Definition}

Focus-focus values, counted with multiplicity, are an invariant of the isomorphism class of an almost-toric system. Henceforth, all focus-focus values are counted with multiplicity unless otherwise stated.

The topology of focus-focus leaves is well-known and is completely determined by the finite number $r \geq 1$ of focus-focus singular points contained in the corresponding leaf. In particular, each focus-focus leaf is homeomorphic to a torus with~$r$ homologous cycles, each collapsed to a point (cf.\ Bolsinov and Fomenko \cite[Chapter~9.8]{bol+fom}). Lying at the heart of this result is the existence of a~vector field tangent to a focus-focus leaf, whose flow is periodic and whose fixed points are precisely the focus-focus points.

\begin{figure}[h] \centering
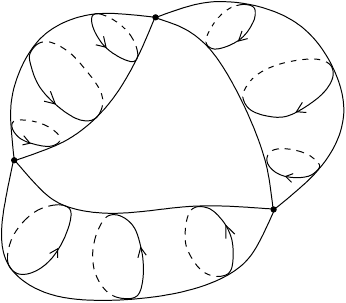
 \caption{A focus-focus leaf containing 3 singular points and the flow of the periodic vector field on it.} \label{LABEL5}
\end{figure}

In fact, the $S^1$-action on any focus-focus leaf can be extended to a suitable neighborhood of the leaf, one that is saturated with respect to the quotient map to the leaf space. To prove this, we start by establishing the existence of this suitable neighborhood.

\begin{Proposition}\label{prop:ff_leaf_sat_nbhd}
 Given an almost-toric system $\is$, any focus-focus leaf $L \subset M$ admits an open neighborhood $V$ satisfying the following properties:
 \begin{itemize}[leftmargin=*]\itemsep=0pt
 \item $L$ is the only singular leaf in $V$;
 \item $V$ is saturated with respect to the quotient map $q \colon M \to \mathcal{L}$;
 \item $V$ contains at most one leaf of any fiber of $\Phi$.
 \end{itemize}
\end{Proposition}

(Proposition \ref{prop:ff_leaf_sat_nbhd} is probably well-known to experts but we could not find a proof in the literature. For this reason, a~proof is provided below.)

\begin{proof} Fix a focus-focus leaf $L$. First, we show that $L$ admits an open neighborhood $Z$ that is saturated with respect to $q$ and in which $L$ is the only focus-focus leaf. Suppose not; then there exists a~sequence of focus-focus points $\{p_n\}$ with the property that $\Phi(p_n) \to \Phi(L)$; since~$\Phi$ is proper onto its image, there exists a convergent subsequence $p_{n_j} \to p$. The limit point~$p$ is necessarily singular, but the local normal form for almost-toric orbits yields a contradiction.

Fix such a neighborhood $Z$ and let $p_1,\ldots,p_r \in L$ denote the focus-focus points in $L$. For $i=1,\ldots, r$, let $V_i \subset Z$ be an open neighborhood of $p_i$ that can be put in local normal form (see Section~\ref{sec:almost-toric-orbits}). Consider the subset
 \begin{gather*}
 \hat{V} := \bigcup\limits_{i=1}^r \bigcup\limits_{\mathbf{t} \in \R^2} \phi^{\mathbf{t}} (V_i ),
 \end{gather*}
where, as in Section \ref{sec:almost-toric-orbits}, $\phi^{\mathbf{t}}$ denotes the Hamiltonian action by $\mathbf{t} \in \R^2$; this is the union of the orbits that intersect at least one~$V_i$. Since, for $i=1,\dots, r$, $V_i$ is open, so is~$\hat{V}$; moreover, $\hat{V}$~contains~$L$ because $L = \bigcup\limits_{i=1}^r \bigcup\limits_{\mathbf{t} \in \R^2} \phi^{\mathbf{t}}(V_i \cap L)$. Next we show that $\hat{V}$ is also saturated with respect to~$q$. To see this, observe that if $p \in V_i \smallsetminus L$, then the leaf passing through $p$ is contained in~$Z$ and is not equal to $L$, thus implying that it is not a~focus-focus leaf. Since~$p$ is regular (by the local normal form for focus-focus points), the leaf through~$p$ is regular and is, therefore, an orbit of the Hamiltonian $\R^2$-action, which is contained in $\hat{V}$ by construction.

The above construction does not necessarily guarantee that $\hat{V}$ contains at most one leaf of any fiber of $\Phi$. However, the local normal form for a focus-focus point and the~$\R^2$ action can be used to determine a possibly smaller neighborhood in which that property holds. Specifically, fix some $i\in \{1,\ldots, r\}$. There exists a~smooth section $\sigma_i$ of $\Phi$ defined near $\Phi(p_i)$ whose image is contained in $V_i$; this implies that $\sigma_i(\Phi(p_i)) \in (V_i \cap L) \smallsetminus \{p_i\}$ (see Remark~\ref{rmk:geom_prop}). The structure of the focus-focus leaf $L$ implies that there exists~$\mathbf{t}_0 \in \R^2$ with the following property
 \begin{itemize}[leftmargin=*]\itemsep=0pt
 \item if $r=1$, then $\sigma_i(\Phi(p_i))$ and $\phi^{\mathbf{t}_0} (\sigma_i(\Phi(p_i)) )$ lie in different connected components of $(V_i \cap L) \smallsetminus \{p_i\}$ (cf.\ \vungoc~\cite{vu_ngoc-bohr}~-- see Fig.~\ref{LABEL6});
 \item if $r > 1$, there exists $j \neq i$ with $\phi^{\mathbf{t}_0} (\sigma_i(\Phi(p_i))) \in V_j$ (cf.\ Bolsinov and Fomenko \cite[Chapter~9.8]{bol+fom}~-- see Fig.~\ref{LABEL7}).
 \end{itemize}
In other words, the section $\sigma_i$ flows out of $V_i$ and into $V_j$ (and if $i=j$, then it approaches $V_i$ from the `opposite' side).

\begin{figure}[h] \centering
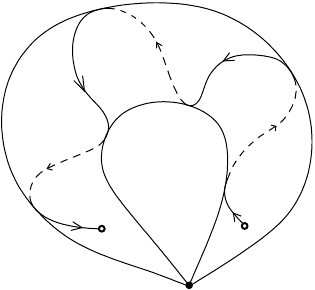
\caption{The flow line of the regular point $\sigma_i(\Phi(p_i))$ for time $\mathbf{t}_0$ in the case $r =1$.} \label{LABEL6}
\end{figure}

\begin{figure}[h] \centering
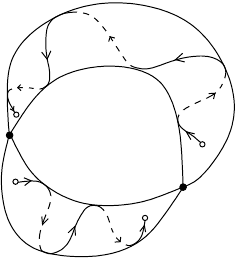
\caption{Flow lines of regular points on a focus-focus leaf containing 2 focus-focus points, displaying heteroclinic behavior and spiraling from one focus-focus point to the other.} \label{LABEL7}
\end{figure}

Let $c \in \R^2$ be sufficiently close to $\Phi(p_i)$ so that $\sigma_i(c)$ is defined; then the above properties show that
 \begin{itemize}[leftmargin=*]\itemsep=0pt
 \item if $r=1$, $\Phi^{-1}(c) \cap \hat{V}$ is connected;
 \item if $r > 1$, the intersections $\Phi^{-1}(c) \cap V_i$ and $\Phi^{-1}(c) \cap V_j$ lie on the same leaf of $\Phi^{-1}(c)$.
 \end{itemize}
In the latter case, using again the structure of the focus-focus leaf~$L$ (cf.\ Bolsinov and Fomenko \cite[Chapter~9.8]{bol+fom}), we can iterate the above argument finitely many times to ensure that, for all~$c$ sufficiently close to $\Phi(L)$ ($ = \Phi(p_i)$ for all $i=1,\ldots,r$), $\Phi^{-1}(c) \cap \hat{V}$ is connected. This shows that~$\hat{V}$ can be shrunk as desired.
\end{proof}

\begin{Corollary}\label{cor:ff_values_discrete} The set of focus-focus values in the leaf space of an almost-toric system is discrete.
\end{Corollary}

Any open neighborhood of a focus-focus leaf as in Proposition \ref{prop:ff_leaf_sat_nbhd} is henceforth referred to as a {\em {\rm ($q$-)}saturated regular neighborhood} of a focus-focus leaf. A saturated regular neighborhood has the necessary $S^1$-symmetry.

\begin{Proposition}\label{prop:S1-action} Given a focus-focus leaf $L$ of an almost-toric system $\is$, any saturated regular neighborhood of $L$ admits a local effective system-preserving Hamiltonian $S^1$-action. Moreover, this action is unique up to sign.
\end{Proposition}

\begin{proof}[Sketch of proof] The ideas behind proving this result are known (cf.\ Bolsinov, Fomenko~\cite[Lemma~9.8]{bol+fom} and Zung \cite[Section 3]{zung-focus-focus}), but the key ideas are provided below for completeness. Let $L$ be a focus-focus leaf and let~$V$ be a~saturated regular neighborhood of~$L$. Let~$p$ be a~focus-focus point on $L$; the local normal form for~$p$ implies that, near $p$, there exists a local system-preserving Hamiltonian $S^1$-action. By construction of~$V$, this action can be extended to the whole of $V$ and is independent (up to sign) of the choice of focus-focus point $p \in L$. Uniqueness of the action up to sign is proved in \cite[Proposition~4]{zung-focus-focus}.
\end{proof}

In fact, a saturated regular neighborhood of a focus-focus leaf (of multiplicity one) is a singular Liouville foliation of (simple) focus-focus type in the sense of \vungoc\ \cite[Definition~2.4]{vu_ngoc_invariant}. Moreover, saturated regular neighborhoods of focus-focus leaves ought to be thought of as analogous to the neighborhoods of elliptic tori that can be put in local normal form. For instance, the following result holds.

\begin{Proposition}\label{prop:faithful_at} Let $V$ be a saturated regular neighborhood of a focus-focus leaf~$L$ in an almost-toric system~$\is$. The subsystem $(V, \omega|_V, \Phi|_V)$ is faithful and almost-toric.
\end{Proposition}
\begin{proof} By construction, the fibers of $\Phi|_V$ are connected and all singular orbits in the subsystem are almost-toric. To prove the result, it suffices to show that $\Phi|_V$ is proper onto its image, for then the subsystem is almost-toric by Definition~\ref{defn:almost-toric} and faithful by Lemma~\ref{lemma:faithful_general}. Set $B_V:= \Phi|_V(V)$; to prove that~$\Phi|_V$ is proper onto its image, it suffices to check that it is proper at every point $c \in B_{V}$. Observe that, by definition of~$V$, $B_V$ contains only one singular value, which equals~$\Phi(L)$. If $c \neq \Phi(L)$, then $\Phi|_V$ is proper at $c$ as $\Phi|_{V \smallsetminus L}$ is a submersion with compact and connected fibers. If $c = \Phi(L)$, then arguing as in the second half of the proof of Lemma~\ref{lemma:faithful_general}, it can be shown that $\Phi|_V$ is also proper at $c$.
\end{proof}

To summarize the above results and motivate subsequent sections, we state the following description of regular neighborhoods of leaves in almost-toric systems.

\begin{Corollary}\label{cor:reg_neigh} Given an almost-toric system $\is$, any leaf $L$ admits an open neighborhood $V$ satisfying the following properties:
 \begin{itemize}[leftmargin=*]\itemsep=0pt
 \item $V$ is saturated with respect to the quotient map $q\colon M \to \mathcal{L}$;
 \item the subsystem $(V, \omega|_V, \Phi|_V)$ is faithful and admits a system-preserving Hamiltonian $S^1$-action.
 \end{itemize}
\end{Corollary}

\subsection{Faithful almost-toric systems}\label{sec:faithful_at}
By Corollary \ref{cor:reg_neigh}, any leaf in an almost-toric system admits an open neighborhood whose corresponding subsystem is faithful almost-toric. Therefore it makes sense to view faithful almost-toric systems as building blocks of almost-toric systems. Throughout this subsection, let $\is$ be a faithful almost-toric system. Moreover, fix the identification between $\mathcal{L}$ and $B=\Phi(M)$ and denote by $B_{\mathrm{ff}}$ the image of~$\mathcal{L}_{\mathrm{ff}}$ under this identification.

\begin{Lemma}\label{rk:ff_discrete} Given a faithful almost-toric system, the set of focus-focus values satisfies $B_{\mathrm{ff}} \subset \operatorname{Int}(B)$ and is discrete in $\operatorname{Int}(B)$.
\end{Lemma}

\begin{proof} Theorem \ref{thm:ff_normal}, the local normal form for focus-focus points, yields readily that $B_{\mathrm{ff}} \subset \operatorname{Int}(B)$. The last statement then follows from Corollary \ref{cor:ff_values_discrete}.
 \end{proof}

The arguments in \vungoc\ \cite[Proposition 3.9]{vu-ngoc} can be used to prove the following result, stated without proof.
\begin{Lemma}\label{lem:breg_connected} Given a faithful almost-toric system on a connected symplectic manifold, the subsets $B_{\mathrm{reg}}$ and $\operatorname{Int}(B)$ of its moment map image $B$ are path-connected.
\end{Lemma}

The restriction on the types of singular orbits, as well as faithfulness, imply the following useful fact for faithful almost-toric systems.

\begin{Lemma}\label{lemma:pc} Given a faithful almost-toric system $\is$ and a continuous path $\gamma \colon [0,1]$ $\to B$, the subset $\Phi^{-1}(\gamma([0,1]))$ is path-connected.
\end{Lemma}

\begin{proof} Fix a path $\gamma$ as in the statement. By the local normal form for almost-toric singular points, given any point $c \in B$, there exists an open, path-connected neighborhood $U \subset B$ of $c$ and a continuous section $\sigma \colon U \to M$. The image of $\gamma$ is contained in the union of finitely many such neighborhoods. Connectedness of the fibers of~$\Phi$ then implies the desired result.
\end{proof}

Lemma \ref{lemma:pc} has the following simple consequence that states that path-connectedness of the fibers of one of the components of the moment map of a faithful almost-toric system can be checked directly on the moment map image.

\begin{Lemma}\label{lemma:connected_J} Let $\vat$ be a faithful almost-toric system. Then the fibers of $J$ are path-connected if and only if, for all $x_0 \in \R$, the intersection $B \cap \{(x,y) \,|\, x = x_0\}$ is either empty or path-connected.
\end{Lemma}

\begin{proof} If the fibers of $J$ are path-connected, then the claimed result holds. Conversely, suppose that $B \cap \{(x,y) \,|\, x = x_0\}\neq \varnothing $ is path-connected; the aim is to prove that $J^{-1}(x_0)$ is path-connected. Let $p_1,p_2 \in J^{-1}(x_0)$ and set $c_i = \Phi(p_i)$, for $i=1,2$. Then the vertical segment joining $c_1$ to $c_2$ is contained in $B$ by assumption; Lemma~\ref{lemma:pc} then implies that there exists a~path joining $p_1$ to $p_2$ contained in the preimage under $\Phi$ of that vertical segment. Since that preimage is contained in $J^{-1}(x_0)$, this shows the desired result.
\end{proof}

The following result describes some {\em smooth} properties of the moment map image of a faithful almost-toric system.

\begin{Lemma}\label{rmk:almost-toric-image} Given a faithful almost-toric system $\is$,
 \begin{itemize}[leftmargin=*]\itemsep=0pt
 \item the structure of smooth manifold with corners on $B_{\mathrm{wt}}$, the weakly toric leaf space, extends to all of $B$;
 \item the image of elliptic tori is precisely $\partial B = \partial B_{\mathrm{wt}}$, where corners and facets $($or curved edges$)$ of~$B$ are the images of elliptic-elliptic and elliptic-regular points, respectively;
 \item the set of focus-focus values $B_{\mathrm{ff}} \subset \operatorname{Int}(B)$ is at most countable, and the set of its limit points in~$\R^2$ is contained in $\mathrm{Bdy}(B) \smallsetminus \partial B$.
 \end{itemize}
\end{Lemma}

\begin{proof} Since $B = B_{\mathrm{wt}} \cup B_{\mathrm{ff}}$ and $B_{\mathrm{wt}}$ has a structure of smooth manifold with corners (see Proposition~\ref{cor:lt_iam}), it suffices to define charts near each point in $B_{\mathrm{ff}}$ that are compatible with the given smooth atlas on $B_{\mathrm{wt}}$. Observe that the smooth atlas on $B_{\mathrm{wt}}$ is compatible with the standard smooth structure of the ambient~$\R^2$, i.e., charts are smooth as maps between subsets of $\R^2$ in the standard sense. This is a~consequence of faithfulness, of the Liouville--Arnol'd theorem (see Theorem~\ref{thm:la}), and of the local normal form for elliptic tori which, in light of Theorem~\ref{thm:comp_ell}, is given (up to isomorphism of integrable systems) by Theorem~\ref{thm:gsm}. Let $c \in B_{\mathrm{ff}}$. By Lemma~\ref{rk:ff_discrete}, $c \in \operatorname{Int}(B)$ and, thus, there exists an open neighborhood $U$ of $c$ in $\R^2$ that is contained in $\operatorname{Int}(B) \subset B$. On this open neighborhood, we define the smooth structure to be the one induced by inclusion into $\R^2$, i.e., the local chart is given by the inclusion $U \hookrightarrow \R^2$. By the above argument, this chart is compatible with the smooth charts on $B_{\mathrm{wt}}$, thus proving the first point. To prove the second bullet point, observe that $B_{\mathrm{ff}} \subset B \smallsetminus \partial B$ so that $\partial B = \partial B_{\mathrm{wt}}$, i.e., $B$ and $B_{\mathrm{wt}}$ have equal boundaries as smooth manifolds with corner. The result then follows from the description of the smooth boundary of $B_{\mathrm{wt}}$ (see Remark~\ref{rmk:loc_tor_faithful}). It remains to prove the last bullet point. Observe that $B$ is second countable and that any discrete subset of a second countable topological space is at most countable. Therefore $B_{\mathrm{ff}}$ is at most countable. Clearly, the set of limit points of $B_{\mathrm{ff}}$ is contained in $\mathrm{Bdy}(B)$. However, by faithfulness and the local normal form for elliptic tori, no point in $\partial B = \partial B_{\mathrm{wt}}$ can be a~limit point of~$B_{\mathrm{ff}}$ as desired.
 \end{proof}

\begin{exm}[{This example is based on Symington \cite[Section~11]{symington} and Zung \cite[Example~4.19]{zung_ii}}]\label{exm:K3} Suppose~$(M,\omega)$ is symplectomorphic to a K3 surface (for example, a smooth quartic hypersurface in~$\C P^3$). Such a symplectic manifold admits singular Lagrangian fibrations over $S^2$ in which each singular fiber has a neighborhood that, with respect to an appropriate coordinate chart on $S^2$, defines a faithful almost-toric system that is a regular saturated neighborhood of a focus-focus fiber with one singular orbit. Suppose $\Pi\colon (M,\omega)\to S^2$ is such a fibration and let $p\in S^2$ be the image of a~regular fiber. Let $N=M\smallsetminus\Pi^{-1}(p)$ and let $\phi:S^2\smallsetminus p\to \R^2$ be an embedding. Then $(N,\omega|_N,\phi\circ\Pi|_N)$ defines a faithful almost-toric system with 24 focus-focus leaves.
\end{exm}

A natural question arising from Lemma \ref{rmk:almost-toric-image} is whether the $\Z$-affine structure $\mathcal{A}_{\mathrm{wt}}$ on~$B_{\mathrm{wt}}$ can be extended to~$B$. The presence of focus-focus fibers prevents this from happening, as the $\Z$-affine structure on any neighborhood of a focus-focus value has non-trivial affine holonomy.
\begin{Theorem}[{Zung \cite[Proposition~3 and Corollary~1]{zung-focus-focus}}]\label{thm:holonomy} Let $\is$ be faithful almost-toric system and let $U \subset B = \Phi(M)$ be an open neighborhood of a~focus-focus value~$c$, sufficiently small such that~$U$ contains no other focus-focus value. The affine holonomy of the $\Z$-affine structure on $U \smallsetminus \{ c\} \subset B_{\mathrm{wt}}$ is given, in a suitable basis, by
 \begin{gather*} 
 \begin{split}
 \pi_1(U) \cong \Z &\to \mathrm{AGL}(2;\Z), \\
 k &\mapsto \left(
 \begin{pmatrix}
 1 & 0 \\
 kr_c & 1
 \end{pmatrix},
 \begin{pmatrix}
 0 \\
 0
 \end{pmatrix}
 \right),
 \end{split}
 \end{gather*}
where $r_c \geq 1$ is the multiplicity of $c$.
\end{Theorem}

\begin{rmk}\label{rmk:eigenline} The eigenspace associated to the above representation reflects the uniqueness (up to sign) of the local effective system-preserving Hamiltonian $S^1$-action in a~neighborhood of a singular fiber containing focus-focus points (see Proposition~\ref{prop:S1-action}). With respect to the local choice of basis in Theorem~\ref{thm:holonomy}, this action is induced by the first integral of the moment map.
\end{rmk}

Another natural question is which almost-toric systems are isomorphic to toric ones? Certainly, the system must not contain focus-focus points.
If the system is faithful, it suffices that there be a $\Z$-affine immersion of the moment map image into $(\R^2,\mathcal{A}_0)$. However, as the next example illustrates, an absence of focus-focus points does not suffice.

\begin{exm}\label{exm:no_carto} The integrable system with elliptic singularities whose total space is diffeomorphic to $S^2 \times \mathbb{T}^2$ constructed in Example~\ref{rmk:not_weaker} is faithful almost-toric but is not isomorphic to a~toric system, as it has a moment map image that is not simply connected (it is the annulus $A$ considered in Example~\ref{rmk:not_weaker}), and hence not homeomorphic to a polygon.
\end{exm}

\begin{rmk}\label{rmk:carto_image} Consider a faithful almost-toric system $\is$ with focus-focus points and suppose that $(f,S)$ is a~cartographic pair. The definition of a cartographic pair (Definition~\ref{defn:rh}) and the local normal form for singular orbits of elliptic type together imply that the cartographic image of curved edges in~$S$ are line segments whose tangent vectors can be chosen to have coprime integer coefficients, and whenever two edges are incident to a corner contained in~$S$ the associated tangent vectors span~$\Z^2 \subset \R^2$.
\end{rmk}

An important question is thus, when does a faithful almost-toric system admit a cartographic homeomorphism? Addressing this problem in full generality is beyond the scope of this paper. However, in light of Corollary~\ref{cor:reg_neigh}, a natural family of faithful almost-toric systems to consider arises: namely, those admitting a~{\em global} system-preserving $S^1$-action.

\newpage

\part{Faithful semitoric systems}\label{part:faithful-semitoric}
This part studies faithful semitoric systems and proves the main results of the paper. As pointed out in the introduction, faithful semitoric systems are closely related to (proper) semitoric systems (cf.\ Pelayo, Ratiu and \vungoc~\cite{pvr_carto,vu-ngoc}). Many of the ideas and proofs
that appear in this section are inspired by the work in {\em op.\ cit.} Section~\ref{sec:defin-basic-prop} introduces faithful semitoric systems and establishes their basic properties, while Section~\ref{sec:dscc-home-boldsymb} establishes that faithful semitoric systems possess cartographic homeomorphisms of a~special type (see Theorem~\ref{prop:rh}). That result allows one to prove that any faithful semitoric system is isomorphic to an {\em $\boldsymbol{\eta}$-cartographic} one (see Definition~\ref{defn:eta-carto} and Theorem~\ref{thm:eta-carto}). These representatives are particularly useful when defining surgeries on faithful semitoric systems (cf.\ the forthcoming~\cite{HSSS-surgeries}).

\section{Basic properties of faithful semitoric systems}\label{sec:defin-basic-prop}
This section introduces, and establishes basic properties of, faithful semitoric systems and proves that such systems can be viewed as `building
blocks' of almost-toric systems (see Proposition~\ref{prop:at-vat} for a~formal statement). Unlike the more general almost-toric systems, faithful semitoric systems possess a global Hamiltonian $S^1$-action, a fact that has some important geometric consequences that are investigated in Section~\ref{sec:geom-impl-s1}. In Section~\ref{sec:dscex-relat-other}, the relation between faithful semitoric systems and (proper) semitoric systems is explained, while Section~\ref{sec:dscp-noti-isom} identifies various notions of isomorphisms for faithful semitoric systems, relating those of\ Pelayo, Ratiu and \vungoc~\cite[Definition~1.5]{pvr_carto} with the notion of isomorphism of Definition~\ref{defn:cihs} (see Remark~\ref{rmk:toric_im} and Proposition~\ref{prop:im_equ}).

\subsection{Definition and a connectedness result}\label{sec:dscd-conn-result}
To the best of our knowledge, there are no general results regarding the existence of cartographic homeomorphisms for faithful almost-toric systems, even if the total space is closed (cf.\ Leung, Symington~\cite{leung_symington} and Symington~\cite{symington}). However, the existence results of Pelayo, Ratiu, \vungoc~\cite{pvr_carto} and \vungoc~\cite{vu-ngoc} hint at the fact that if the first integral in an almost-toric system $\vat$ is the moment map for an effective Hamiltonian $S^1$-action, then some control on $J$ suffices.

\begin{Definition}\label{defn:vat} A {\em faithful semitoric system} is an integrable system $\vat$ on a~connected 4-dimensional symplectic manifold satisfying:
\begin{enumerate}[label=(F\arabic*), ref=(F\arabic*), leftmargin=*]\itemsep=0pt
 \item \label{item:23} the moment map $\Phi$ is proper onto its image;
 \item \label{item:24} all singular orbits of $\Phi$ are almost-toric, i.e., they are compact, non-degenerate and without hyperbolic blocks (see Definition \ref{defn:almost_toric});
 \item \label{item:3} the first integral $J$ is the moment map of an effective Hamiltonian $S^1$-action;
 \item \label{item:4} any fiber of $J$ is connected;
 \item \label{item:100} the set of critical values of $J$ does not contain any limit points in $J(M)$; 	
 \item \label{item:101} any fiber of $J$ contains at most finitely many isolated fixed points of the $S^1$-action.
 \end{enumerate}
\end{Definition}

\begin{rmk}\label{rmk:defn_fst} It is worth observing that properties~\ref{item:23} and~\ref{item:24} imply that the systems of Definition~\ref{defn:vat} are almost-toric in the sense of Definition~\ref{defn:almost-toric}, while property \ref{item:3} yields the existence of a global system-preserving Hamiltonian $S^1$-action. Properties~\ref{item:4}--\ref{item:101} control that $S^1$-action, imposing conditions that are automatic in the case in which $J$ is a proper map (cf.\ \vungoc~\cite{vu-ngoc}). In particular, property~\ref{item:4} is equivalent
 to path-connectedness of the fibers of $J$ by the local normal form for a Hamiltonian $S^1$-action (cf.\ Guillemin, Sternberg and Marle~\cite{GS2,marle}).
\end{rmk}

Properties \ref{item:23}--\ref{item:3} justify using the terminology `semitoric' for the systems of Definition~\ref{defn:vat} (extending the original terminology introduced by \vungoc\ in~\cite{vu-ngoc}). To justify the use of the adjective `faithful' in Definition~\ref{defn:vat}, we state the following connectedness result.

\begin{Theorem}\label{thm:connected} Let $\vat$ be a faithful semitoric system. Then the fibers of $\Phi$ are connected.
\end{Theorem}

Assuming Theorem \ref{thm:connected}, Lemma \ref{lemma:faithful_general} yields that the systems of Definition~\ref{defn:vat} are faithful in the sense of Definition~\ref{defn:faithful}. The proof of Theorem \ref{thm:connected} relies heavily on arguments of Pelayo, Ratiu and \vungoc\
 underlying \cite[proof of Theorem 4.7]{pvr}. It is included below for completeness and to highlight the fact that, in our context, it suffices to assume that the moment map is proper {\em onto its image}.

\begin{proof}[Proof of Theorem \ref{thm:connected}] Fix a faithful semitoric system $\vat$ and a point $(x_0,y_0)$ $\in \Phi(M)$. Then either $x_0$ is a regular value of $J$ or it is singular. First, suppose that $x_0$ is a~regular value of~$J$. By property~\ref{item:4}, $J^{-1}(x_0)$ is a connected, embedded submanifold of $M$ and, as in \cite[proof of Theorem~4.7]{pvr}, let $H_{x_0}$ denote the function $H|_{J^{-1}(x_0)}\colon J^{-1}(x_0) \to \R$. The condition that the singular orbits of $\Phi$ be almost-toric (property~\ref{item:24}) implies that $H_{x_0}$ is a Morse--Bott function whose critical manifolds have index equal to either~0 or~2 (cf.\ Pelayo, Ratiu and \vungoc\ \cite[Case~1A in proof of Theorem~4.7]{pvr}). Moreover, $H_{x_0}$ is proper onto its image: if $K \subset H_{x_0}\left(J^{-1}(x_0)\right)$ is a~compact subset, then $\{x_0\} \times K$ is a compact subset contained in~$\Phi(M)$. Since $\Phi$ is proper onto its image (by property~\ref{item:23}), $\Phi^{-1}(\{x_0\} \times K) = H_{x_0}^{-1}(K)$ is compact, as desired. If $H_{x_0}\big(J^{-1}(x_0)\big) \subset \R$ is closed, then~$H_{x_0}$ is proper and it is possible to argue as in \cite[Case~1A in proof of Theorem~4.7]{pvr} to deduce that the fibers of~$H_{x_0}$ are connected and, in particular, so is $\Phi^{-1}(x_0,y_0) = H_{x_0}^{-1}(y_0)$. Therefore, suppose that the image of~$H_{x_0}$ is not closed. Since $J^{-1}(x_0)$ is connected, it follows that $H_{x_0}\big(J^{-1}(x_0)\big)$ is either an open interval or a half-open interval. In the former case, compose~$H_{x_0}$ with a diffeomorphism sending the open interval to~$\R$, while, in the latter, compose $H_{x_0}$ with a~diffeomorphism sending the half-open interval to $\R_{\geq 0}$. In both cases, the result is a proper, smooth Morse--Bott function on $J^{-1}(x_0)$ whose fibers agree with those of $H_{x_0}$ and whose critical manifolds and indices are precisely those of $H_{x_0}$, thus reducing this situation to the previous case. In particular, this proves $\Phi^{-1}(x_0,y_0)$ is connected if $x_0$ is a regular value of $J$.

To complete the proof, it is possible to use the strategy of Pelayo, Ratiu and \vungoc\ \cite[Steps~1B,~2 and~3 of the proof of Theorem~4.7]{pvr}, namely:
 \begin{itemize}[leftmargin=*]\itemsep=0pt
 \item show that the fibers of regular values of $\Phi$ are connected,
 \item show that there is no critical value in the interior of $\Phi(M)$ except for focus-focus values, and
 \item show that any singular value of $\Phi$ has an open neighborhood that intersects the set of regular values of $\Phi$ in a connected set.
 \end{itemize}
The above three points allow one to conclude, by Pelayo, Ratiu and \vungoc\ \cite[Lemma~4.6]{pvr}, that all fibers of $\Phi$ are connected. (Observe that the latter holds under the weaker assumption of properness onto the image. Furthermore, it is important to remark that connectedness of~$M$ and property~\ref{item:100} are used crucially in \cite[Steps~2 and~3 of the proof of Theorem~4.7]{pvr}.)
\end{proof}

\subsection{The moment map image of a faithful semitoric system}\label{sec:dscthe-moment-map}
Henceforth, let $\vat$ be a faithful semitoric system whose moment map image is denoted by~$B$. The path-connectedness of the intersection of~$B$ with any vertical line in turn implies that intersections of vertical lines with the interior of~$B$ are path-connected.

\begin{Lemma}\label{lemma:int_conn_vert} Let $\vat$ be a faithful semitoric system whose moment map image is denoted by $B$. For any $x_0 \in \operatorname{pr}_1 (\operatorname{Int}(B) )$, the set $\operatorname{Int}(B) \cap \{(x,y) \,|\, x = x_0\}$ is path-connected.
\end{Lemma}

\begin{proof} Suppose not, then there exists $x_0 \in \operatorname{pr}_1(\operatorname{Int}(B))$ such that $\operatorname{Int}(B) \cap \{(x,y) \,|\, x = x_0\}$ is not path-connected. Since the intersection $B\cap \{(x,y) \,|\, x = x_0\}$ is nonempty and path-connected, it follows that there exist $y_1 < y_0 < y_2$ such that $(x_0,y_0) \in \partial B \subset \mathrm{Bdy}(B)$ and $(x_0, y_i) \in \operatorname{Int}(B)$, for $i=1,2$. However this is impossible because it forces a disconnection of the intersection of~$B$ and a vertical line as follows.

For $i=1,2$, because the point $(x_0,y_i)$ is in $\operatorname{Int}(B)$, there is open disk of radius $r_i$ centered at $(x_0, y_i)$ that is a subset of~$\operatorname{Int}(B)$. And because $(x_0,y_0) \in \partial B \subset \mathrm{Bdy}(B)$, there exists a sequence of points $\{(x'_n,y'_n)\}$ such that each $(x'_n,y'_n)$ is in the open ball of radius $\frac{1}{n}$ centered at $(x_0,y_0)$ but $(x'_n,y'_n)\notin B$.
 Consequently, for each $n$ such that $\frac{1}{n}<\min (r_1,r_2)$ the intersection $B\cap \{(x,y) \,|\, x = x_n \}$ is disconnected because
 $y_1< y_0-\frac{1}{n} < y'_n<y_0+\frac{1}{n} < y_2$, with
 $(x_n',y_1), (x_n',y_2)\in B$ and $(x_n',y_n')\notin B$.
\end{proof}

Using Lemma \ref{lemma:int_conn_vert}, one obtains the following property of moment map images of faithful semitoric systems, which generalizes \vungoc\ \cite[Part~2 of Theorem~3.4]{vu-ngoc} and Pelayo, Ratiu and \vungoc\ \cite[Theorem~C]{pvr_carto}. (Observe that the arguments of {\em op.~cit.} cannot be used without some adjustments as the moment map of a faithful semitoric system is only proper {\em onto its image}.)

\begin{Corollary}\label{cor:contractible_mom_map} The moment map image of a faithful semitoric system is contractible.
\end{Corollary}

\begin{proof} Fix a faithful semitoric system with moment map image $B$. As noted in Lemma~\ref{rmk:almost-toric-image}, $B$ is a manifold with corners. The inclusion $\operatorname{Int}(B) = B \smallsetminus \partial B \hookrightarrow B$ is a~homotopy equivalence because~$B$ is homeomorphic to a smooth manifold with boundary, which is homotopy equivalent to the complement of its boundary. Thus it suffices to prove that $\operatorname{Int}(B)$ is contractible. By Lemma~\ref{lem:breg_connected}, $\operatorname{Int}(B)$ is path-connected and so is $\operatorname{pr}_1 (\operatorname{Int}(B) )$. Moreover, $\operatorname{Int}(B) \subset \R^2$ is open and $\operatorname{pr}_1$ is an open map. Therefore $\operatorname{pr}_1\colon \operatorname{Int}(B) \to \operatorname{pr}_1(\operatorname{Int}(B))$ is a surjective submersion whose fibers are diffeomorphic to~$\R$ by Lemma~\ref{lemma:int_conn_vert}. Thus it is a fiber bundle (cf.\ Meigniez \cite[p.~3778]{meigniez}). Since $\operatorname{pr}_1 (\operatorname{Int}(B) )$ is an interval, the bundle is trivial. Since both the base and the fiber of this trivial bundle are contractible, so is the total space.
\end{proof}

\begin{Corollary}\label{cor:vat+toric=delzant} A faithful semitoric system is visible toric if and only if it is toric.
\end{Corollary}

\subsection{Subsystems of faithful semitoric systems}\label{sec:dscs-faithf-semi}

Faithful semitoric systems behave well with respect to taking certain `vertical subsystems':

\begin{Proposition}\label{prop:sub_mostly_vat} Let $(M,\om,\Phi)$ be a faithful semitoric system and let $U\subset B$ be open and path-connected. Then the subsystem of $(M,\om,\Phi)$ relative to $U$ is a~faithful almost-toric system satisfying properties~{\rm \ref{item:23}--\ref{item:3}} and~{\rm \ref{item:100}--\ref{item:101}} of a~faithful semitoric system. Moreover, if for all $x_0 \in \R$, $\{(x,y) \,|\, x = x_0\} \cap U$ is either empty or path-connected, then the subsystem relative to $U$ is, in fact, a faithful semitoric system.
\end{Proposition}

\begin{proof} The subsystem relative to $U$ is faithful by Corollary~\ref{cor:faithfulness_preserved} and almost-toric by Remark~\ref{rmk:sat_subsyst_at}. The total space is connected because~$U$ is connected and the subsystem is faithful. Property~\ref{item:3} is satisfied because~$\Phi^{-1}(U)$ is a union of orbits of its first integral. Finally, properties~\ref{item:100} and~\ref{item:101} are preserved under taking subsystems. This proves the first assertion.

Assume that the intersection of $U$ with any vertical line is either empty or path-connected. Because the subsystem relative to~$U$ is faithful almost-toric, Lemma~\ref{lemma:connected_J} implies the fibers of~$J|_U$ are path-connected, as desired.
\end{proof}

\subsection{Relation to other families of semitoric systems}\label{sec:dscex-relat-other}

Faithful semitoric systems are very closely related to the proper semitoric systems introduced in Pelayo, Ratiu and \vungoc~\cite{pvr_carto}. Comparing Pelayo, Ratiu, \vungoc\ \cite[Definition~1.1]{pvr_carto} and Definition~\ref{defn:vat} above, the only difference between faithful and proper semitoric systems is that the latter are required to have a~{\em proper} moment map, while the former assume only that the map be {\em proper onto its image}. Accordingly, faithful semitoric systems can be thought of as subsystems of proper semitoric systems. In both cases, the conditions on the systems imply that the fibers are connected (see Theorem~\ref{thm:connected} for faithful semitoric systems, and cf.\ Pelayo, Ratiu, \vungoc\ \cite[Theorem~4.7]{pvr} and Pelayo, Ratiu, \vungoc\ \cite[Remark~1.3]{pvr_carto}). The main reason to introduce faithful semitoric systems is to have a conceptual framework to study integrable surgeries on (proper) semitoric systems that are defined by taking suitable subsystems of (proper) semitoric systems, cf.\ the forthcoming~\cite{HSSS-surgeries}.

\begin{exm}\label{exm:relation} The aim of this example is to illustrate how semitoric systems in the sense of \vungoc~\cite{vu-ngoc}, proper semitoric systems in the sense of Pelayo, Ratiu and \vungoc~\cite{pvr_carto}, and faithful semitoric systems in the sense of Definition~\ref{defn:vat} are related. First, observe that a~semitoric system is necessarily proper semitoric and that a proper semitoric system is necessarily faithful semitoric. Pelayo, Ratiu and \vungoc~\cite{pvr_carto} provide an example of a proper semitoric that fails to be semitoric, namely the spherical pendulum. Next we construct an example of a~faithful semitoric system that fails to be proper semitoric. This is based on Pelayo, Ratiu and \vungoc\ \cite[Section~7]{pvr_carto}. Let~$\omega_{S^2}$ denote the standard symplectic form on~$S^2$ and let $\big(S^2 \times S^2, \omega:= \operatorname{pr}^*_1 \omega_{S^2}+\operatorname{pr}^*_2 \omega_{S^2}, \Phi\big)$ be the toric system whose underlying toric manifold has moment map image equal to the square $[-1,1] \times [-1,1]$. Since the total space is compact, this system is semitoric. As in Pelayo, Ratiu and \vungoc\ \cite[Step~2 of Section~7]{pvr_carto}, consider the open subset $U \subset [-1,1] \times [-1,1]$ that is the complement of $\{0\} \times [0,1]$ (cf.\ Pelayo, Ratiu and \vungoc\ \cite[Fig.~7]{pvr_carto}). The subset $U$ satisfies all the hypotheses of Proposition~\ref{prop:sub_mostly_vat}. Therefore, the subsystem of $\big(S^2 \times S^2, \omega, \Phi\big)$ relative to~$U$ is faithful semitoric. However, observe that it is {\em not} proper semitoric as the restriction of $\Phi$ fails to be proper. (In Pelayo, Ratiu and \vungoc\ \cite[Step~6 of Section~7]{pvr_carto} it is shown to be {\em isomorphic} to a proper semitoric system.)
\end{exm}

While faithful semitoric systems are more restrictive than almost-toric systems (see Definition~\ref{defn:almost-toric}), they capture the semi-local structure of almost-toric systems.

\begin{Proposition}\label{prop:at-vat} Given an almost-toric system $\is$, any leaf $L$ admits an open neighborhood $V$ such that $(V,\omega|_V,\Phi|_V)$ is isomorphic to a faithful semitoric system.
\end{Proposition}

\begin{proof} A leaf of an almost-toric system is either weakly toric or is a~focus-focus leaf. In the former case, the local normal forms for elliptic tori yield the result. Thus, suppose that $L$ is a~focus-focus leaf and let $V$ be a saturated regular neighborhood of $L$. By construction, $V$ is connected; moreover, Proposition~\ref{prop:faithful_at} gives that $(V,\omega|_V,\Phi|_V)$ is faithful almost-toric. Proposition~\ref{prop:S1-action} gives that there exists a system-preserving effective Hamiltonian $S^1$-action on $(V,\omega|_V,\Phi|_V)$; in other words, that subsystem is isomorphic to an almost-toric system whose first integral is the moment map of an effective Hamiltonian $S^1$-action. Thus, without loss of generality, it may be assumed that the first component of $\Phi|_V$ is the moment map of an effective Hamiltonian $S^1$-action. By construction of~$V$, properties~\ref{item:100} and~\ref{item:101} of a faithful semitoric system hold for $(V,\omega|_V,\Phi|_V)$. Set $c = \Phi(L)$. Let $U \subset \Phi|_V(V)$ be an open neighborhood of~$c$ with the property that its intersection with any vertical line is either empty or path-connected. Using Proposition~\ref{prop:sub_mostly_vat}, the subsystem relative to $U$ of $(V,\omega|_V,\Phi|_V)$ is faithful semitoric as desired.
\end{proof}

\begin{rmk}\label{rmk:not_generalized} The faithful semitoric systems of Proposition \ref{prop:at-vat} are not necessarily proper semitoric because their moment map images are not necessarily closed in~$\R^2$, while those of proper semitoric systems are.
\end{rmk}

\subsection{Notions of isomorphism of faithful semitoric systems}\label{sec:dscp-noti-isom}

The notions in the literature for isomorphisms of (proper) semitoric systems (cf.\ Pelayo, \vungoc\ \cite[Section~2]{pelayo-vu-ngoc-inventiones} and Pelayo, Ratiu, \vungoc\ \cite[Definition~1.5]{pvr_carto}) are, {\em a~priori}, more restrictive than the notion of isomorphism of integrable systems introduced in Definition~\ref{defn:cihs}: they are based on the fact that (proper) semitoric systems can be viewed as symplectic 4-manifolds endowed with a Hamiltonian $S^1 \times \R$-action. Seeing as the present paper is written from the point of view of integrable systems but with a view to study the relation between semitoric systems and Hamiltonian $S^1$-spaces via integrable surgeries (cf.\ the forthcoming~\cite{HSSS-surgeries,HSSS-lifting}), we introduce two other notions of isomorphism and illustrate how the four are related.

\begin{rmk}\label{rmk:im_faithful} The term {\em complexity one} in Definition \ref{defn:im_faithful_st} comes from Karshon and Tolman~\cite[Definition~1.2]{kt1} where a~four-dimensional {\em complexity one space} is a~connected symplectic four-manifold equipped with an effective Hamiltonian $S^1$-action whose moment map is proper onto its image. So, given a faithful semitoric system~$\vat$, forgetting the second integ\-ral~$H$ yields a~complexity one space.
 \end{rmk}

 \begin{Definition}\label{defn:im_faithful_st} Two faithful semitoric systems $(M_1,\omega_1,\Phi_1)$ and $(M_2,\omega_2,\Phi_2)$ are
 \begin{itemize}[leftmargin=*]\itemsep=0pt
 \item {\em isomorphic as integrable systems} if they are isomorphic as in Definition~\ref{defn:cihs},
 \item {\em isomorphic as complexity one systems} if they are isomorphic as integrable systems via a pair $(\Psi,\psi)$, where $\psi \colon B_1 \to B_2$ is of the form $\psi = \big(\psi^{(1)}, \psi^{(2)}\big)$ with $\psi^{(1)}(x,y) = \zeta x + a$, for some $\zeta \in \{+1,-1\}$ and $a \in \R$,
 \item {\em strictly isomorphic as complexity one systems} if they are isomorphic as complexity one systems and, in addition, $\zeta = +1$ and $a=0$,
 \item {\em isomorphic as semitoric systems} if they are strictly isomorphic as complexity one systems and, in addition, $\psi$ is orientation-preserving.
 \end{itemize}
 \end{Definition}

Note that the notion of isomorphism as semitoric systems is the one commonly used in the literature for proper semitoric systems (cf.\ Pelayo, Ratiu and \vungoc\ \cite[Definition~1.5]{pvr_carto}).

It is clear that the notions of isomorphism in Definition~\ref{defn:im_faithful_st} are listed from the coarsest (isomorphism as integrable systems) to the finest (isomorphism as semitoric systems). The next remark explains how, in general, a coarser isomorphism does not imply the existence of any finer~one.

 \begin{rmk}\label{rmk:toric_im} Suppose that $(M_1,\omega_1,\Phi_1)$ and $(M_2,\omega_2,\Phi_2)$ are closed toric systems on connected symplectic 4-manifolds. Then, viewing these systems as being faithful semitoric, Proposition~\ref{prop:ds} yields that
 \begin{itemize}[leftmargin=*]\itemsep=0pt
 \item they are isomorphic as integrable systems if and only if there exists an element $h \in \mathrm{AGL}(2;\Z)$ such that $\Phi_2 = h \circ \Phi_1$,
 \item \looseness=-1 they are isomorphic as complexity one systems if and only if there exists an element $h \in \mathrm{L}$ such that $\Phi_2 = h \circ \Phi_1$, where $\mathrm{L} \subset \mathrm{AGL}(2;\Z)$ is the subgroup consisting of
 elements of the form
 \begin{gather*}
 \left( \begin{pmatrix}
 \eta_1 & 0 \\
 k & \eta_2
 \end{pmatrix},
 \begin{pmatrix}
 a \\
 b
 \end{pmatrix}\right),
 \end{gather*}
 where, for $i=1,2$, $\eta_i \in \{+1,-1\}$, $k \in \Z$, and $a,b \in \R$,
 \item they are strictly isomorphic as complexity one systems if and only if there exists an element $h \in \operatorname{Vert}(2;\Z)$ such that $\Phi_2 = h \circ \Phi_1$, where $\operatorname{Vert}(2;\Z) \subset
 \mathrm{AGL}(2;\Z)$ is the subgroup consisting of elements of the form{\samepage
 \begin{gather*}
 \left( \begin{pmatrix}
 1 & 0 \\
 k & \pm 1
 \end{pmatrix},
 \begin{pmatrix}
 0 \\
 b
 \end{pmatrix}\right),
 \end{gather*}
 where $k \in \Z$ and $b \in \R$, and}

 \item they are isomorphic as semitoric systems if and only if there exists an element $h \in \operatorname{Vert}^+(2;\Z)$ such that $\Phi_2 = h \circ \Phi_1$, where $\operatorname{Vert}^+(2;\Z) \subset
 \operatorname{Vert}(2;\Z)$ is the normal subgroup consisting of elements whose linear part is orientation-preserving.
 \end{itemize}
 Accordingly, given any of the coarser notions of isomorphism, it is possible to construct examples of closed toric systems that are isomorphic in the given sense but not in any finer sense.
 \end{rmk}

However, in the presence of focus-focus points isomorphisms as integrable systems {\em are necessarily} isomorphisms as complexity one spaces.

\begin{Proposition}\label{prop:im_equ} For $i=1,2$, let $(M_i,\omega_i,\Phi_i)$ be a faithful semitoric system and suppose that $(M_1,\omega_1,\Phi_1)$ contains at least one focus-focus point. Then $(M_1,\omega_1,\Phi_1)$ and $(M_2,\omega_2,\Phi_2)$ are isomorphic as integrable systems if and only if they are isomorphic as complexity one spaces.
\end{Proposition}

The above result is not hard to prove and is probably well-known amongst experts. However, we could not find it anywhere in the literature and it shows how the choice of isomorphisms of (proper) semitoric systems in Pelayo, \vungoc\ \cite[Section~2]{pelayo-vu-ngoc-inventiones} and Pelayo, Ratiu, \vungoc\ \cite[Definition~1.5]{pvr_carto} fit naturally into the broad problem of classifying integrable systems. Its proof is postponed to Section~\ref{sec:cuts-cart-home}.

\subsection[Geometric implications of a global $S^1$-action]{Geometric implications of a global $\boldsymbol{S^1}$-action}\label{sec:geom-impl-s1}
Fix a faithful semitoric system $\vat$ and denote the set of fixed points of the Hamiltonian $S^1$-action, one of whose moment maps is~$J$, by~$M^{S^1}$. Its connected components are either isolated fixed points or {\em symplectic fixed surfaces}, i.e., symplectic submanifolds of dimen\-sion~2 that are fixed under the $S^1$-action. This is a~consequence of the Marle--Guillemin--Sternberg local normal form (cf.\ Guillemin, Sternberg~\cite{GS2} and Marle~\cite{marle}).

\begin{Proposition}\label{prop:top_fixed_surfaces} Let $\vat$ be a faithful semitoric system such that $M^{S^1}$ contains a fixed surface~$\Sigma$. Then~$J(\Sigma)$ is a global extremum of~$J$.
\end{Proposition}

\begin{proof} To show that $J(\Sigma)$ is a global extremum of $J$ it suffices to show that it cannot lie in the interior of the interval $J(M)$. Assume the contrary: then $J^{-1}(J(\Sigma)) = \Sigma$ as the fibers of~$J$ are connected by assumption. Since $J(\Sigma) \in \operatorname{Int}(J(M))$, it follows that $J(M) \smallsetminus J(\Sigma)$ is disconnected, thus implying that $M \smallsetminus \Sigma = J^{-1}(J(M) \smallsetminus J(\Sigma))$ is disconnected. However, this is absurd, since~$M$ is connected and $\Sigma \subset M$ is a submanifold of codimension~2.
\end{proof}

In fact, given a faithful semitoric system $\vat$ for which $M^{S^1}$ contains a fixed surface $\Sigma$, each point of $\Sigma$ belongs to some singular orbit of elliptic type, so $\Phi(\Sigma) \subset \partial B$. Moreover $\Phi(\Sigma) \subset \{ (x,y) \,|\, x = J(\Sigma) \}$. The image $\Phi(\Sigma) \subset B$ is henceforth referred to as a~{\em vertical edge} of~$B$. The following result provides a~characterization of vertical edges.

\begin{Proposition}\label{prop:vert_edge} Let $\vat$ be a faithful semitoric system and suppose that there exists $x_0 \in J(M)$ and distinct points $c^{\infty}_1, c^{\infty}_2, c^{\infty}_3 \in \partial B \cap \{ (x,y) \,|\, x = x_0\}$. Then $\partial B \cap \{ (x,y) \,|$ $x = x_0\}$ is a vertical edge.
\end{Proposition}

\begin{proof} For $i=1,2,3$, set $c^{\infty}_i = (x_0,y_i)$ and assume, without loss of generality, that $y_1 < y_2 < y_3$. Since the fibers of $J$ are connected by property~\ref{item:4}, it follows that
 \begin{gather*}
\{ (x,y) \,|\, x=x_0\, , \, y_1 \leq y \leq y_3 \} \subset B.
 \end{gather*}
The local normal forms for almost-toric singular orbits in the presence of a system-preserving Hamiltonian $S^1$-action (see Remark~\ref{rmk:lt}), together with faithfulness of $\vat$, force $\Phi^{-1}\big(c^{\infty}_2\big)$ to be a singular orbit of elliptic-regular type, all of whose points are critical for~$J$. Then $\Phi^{-1}\big(c^{\infty}_2\big)$ lies on a fixed surface and Proposition~\ref{prop:top_fixed_surfaces}, together with connectedness of the fibers of~$J$, completes the proof.
\end{proof}

Combining Propositions \ref{prop:top_fixed_surfaces} and~\ref{prop:vert_edge}, we obtain the following result.

\begin{Corollary}\label{cor:int_interior} Let $\vat$ be a faithful semitoric system and consider a point $x_0 \in \operatorname{Int}(J(M))$. Then the intersection $\partial B \cap \{ (x,y) \,|\, x = x_0\}$ consists of at most two points.
\end{Corollary}

To conclude this section, we state the following result, which isentirely analogous to Hohloch, Sabatini and Sepe \cite[Lemma~3.3]{HSS}.

\begin{Proposition}\label{lemma:iso_fixed_points}
 Let $\vat$ be a faithful semitoric system. Then the isolated fixed points in $M^{S^1}$ are either
 \begin{itemize}[leftmargin=*]\itemsep=0pt
 \item focus-focus singular orbits, or
 \item elliptic-elliptic singular orbits whose image is not a corner adjacent to a vertical edge.
 \end{itemize}
\end{Proposition}

\section[Cartographic homeomorphisms and $\eta$-cartographic representatives]{Cartographic homeomorphisms\\ and $\boldsymbol{\eta}$-cartographic representatives}\label{sec:dscc-home-boldsymb}
This section proves that faithful semitoric systems admit cartographic homeomorphisms constructed by choosing suitable vertical cuts of their moment map images (see Section~\ref{sec:cuts-cart-home} and Theorem~\ref{prop:rh}). The set of all cartographic homeomorphisms of a given faithful semitoric system is described in Section~\ref{sec:set-cart-home}, which generalizes \vungoc~\cite[Section~4]{vu-ngoc}. This set is used to construct an invariant of faithful semitoric systems up to various notions of isomorphisms (see Lemmas~\ref{lemma:im_invariant} and~\ref{lemma:im_inv_full}). These results can be seen as first steps toward classifying faithful semitoric systems up to any of the notions of isomorphism of Section~\ref{sec:dscp-noti-isom}. Finally, Section~\ref{sec:choos-an-appr} shows that cartographic homeomorphisms can be made smooth by modifying them on arbitrarily small neighborhoods of the corresponding cuts (Theorem~\ref{thm:rect-embedding}). This result provides representatives (which we call {\em $\boldsymbol{\eta}$-cartographic}) in the isomorphism class of any faithful semitoric system (for any
 notion of isomorphism given in Section~\ref{sec:dscp-noti-isom}).

\subsection{Vertical cuts and existence of cartographic homeomorphisms}\label{sec:cuts-cart-home}
A fundamental property of faithful semitoric systems is that, as proved below, assuming the following mild restriction, they admit cartographic homeomorphisms (see Theorem~\ref{prop:rh}).

 \begin{assumption} Any faithful semitoric system is henceforth assumed to satisfy the following pro\-per\-ty:
 \begin{enumerate}[label={\rm (F\arabic*)}, ref=(F\arabic*), leftmargin=*, start=7]\itemsep=0pt
 \item \label{item:7} Any focus-focus value has multiplicity~$1$.
 \end{enumerate}
 \end{assumption}

The above condition is generic according to Zung \cite{zung-symplectic}. Moreover, it is invariant under any of the notions of isomorphism of faithful semitoric systems of Definition~\ref{defn:im_faithful_st}, and descends to saturated subsystems satisfying all the hypotheses of Proposition~\ref{prop:sub_mostly_vat}.

While the existence of cartographic homeomorphisms is expected to hold without imposing property~\ref{item:7}, proofs would require a more detailed understanding of neighborhoods of focus-focus fibers with more than one focus-focus point, which is beyond the scope of this paper (cf.\ \vungoc\ \cite[Section~7]{vu_ngoc_invariant} and Pelayo, Tang~\cite{pel_tang} for results in this direction). To the best of our knowledge, all existing proofs of the existence of cartographic homeomorphisms assume, either tacitly or explicitly, property~\ref{item:7} (cf.\ \vungoc\ \cite[Step~4 of the proof of Theorem~3.8]{vu-ngoc} and Pelayo, Ratiu, \vungoc\ \cite[Step~4 of the proof of Theorem~B]{pvr_carto}).

\begin{rmk}\label{rmk:simple} Note that property \ref{item:7} is weaker than the notion of {\em simple} used in the literature, as the latter means that there exists at most one focus-focus point on any fiber of~$J$ (cf.\
 Pelayo and \vungoc\ \cite[Definition~3.3]{pelayo-vu-ngoc-inventiones}).
\end{rmk}

The aim of this section is to prove that any faithful semitoric system satisfying property~\ref{item:7} admits a cartographic homeomorphism that, loosely speaking, encodes the affine holonomy of the $\Z$-affine structure on
the weakly toric part of the leaf space (see Theorem~\ref{prop:rh} for a precise statement). It is important to remark that there are proofs of this result for special families of faithful semitoric systems (cf.\ \vungoc\ \cite[Theorem~3.8]{vu-ngoc} and Pelayo, Ratiu, \vungoc\ \cite[Theorem~B]{pvr_carto} for semitoric and proper semitoric systems respectively). Those proofs are utilized and adjusted as needed in what follows.

\begin{Lemma}\label{lemma:vat_toric} Let $\vat$ be a faithful semitoric system without focus-focus points whose moment map image is denoted by $B$. Then there exists a cartographic pair $(f,B)$, where~$f$ is of the form
 \begin{gather} \label{eq:8}
 f(x,y) = \big(f^{(1)}, f^{(2)} \big)(x,y) = \big(x,f^{(2)}(x,y)\big).
 \end{gather}
In particular, $(M,\omega, f \circ \Phi)$ is a visible toric system and $f(B) \subset \R^2$ is locally convex.
\end{Lemma}

\begin{proof} The proof is analogous to Pelayo, Ratiu and \vungoc\ \cite[Step~2 of Theorem~B]{pvr_carto}, but is included in this paper for completeness.

The lack of focus-focus points implies $B = B_{\mathrm{wt}}$, thus $B$ inherits a $\Z$-affine structure. By Corollary~\ref{cor:contractible_mom_map}, $B$ is contractible, so there exists a developing map $f \colon B \cong \tilde{B} \to \R^2$. By definition of the $\Z$-affine structure on $B$, since the first integral~$J$ of~$\vat$ is the moment map of an effective Hamiltonian $S^1$-action, one can choose the above developing map to be of the form $f(x,y) = \big(x, f^{(2)}(x,y)\big)$ for some smooth function $f^{(2)}\colon B \to \R$. Fix such a choice.

To show that $f$ is the required cartographic homeomorphism, it suffices to show that~$f$ is injective. If $f(x_0,y_0) = f(x_1,y_1)$, one gets immediately $x_0 = x_1$. The map $f^{(2)}(x_0, \cdot) \colon B \cap \{(x,y) \,|\, x = x_0 \} \to \R$ is strictly monotone as $\frac{\partial f^{(2)}}{\partial y}$ does not vanish on~$B$, because $f$ is locally a~diffeomorphism. This implies that $y_0 = y_1$ as required.

Since $\vat$ is faithful semitoric and because of the form of $f$, $(M, \omega, f\circ \Phi)$ is faithful semitoric. In fact, it is toric because $f \circ \Phi \colon \sm \to \R^2$ is the moment map of an effective Hamiltonian $\mathbb{T}^2$-action. By Corollary~\ref{cor:vat+toric=delzant}, $(M, \omega, f \circ \Phi)$ is visible toric and by Corollary~\ref{cor:delz_sys_loc_convex}, $f (B)$ is locally convex.
\end{proof}

{\bf Intermezzo: classification of faithful semitoric systems with no focus-focus points via their cartographic images.}
Lemma \ref{lemma:vat_toric} allows one to classify faithful semitoric systems with no focus-focus points up to isomorphisms of integrable systems. This is achieved by understanding their sets of cartographic pairs and corresponding images.

\begin{Lemma}\label{lemma:max_carto} Let $\is$ be a faithful semitoric system with no focus-focus points, let $B$ denote its moment map image and fix a cartographic pair $(f,B)$ for $\is$ as in Lemma~{\rm \ref{lemma:vat_toric}}. The set of cartographic pairs of $\is$ whose second component is $B$ is
 \begin{gather*}
\{ (h \circ f,B ) \,|\, h \in \mathrm{AGL}(2;\Z) \}.
 \end{gather*}
 \end{Lemma}
\begin{proof} Corollary \ref{lem:infinitely_many_carto} implies that, for any $h \in \mathrm{AGL}(2;\Z)$, $(h \circ f, B)$ is also a cartographic pair for~$\is$. Conversely, suppose that $(\hat{f},B)$ is a cartographic pair for~$\is$. The map $\hat{f} \circ f^{-1} \colon f(B) \to \hat{f}(B)$ is a $\Z$-affine isomorphism, where $f(B)$ and $\hat{f}(B)$ are endowed with the restriction of the standard $\Z$-affine structure on $\R^2$. Since $f(B)$ and $\hat{f}(B)$ are connected, it follows that there exists $h \in \mathrm{AGL}(2;\Z)$ such that $h = \hat{f} \circ f^{-1}$, as desired.
 \end{proof}

\begin{rmk}\label{rmk:freedom_no_ff} With the hypotheses of Lemma \ref{lemma:max_carto}, a~cartographic homeomorphism $(\hat{f},B)$ of $\is$ is of the form given by equation~\eqref{eq:8} if and only if the element $h \in \mathrm{AGL}(2;\Z)$ constructed in the above proof belongs to the subgroup $\operatorname{Vert}(2;\Z)$ consisting of $\Z$-affine transformations that fix all vertical lines in~$\R^2$ (see Remark~\ref{rmk:toric_im}).
\end{rmk}

Cartographic pairs of $\is$ whose second component is $B$ are henceforth referred to as being {\em maximal} and so are the corresponding cartographic images (see Remark~\ref{rmk:smooth_locus}). Before stating the classification of faithful semitoric systems with no focus-focus points (up to isomorphisms of integrable systems, see Lemma~\ref{lemma:class_no_ff}), observe that the number of focus-focus points is an invariant of the isomorphism class of a faithful semitoric system (see Section~\ref{sec:at_generalities}).

 \begin{Lemma}\label{lemma:class_no_ff} For $i=1,2$ let $(M_i,\omega_i,\Phi_i)$ be faithful semitoric systems with no focus-focus points. Then $(M_1,\omega_1,\Phi_1)$ and $(M_2,\omega_2,\Phi_2)$ are isomorphic as integrable
 systems if and only if their sets of maximal cartographic images are equal.
 \end{Lemma}
\begin{proof} Suppose first that $(M_1,\omega_1,\Phi_1)$ and $(M_2,\omega_2,\Phi_2)$ are isomorphic as integrable systems. Then Corollary~\ref{cor:pull_back_iso_carto} implies that the set of maximal cartographic images of the former is included in that of the latter. Reversing the above argument, we obtain the desired equality. Conversely, suppose that their sets of maximal cartographic images are equal. Then, for $i=1,2$, there exists a~cartographic pair $(f_i,B_i)$ with the property that $f_1(B_1) = f_2 (B_2 )$. Lemmas~\ref{lemma:vat_toric} and~\ref{lemma:max_carto} imply that, for $i=1,2$, $(M_i,\omega_i,f_i\circ \Phi_i)$ are visible toric systems; moreover, the above assumption gives that their moment map images are equal. By Proposition~\ref{prop:ds}, this implies that there exists a symplectomorphism $\Psi \colon (M_1,\omega_1) \to (M_2,\omega_2)$ such that $f_2 \circ \Phi_2 \circ \Psi = f_1 \circ \Phi_1$. Set $\psi:=f_2^{-1} \circ f_1 \colon B_1 \to B_2$. The map $\psi$ is a diffeomorphism since both $f_1$ and $f_2$ are diffeomorphisms and it satisfies $\Phi_2 \circ \Psi = \psi \circ \Phi_1$. Thus $(\Psi,\psi)$ is an isomorphism between $(M_1,\omega_1,\Phi_1)$ and $(M_2,\omega_2,\Phi_2)$ as desired.
 \end{proof}

Using Lemma \ref{lemma:class_no_ff}, it is possible to deduce the classification of faithful semitoric systems with no focus-focus points up to any notion of isomorphism given in Definition~\ref{sec:dscp-noti-isom} (see Pelayo, Ratiu and \vungoc~\cite[Lemma~3.3]{pvr_carto} for a proof of the classification of proper semitoric systems according to the finest isomorphism type).

\medskip

\begin{figure}[h] \centering
 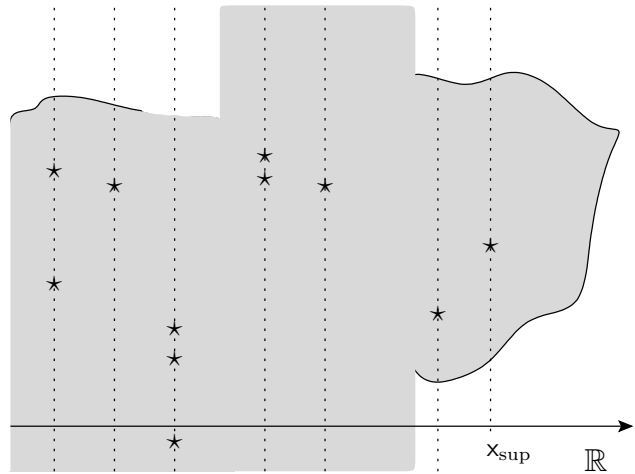
 \caption{The image of the moment map (gray) with the focus-focus values (marked by $\star$) and their projection onto the first component.}\label{imageMomentMap}
\end{figure}

{\bf Cartographic homeomorphisms for faithful semitoric systems with at least one focus-focus point.} Let $\vat$ be a faithful semitoric system with $B =\Phi(M)$ that contains at least one focus-focus point. An example is
shown in Fig.~\ref{imageMomentMap}. We introduce vertical `cuts' at the focus-focus values, along what Symington~\cite{symington} refers to as `eigenrays' in the more general context of (faithful) almost-toric systems. Our terminology is motivated by \vungoc~\cite{vu-ngoc}.

Let $B_{\mathrm{ff}} \subset\operatorname{Int}(B)$ denote the set of focus-focus values. By Lemma~\ref{rmk:almost-toric-image}, it is a countable subset. To order the elements of $B_{\mathrm{ff}}$ we fix the following convention for the indexing set of $B_{\mathrm{ff}}$. Set
\begin{gather} \label{eq:9}
\mathsf{x}_{\sup} := \sup \{\operatorname{pr}_1(c) \,|\, c \in B_{\mathrm{ff}}\}, \qquad \mathsf{x}_{\inf} := \inf \{\operatorname{pr}_1(c) \,|\, c \in B_{\mathrm{ff}}\},
\end{gather}
where $\operatorname{pr}_1 \colon \R^2 \to \R$ is projection onto the first component. By property~\ref{item:100}, this sup\-re\-mum~$\mathsf{x}_{\sup}$ (respectively infimum~$\mathsf{x}_{\inf}$) is either attained as a maximum (respectively as a minimum) or does not lie in~$J(M)$. Set
\begin{gather} \label{eq:2}
 I:= \begin{cases}
 \{1,2,\ldots, \lvert B_{\mathrm{ff}} \rvert\} & \text{if } \lvert B_{\mathrm{ff}} \rvert < \infty, \\
 \{1,2,\ldots\} & \text{if } \lvert B_{\mathrm{ff}} \rvert = \lvert \N\rvert \text{ and } \mathsf{x}_{\inf} \in J(M), \\
 \{0, -1, -2,\ldots\} & \text{if } \lvert B_{\mathrm{ff}} \rvert = \lvert \N\rvert \text{ and } \mathsf{x}_{\sup} \in J(M), \\
 \Z & \text{otherwise.}
 \end{cases}
\end{gather}
By construction, the cardinality of $I$ equals that of~$B_{\mathrm{ff}}$ and thus we think of the elements of the latter as being indexed by~$I$. Order the elements of $B_{\mathrm{ff}}$ as follows. For $i \in I$, set $c_i =
(x_i,y_i)$. Then require that $i<j$ implies either $x_i < x_j$, or $x_i = x_j$ and $y_i < y_j$; moreover, if $0,1 \in I$, require that $x_0 < x_1$. (If $I = \Z$, the above ordering is unique up to the choice of which focus-focus value is labeled with $0$.)

\begin{figure}[h] \centering
 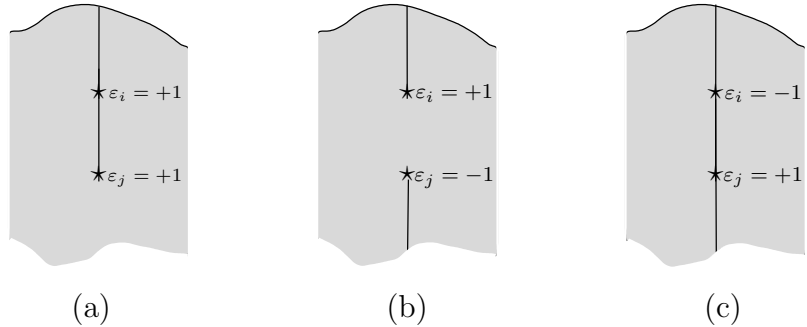
 \caption{Image of the moment map (gray) with cuts emanating from the focus-focus values (marked by~$\star$). The choice of cuts in (a) and (b) leads to a simply connected set whereas the choice in (c) yields two connected components.} \label{cut}
\end{figure}

For each $i \in I$ choose a sign $\ep_i \in \{+1,-1\}$, and denote the associated vertical {\em cut} in $B$ at~$c_i$ by
\begin{gather*}
 l^{\ep_i} := \big\{(x,y) \in \R^2 \,|\, x =x_i ,\, \ep_iy \geq \ep_i y_i \big\} \cap B.
\end{gather*}
When $\epsilon_i = +1$ (respectively $-1$), the cut $l^{\ep_i}$ is simply the intersection of $B$ with the vertical half-line starting at $c_i$ going `up' (respectively `down'), see Fig.~\ref{cut}. Therefore the former is referred to as being {\em upward}, while the latter as being {\em downward}. For a~fixed $\boldsymbol{\ep} \in \{+1,-1\}^I$, denote the union of the cuts by $l^{\boldsymbol{\ep}}$ and set $S_{\boldsymbol{\ep}}:= B
\smallsetminus l^{\boldsymbol{\ep}}$. Moreover, to each element $(x,y) \in B_{\mathrm{wt}} = B \smallsetminus B_{\mathrm{ff}}$,
associate the integer
\begin{gather*}
 j_{\boldsymbol{\ep}}(x,y) :=\sum\limits_{\{i \in I \,|\, (x,y) \in l^{\ep_i}\}} \ep_i,
\end{gather*}
with the convention that $j_{\boldsymbol{\ep}}(x,y) = 0$ for $(x,y) \in S_{\boldsymbol{\ep}}$. Finiteness of $j_{\boldsymbol{\ep}}(x,y)$ follows from pro\-per\-ty~\ref{item:101} and Proposition~\ref{lemma:iso_fixed_points}. The quantity $j_{\boldsymbol{\ep}}(x,y)$ is a signed count of the number of cuts that pass through $(x,y)$, where upward cuts are counted positively and downward cuts negatively.

\begin{Corollary}\label{cor:complement_open} The subset $S_{\boldsymbol{\ep}}$ is open and dense in $B$.
\end{Corollary}

\begin{proof} Density of $S_{\boldsymbol{\ep}}$ in $B$ is trivial, so it remains to prove its openness in $B$ and, to this end, it suffices to prove that $l^{\boldsymbol{\ep}}$ is closed in~$B$. Let $\{(x_n,y_n)\} \subset l^{\boldsymbol{\ep}}$ be a sequence that converges to $(x_0,y_0) \in B$. This implies that the sequence $\{x_n\} = \operatorname{pr}_1(\{(x_n,y_n)\}) \subset J(M)$ converges to $x_0 \in J(M)$. By construction and by Proposition~\ref{lemma:iso_fixed_points}, $\{x_n\}$ is contained in the subset of critical values of~$J$, which does not contain any limit points in~$J(M)$ by property~\ref{item:100}. Therefore, for all but finitely many~$n$, $x_n = x_0$, which, in turn, implies that $(x_n,y_n) \in \{ (x,y) \,|\, x = x_0\}$ for all but finitely many~$n$. By property~\ref{item:101} and Proposition~\ref{lemma:iso_fixed_points}, the vertical line $\{ (x,y) \,|\, x = x_0\}$ contains finitely many focus-focus values and, therefore, finitely many cuts. Seeing as each cut is a closed subset, then the union of all cuts contained on $\{ (x,y) \,|\, x = x_0\}$ is closed. Therefore, $(x_0,y_0) \in l^{\boldsymbol{\ep}}$ as required.
\end{proof}

\begin{rmk}\label{rmk:not_closed} In general, it is not true that $l^{\boldsymbol{\ep}}$ is closed in $\R^2$, for $\mathsf{x}_{\sup}, \mathsf{x}_{\inf}$ may belong to $\R \smallsetminus J(M)$.
\end{rmk}

The notation $S_{\boldsymbol{\ep}}$ is suggestive of the fact that there exists a cartographic homeomorphism $f_{\boldsymbol{\ep}} \colon B \to \R^2$ such that $(f_{\boldsymbol{\ep}}, S_{\boldsymbol{\ep}})$ is a~cartographic pair for~$\vat$. Before stating and proving the precise existence statement (Theorem~\ref{prop:rh} below), we prove some further properties of~$S_{\boldsymbol{\ep}}$ (see Fig.~\ref{cut}).

\begin{Lemma}\label{lemma:S_ep_path_conn} The subset $S_{\boldsymbol{\ep}}$ is path-connected if and only if $\ep_i \geq \ep_j$ for all $i >j $ with $x_i = x_j$.
\end{Lemma}

\begin{proof} Suppose first that $S_{\boldsymbol{\ep}}$ is path-connected and let $i >j \in I$ be such that $x_i = x_j$. Let $(x_1,y_1), (x_2,y_2) \in S_{\boldsymbol{\ep}}$ be points with $x_1 < x_i < x_2$. Such points exist because focus-focus values are contained in $\operatorname{Int}(B)$. Since $S_{\boldsymbol{\ep}}$ is path-connected, there exists a path in $S_{\boldsymbol{\ep}}$ starting at $(x_1,y_1)$ and ending at $(x_2,y_2)$. Therefore, there exists a point $(x_i,y') \in S_{\boldsymbol{\ep}}$ and thus $\ep_i \geq \ep_j$.

Conversely, suppose that $\boldsymbol{\ep}$ satisfies the condition that $\ep_i \geq \ep_j$ for all $i >j $ with $x_i = x_j$. First we show that, for all $\mathsf{x}_1 \in J(M)$, $\{(x,y) \,|\, x= \mathsf{x}_1 \} \cap S_{\boldsymbol{\ep}} \neq \varnothing$ and that the set is path-connected. If $\mathsf{x}_1 \notin \operatorname{pr}_1(B_{\mathrm{ff}})$, one obtains
 \begin{gather*}
\{(x,y) \,|\, x= \mathsf{x}_1 \} \cap S_{\boldsymbol{\ep}} = \{(x,y) \,|\, x= \mathsf{x}_1 \} \cap B
 \end{gather*}
and the result follows from the fact that $\vat$ is a faithful semitoric system. Suppose, therefore, that $\mathsf{x}_1 \in \operatorname{pr}_1(B_{\mathrm{ff}})$. By property~\ref{item:101} and Proposition~\ref{lemma:iso_fixed_points}, there are finitely many focus-focus values $(\mathsf{x}_1,y_{i_1}), \ldots, (\mathsf{x}_1,y_{i_N})$ lying on the vertical line $\{(x,y) \,|\, x= \mathsf{x}_1 \}$. Set
 \begin{gather*}
 y_{+} := \inf \{ y_{i_k} \,|\, \ep_{i_k} = +1\}, \qquad y_{-} := \sup \{ y_{i_k} \,|\, \ep_{i_k} = -1\}.
 \end{gather*}
Since $\boldsymbol{\ep}$ satisfies the condition in the statement, it follows that $y_+ > y_-$ and therefore,
 \begin{gather*}
\{(x,y) \,|\, x= x_1 \} \cap S_{\boldsymbol{\ep}} = \{ (x,y) \,|\, x=\mathsf{x}_1\, , \, y_{+} > y > y_{-}\},
 \end{gather*}
which shows that $ \{(x,y) \,|\, x= \mathsf{x}_1 \} \cap S_{\boldsymbol{\ep}}$ is path-connected. By Corollary~\ref{cor:complement_open}, $S_{\boldsymbol{\ep}}$ is open in~$B$. Thus~$S_{\boldsymbol{\ep}}$ satisfies all the hypotheses of Proposition~\ref{prop:sub_mostly_vat} and, therefore, the subsystem of~$\vat$ relative to~$S_{\boldsymbol{\ep}}$ is faithful semitoric. By Corollary~\ref{cor:contractible_mom_map}, $S_{\boldsymbol{\ep}}$ is
 contractible and, in particular, path-connected.
\end{proof}

\begin{Corollary}\label{cor:path_conn_iff_contr} The subset $S_{\boldsymbol{\ep}}$ is path-connected if and only if it is contractible.
\end{Corollary}

\begin{Corollary}\label{cor:there_exists_choice} There exists a choice of $\boldsymbol{\ep}$ making $S_{\boldsymbol{\ep}}$ path-connected.
\end{Corollary}
\begin{proof} The choice of $\ep_i = +1$ for all $i \in I$ satisfies the condition of Lemma~\ref{lemma:S_ep_path_conn}.
\end{proof}

Having established the above preliminary results, we can state and prove existence of cartographic homeomorphisms for faithful semitoric systems.

\begin{Theorem}\label{prop:rh} Let $\vat$ be a faithful semitoric system with $B_{\mathrm{ff}} = \{c_i\}_{i \in I} \neq \varnothing$. For any $\boldsymbol{\ep} \in \{+1,-1\}^I$, there exists a~cartographic pair $(f_{\boldsymbol{\ep}},S_{\boldsymbol{\ep}})$, where $S_{\boldsymbol{\epsilon}} = B \smallsetminus l^{\boldsymbol{\epsilon}}$ is the complement of the vertical cuts determined by $\boldsymbol{\ep}$ and $f_{\boldsymbol{\ep}}$
 is of the form $f_{\boldsymbol{\ep}}(x,y)= \big(x,f_{\boldsymbol{\ep}}^{(2)}(x,y)\big)$, satisfying the following properties
 \begin{enumerate}[label={\rm (C\arabic*)}, ref=(C\arabic*), leftmargin=*]\itemsep=0pt
 \item \label{item:8} for all $(x,y) \in S_{\boldsymbol{\ep}}$, the quantity $\operatorname{sgn}\big(\frac{\partial f^{(2)}_{\boldsymbol{\ep}}}{\partial y} (x,y) \big)=: \operatorname{sgn}(f_{\boldsymbol{\ep}})$ is
 constant;
 \item \label{item:17} for all $(\mathsf{x},\mathsf{y}) \in B_{\mathrm{wt}}$,
 \begin{gather*} 
 \lim_{\stackrel{(x, y) \to (\mathsf{x},\mathsf{y})}{x <\mathsf{x}}} D f_{\boldsymbol{\ep}} (x, y)= \begin{pmatrix}
 1 & 0 \\ \operatorname{sgn}(f_{\boldsymbol{\ep}}) j_{\boldsymbol{\ep}}(\mathsf{x},\mathsf{y}) & 1
 \end{pmatrix}
 \lim_{\stackrel{(x, y) \to (\mathsf{x},\mathsf{y})}{x >\mathsf{x}}} D f_{\boldsymbol{\ep}}(x, y).
 \end{gather*}
 \end{enumerate}
 In particular, $f_{\boldsymbol{\ep}}(B)$ is locally convex.
\end{Theorem}

Any cartographic homeomorphism $f_{\boldsymbol{\ep}} \colon B \to \R^2$ satisfying the properties in Theorem~\ref{prop:rh} is said to be {\em associated} to~$\boldsymbol{\ep}$.

The proof of Theorem \ref{prop:rh} is split into two cases: when $S_{\boldsymbol{\ep}}$ is path-connected and when it is not.

\begin{proof}[Proof of Theorem \ref{prop:rh} if $S_{\boldsymbol{\ep}}$ is path-connected] A~choice of $\boldsymbol{\ep}$ whose associated set $S_{\boldsymbol{\ep}}$ is path-connected exists by Corollary~\ref{cor:there_exists_choice}. Suppose that $S_{\boldsymbol{\ep}}$ is path-connected. The idea is to argue as in Pelayo, Ratiu and \vungoc\ \cite[Steps~2--4 in the proof of Theorem~B]{pvr_carto}, recalling and adjusting the argument therein as much as necessary for our purposes.

By construction $S_{\boldsymbol{\ep}} \subset B_{\mathrm{wt}}$. Now let $\mathsf{q} \colon \tilde{B}_{\mathrm{wt}} \to B_{\mathrm{wt}}$ denote the universal covering. By Corollary~\ref{cor:path_conn_iff_contr}, $S_{\boldsymbol{\ep}}$ is contractible and, in particular, simply connected. Therefore, there exists a~smooth section $\sigma \colon S_{\boldsymbol{\ep}} \to \tilde{B}_{\mathrm{wt}}$ of $\mathsf{q}$. Consider a developing map $\mathrm{dev} \colon \tilde{B}_{\mathrm{wt}} \to \R^2$ constructed by fixing basepoints $\mathbf{x}_0 \in S_{\boldsymbol{\ep}}$ and $\tilde{\mathbf{x}}_0 \in \sigma( S_{\boldsymbol{\ep}})$. Set $f_{\boldsymbol{\ep}} := \mathrm{dev} \circ \sigma \colon S_{\boldsymbol{\ep}} \to \R^2$. Arguing as in the proof of Lemma~\ref{lemma:vat_toric}, it is possible to choose $\mathrm{dev}$ so that $f_{\boldsymbol{\ep}}(x,y) = \big(x,f_{\boldsymbol{\ep}}^{(2)}(x,y)\big)$ for any $(x,y) \in S_{\boldsymbol{\ep}}$. Fix such a choice. Following the arguments in Pelayo, Ratiu and \vungoc\ \cite[Step~4 of the proof of Theorem~B]{pvr_carto}, $f_{\boldsymbol{\ep}}$ can be extended to an embedding $B \to \R^2$ which, by abuse of notation, is also denoted by $f_{\boldsymbol{\ep}}$. By construction and by density of $S_{\boldsymbol{\ep}} \subset B$, $(f_{\boldsymbol{\ep}},S_{\boldsymbol{\ep}})$ is a~cartographic pair with $f_{\boldsymbol{\ep}}(x,y) = \big(x,f^{(2)}_{\boldsymbol{\ep}}(x,y)\big)$. Thus for all $(x,y) \in S_{\boldsymbol{\ep}}$, $\frac{\partial
 f^{(2)}_{\boldsymbol{\ep}}}{\partial y} (x,y) \neq 0$. Since $S_{\boldsymbol{\ep}}$ is path-connected, property~\ref{item:8} follows.

To complete the proof, there are two cases to consider, depending on whether $\operatorname{sgn}(f_{\boldsymbol{\ep}}) = +1$ or $\operatorname{sgn}(f_{\boldsymbol{\ep}}) = -1$. In the first case, property~\ref{item:17} and local convexity of $f_{\boldsymbol{\ep}}(B)$ can be proved as in \vungoc~\cite[Steps~5 and~6 of the proof of Theorem~3.8]{vu-ngoc}. Thus suppose that $\operatorname{sgn}(f_{\boldsymbol{\ep}}) = -1$. Setting $\hat{f}_{\boldsymbol{\ep}} :=\left(\begin{smallmatrix} 1 & 0 \\ 0 & -1 \end{smallmatrix} \right) \circ f_{\boldsymbol{\ep}}$, $\big(\hat{f}_{\boldsymbol{\ep}},S_{\boldsymbol{\ep}}\big)$ is a cartographic pair which can be constructed as above satisfying $\operatorname{sgn}(\hat{f}_{\boldsymbol{\ep}}) = +1$. (This corresponds to adjusting the above choice of developing map by composing on the left with the map $\left(\left(\begin{smallmatrix} 1 & 0 \\ 0 & -1 \end{smallmatrix} \right), \left(\begin{smallmatrix} 0 \\ 0 \end{smallmatrix} \right)\right) \in \operatorname{Vert}(2;\Z)$; see Remark~\ref{rmk:freedom_no_ff}.) Fix $(\mathsf{x},\mathsf{y}) \in B_{\mathrm{wt}}$. Then, using property~\ref{item:17} for $\hat{f}_{\boldsymbol{\ep}}$ and the fact that $\operatorname{sgn}(f_{\boldsymbol{\ep}}) = -1$,{\samepage
 \begin{align*}
 \lim_{\stackrel{(x, y) \to (\mathsf{x},\mathsf{y})}{x <\mathsf{x}}} D f_{\boldsymbol{\ep}} (x, y) &=
 \begin{pmatrix}
 1 & 0 \\
 0 & -1
 \end{pmatrix}\lim_{\stackrel{(x, y) \to (\mathsf{x},\mathsf{y})}{x <\mathsf{x}}} D \hat{f}_{\boldsymbol{\ep}} (x, y) \\
 & = \begin{pmatrix}
 1 & 0 \\
 0 & -1
 \end{pmatrix}\begin{pmatrix} 1 & 0 \\ j_{\boldsymbol{\ep}}(\mathsf{x},\mathsf{y}) & 1
 \end{pmatrix}
 \lim_{\stackrel{(x, y) \to (\mathsf{x},\mathsf{y})}{x > \mathsf{x}}} D \hat{f}_{\boldsymbol{\ep}}(x, y) \\
 & = \begin{pmatrix}
 1 & 0 \\ \operatorname{sgn}(f_{\boldsymbol{\ep}}) j_{\boldsymbol{\ep}}(\mathsf{x},\mathsf{y}) & 1
 \end{pmatrix}\begin{pmatrix}
 1 & 0 \\ 0 & - 1
 \end{pmatrix}\lim_{\stackrel{(x, y) \to (\mathsf{x},\mathsf{y})}{x > \mathsf{x}}} D \hat{f}_{\boldsymbol{\ep}}(x,
 y) \\
 & =\begin{pmatrix} 1 & 0 \\ \operatorname{sgn}(f_{\boldsymbol{\ep}}) j_{\boldsymbol{\ep}}(\mathsf{x},\mathsf{y}) & 1
 \end{pmatrix}\lim_{\stackrel{(x, y) \to (\mathsf{x},\mathsf{y})}{x > \mathsf{x}}} D f_{\boldsymbol{\ep}}(x, y),
 \end{align*}
This proves property \ref{item:17} in general.}

Finally, observe that $\hat{f}_{\boldsymbol{\ep}}(B)$ is locally convex and that $\Z$-affine maps preserve this property. Thus $f_{\boldsymbol{\ep}}(B)$ is locally convex as required.
\end{proof}

Now we turn to the case of $S_{\boldsymbol{\ep}}$ not beingpath-connected. There exist proofs for such cases in the literature (see for instance Pelayo, Ratiu and \vungoc\ \cite[Step~5 of the proof of Theorem~B]{pvr_carto}). The argument presented below, however, uses different techniques.

Before delving into the proof, we introduce some useful notions and notation. For any \smash{$\mathsf{x} \in J(M)$}, set
\begin{gather*} 
 N_{\mathsf{x}} : = \lvert \{ i \in I \,|\, x_i = \mathsf{x} \} \rvert.
\end{gather*}
By property \ref{item:101} and Proposition \ref{lemma:iso_fixed_points}, $N_{\mathsf{x}}$ is finite for any $\mathsf{x} \in J(M)$. Moreover, for a fixed $\boldsymbol{\ep} \in \{ +1,-1\}^{I}$, set
\begin{gather*} 
 N_{\mathsf{x}}^{\pm}(\boldsymbol{\ep}) : = \pm \lvert \{ i \in I \,|\, x_i = \mathsf{x} \text{ and } \ep_i = \pm 1\} \rvert.
\end{gather*}
Observe that, for any $\mathsf{x} \in J(M)$ and any $\boldsymbol{\ep}\in \{ +1, -1\}^I$,
\begin{gather*}
 N_{\mathsf{x}} = N_{\mathsf{x}}^{+}(\boldsymbol{\ep}) - N_{\mathsf{x}}^{-}(\boldsymbol{\ep}).
\end{gather*}
Moreover, for any $(\mathsf{x},\mathsf{y}) \in B_{\mathrm{wt}}$ and any $\boldsymbol{\ep} \in \{ +1, -1\}^I$,
\begin{gather*}
 j_{\boldsymbol{\ep}}(\mathsf{x},\mathsf{y}) = N_{\mathsf{x}}^{+}(\boldsymbol{\ep}) + N_{\mathsf{x}}^{-}(\boldsymbol{\ep}).
\end{gather*}
Fix $\mathsf{x} \in J(M)$ with $N_{\mathsf{x}} \neq 0$. Then, by property~\ref{item:101} and Proposition~\ref{lemma:iso_fixed_points}, there exist finitely many indices $i_1 < i_2 < \cdots < i_{N_{\mathsf{x}}}$ in $I$ with $x_{i_j} = \mathsf{x}$. Observe that, by definition of the ordering on $B_{\mathrm{ff}}$,
\begin{gather*}
\{(x,y) \,|\, x = \mathsf{x} \} \cap B_{\mathrm{wt}} \subset \{(x,y) \,|\, x = \mathsf{x} \text{ and } y < y_{i_1} \} \\
\qquad {} \cup \bigcup\limits_{j=1}^{N_{\mathsf{x}} - 1} \{(x,y) \,|\, x = \mathsf{x} \text{ and } y_{i_j} < y < y_{i_{j+1}} \} \cup \{(x,y) \,|\, x = \mathsf{x} \text{ and } y > y_{i_{N_{\mathsf{x}}}}\}.
\end{gather*}

\begin{Lemma}\label{lem:function} For any $\boldsymbol{\ep} \in \{+1,-1\}^I$ and for all $\mathsf{x} \in J(M)$ with $N_{\mathsf{x}} \neq 0$, the function $j_{\boldsymbol{\ep}}(\mathsf{x},\cdot) \colon \{(x,y) \,|\, x =
 \mathsf{x} \}\cap B_{\mathrm{wt}} \to \Z$ satisfies
 \begin{itemize}[leftmargin=*]\itemsep=0pt
 \item $j_{\boldsymbol{\ep}}(\mathsf{x},\mathsf{y}) = N_{\mathsf{x}}^-(\boldsymbol{\ep})$ for all $(\mathsf{x},\mathsf{y}) \in \{(x,y) \,|\, x = \mathsf{x} \text{ and } y < y_{i_1} \}$;
 \item $j_{\boldsymbol{\ep}}(\mathsf{x},\mathsf{y}) = N_{\mathsf{x}}^-(\boldsymbol{\ep}) + k = N^+_{\mathsf{x}}(\boldsymbol{\ep}) - N_{\mathsf{x}} + j$ for all $k =1,\ldots, N_{\mathsf{x}}-1$ and for all $(\mathsf{x},\mathsf{y}) \in \{(x,y) \,|$ $x = \mathsf{x} \text{ and } y_{i_k} < y < y_{i_{k+1}} \}$;
 \item $j_{\boldsymbol{\ep}}(\mathsf{x},\mathsf{y}) = N_{\mathsf{x}}^+(\boldsymbol{\ep})$ for all $(\mathsf{x},\mathsf{y}) \in \{(x,y) \,|\, x = \mathsf{x} \text{ and } y > y_{i_{N_{\mathsf{x}}}} \}$.
 \end{itemize}
Hereby, $i_1<i_2 <\cdots < i_{N_{\mathsf{x}}}$ are the elements of $I$ such that $x_{i_j} = \mathsf{x}$. In particular, the function $j_{\boldsymbol{\ep}}(\mathsf{x},\cdot)$ only depends on $N_{\mathsf{x}}^{\pm}(\boldsymbol{\ep})$.
\end{Lemma}

\begin{proof} Fix $\boldsymbol{\ep} \in \{+1,-1\}^I$ and consider $(\mathsf{x},\mathsf{y}) \in B_{\mathrm{wt}}$ such that $N_{\mathsf{x}} \neq 0$. Suppose first that $\mathsf{y} < y_{i_1}$. This means that $(\mathsf{x},\mathsf{y})$ is `below' all focus-focus values on the vertical line $\{ (x,y) \,|\, x = \mathsf{x} \}$. By definition of the ordering on $B_{\mathrm{ff}}$ and of the cuts associated to $\boldsymbol{\ep}$, for all $k \in \{1,\ldots ,N_{\mathsf{x}}\}$, if $\ep_{i_k} = -1$ then $(\mathsf{x},\mathsf{y}) \in l^{\ep_{i_k}}$, while if $\ep_{i_k} = +1$, then $(\mathsf{x},\mathsf{y}) \notin l^{\ep_{i_k}}$. Thus
 \begin{gather*}
 j_{\boldsymbol{\ep}}(\mathsf{x},\mathsf{y}) = - \lvert \{ i \in I \,|\, x_i = \mathsf{x} \text{ and } \ep_i = - 1\} \rvert = N^{-}_{\mathsf{x}}(\boldsymbol{\ep})
 \end{gather*}
as required. Similarly, if $\mathsf{y} > y_{i_{N_{\mathsf{x}}}}$, then $j_{\boldsymbol{\ep}}(\mathsf{x},\mathsf{y}) = N^{+}_{\mathsf{x}}(\boldsymbol{\ep})$, for $(\mathsf{x},\mathsf{y})$ is `above' all focus-focus
 values on the vertical line $\{ (x,y) \,|\, x = \mathsf{x} \}$.

It remains to prove the intermediate cases for which we proceed by induction on $k$. The base case is $\mathsf{y} < y_{i_1}$, which has already been proved. Suppose that the required statement holds for all $m < k$ and let $y_{i_k} < \mathsf{y} < y_{i_{k+1}}$. Set, for any $(\bar{\mathsf{x}},\bar{\mathsf{y}}) \in B_{\mathrm{wt}}$,
 \begin{gather*}
 j^{\pm}_{\boldsymbol{\ep}}(\bar{\mathsf{x}},\bar{\mathsf{y}}) := \sum\limits_{\stackrel{i \in I\, , \,(\bar{\mathsf{x}},\bar{\mathsf{y}}) \in l^{\ep_i},}{ \ep_i = \pm 1}} \ep_i.
 \end{gather*}
Clearly, for any $(\bar{\mathsf{x}},\bar{\mathsf{y}}) \in B_{\mathrm{wt}}$, $j_{\boldsymbol{\ep}}(\bar{\mathsf{x}},\bar{\mathsf{y}}) = j_{\boldsymbol{\ep}}^+(\bar{\mathsf{x}},\bar{\mathsf{y}}) + j_{\boldsymbol{\ep}}^-(\bar{\mathsf{x}},\bar{\mathsf{y}})$. Fix some $(\mathsf{x},\mathsf{y}') \in B_{\mathrm{wt}}$ with $y_{i_{k-1}} < \mathsf{y}' < y_{i_{k}}$. The inductive hypothesis implies $j_{\boldsymbol{\ep}}(\mathsf{x},\mathsf{y}') = N_{\mathsf{x}}^-(\boldsymbol{\ep}) +k-1$. There are two cases to consider, depending on whether $\ep_{i_k} = +1$ or $\ep_{i_k} = -1$. In the former case, observe that $j^+_{\boldsymbol{\ep}}(\mathsf{x},\mathsf{y}) = j^+_{\boldsymbol{\ep}} (\mathsf{x},\mathsf{y}') + 1$, while $j^-_{\boldsymbol{\ep}}(\mathsf{x},\mathsf{y}) = j^-_{\boldsymbol{\ep}} (\mathsf{x},\mathsf{y}')$. Thus, $j_{\boldsymbol{\ep}}(\mathsf{x},\mathsf{y}) = j_{\boldsymbol{\ep}}(\mathsf{x},\mathsf{y}') + 1 = N^-_{\mathsf{x}}(\boldsymbol{\ep}) + k$ as required. The latter case is proved analogously, swapping the roles of $j^+_{\boldsymbol{\ep}}(\mathsf{x},\mathsf{y})$ and $j^-_{\boldsymbol{\ep}}(\mathsf{x},\mathsf{y})$.
\end{proof}

With the above results at hand, we finish the proof of Theorem \ref{prop:rh}.

\begin{proof}[Proof of Theorem \ref{prop:rh} if $S_{\boldsymbol{\ep}}$ is not path-connected] Suppose that $S_{\boldsymbol{\ep}}$ is not path-con\-nec\-ted. The idea is to reduce this situation to the path-connected case by appealing to the following result.

 \begin{Lemma}~\label{lem:reduce} There exists a unique $\hat{\boldsymbol{\ep}} \in \{+1,-1\}^I$ such that
 \begin{itemize}[leftmargin=*]\itemsep=0pt
 \item $S_{\boldsymbol{\ep}} \subset S_{\hat{\boldsymbol{\ep}}}$;
 \item for all $(\mathsf{x},\mathsf{y}) \in B_{\mathrm{wt}}$, $j_{\boldsymbol{\ep}}(\mathsf{x},\mathsf{y}) = j_{\hat{\boldsymbol{\ep}}}(\mathsf{x},\mathsf{y})$;
 \item $S_{\hat{\boldsymbol{\ep}}}$ is path-connected.
 \end{itemize}
 \end{Lemma}

Assume Lemma~\ref{lem:reduce}, whose proof is below, and let $\hat{\boldsymbol{\ep}}$ be as in Lemma~\ref{lem:reduce}. Let $f_{\hat{\boldsymbol{\ep}}} \colon B \to \R^2$ be a cartographic homeomorphism associated to $\hat{\boldsymbol{\ep}}$. Set $f_{\boldsymbol{\ep}}:= f_{\hat{\boldsymbol{\ep}}}$. Since $S_{\boldsymbol{\ep}} \subset S_{\hat{\boldsymbol{\ep}}}$, Corollary~\ref{cor:complement_open} implies that $(f_{\boldsymbol{\ep}},S_{\boldsymbol{\ep}})$ is a~cartographic pair; property~\ref{item:8} holds by construction. Property~\ref{item:17} holds because Lemma~\ref{lem:reduce} implies that $j_{\boldsymbol{\ep}}(\mathsf{x},\mathsf{y}) = j_{\hat{\boldsymbol{\ep}}}(\mathsf{x},\mathsf{y})$ for all $(\mathsf{x},\mathsf{y}) \in B_{\mathrm{wt}}$. Moreover, $\operatorname{sgn}(f_{\boldsymbol{\ep}}) = \operatorname{sgn}(f_{\hat{\boldsymbol{\ep}}})$ holds by definition. Local convexity of $f_{\boldsymbol{\ep}}(B) = f_{\hat{\boldsymbol{\ep}}}(B)$ is also true as $f_{\hat{\boldsymbol{\ep}}}$ is associated to $\hat{\boldsymbol{\ep}}$ in the sense of Proposition~\ref{prop:rh}. This finishes the proof of Theorem~\ref{prop:rh} for the case that~$S_{\boldsymbol{\ep}}$ is not path-connected.
\end{proof}

\begin{proof}[Proof of Lemma~\ref{lem:reduce}] For a fixed $x_i$, the map $y \mapsto j_{\boldsymbol{\epsilon}}(x_i,y)$ is an integer-valued map whose image is all integers between (and including) $N^-_{x_i}(\boldsymbol{\ep})$ and $N^+_{x_i}(\boldsymbol{\ep})$. Moreover, it is easily seen that the above map is constant along the open segments lying between consecutive focus-focus values. Therefore there exists a segment on the vertical line $\{(x,y) \,|\, x = x_i\}$ along which $j_{\boldsymbol{\epsilon}}(x_i ,\cdot)$ is equal to zero. Consider $\bar{\boldsymbol{\epsilon}} = \{+1,-1\}^I$ that is equal to $\boldsymbol{\ep}$ except for those indices corresponding to focus-focus values lying on $\{(x,y) \,|\, x = x_i\}$; for those indices, set $\bar{\epsilon}_j = \pm1$ according to whether~$c_j$ lies above or below the segment on which $j_{\boldsymbol{\epsilon}}(x_i ,\cdot)$ is equal to zero. It can be checked that, for all $(\mathsf{x},\mathsf{y}) \in B_{\mathrm{wt}}$, $j_{\boldsymbol{\epsilon}}(\mathsf{x},\mathsf{y}) = j_{\bar{\boldsymbol{\epsilon}}}(\mathsf{x},\mathsf{y})$. Performing this operation for all $x_i$ yields $\hat{\boldsymbol{\ep}}$ satisfying the requirements of Lemma~\ref{lemma:S_ep_path_conn} for $S_{\hat{\boldsymbol{\ep}}}$ to be path-connected. Moreover, by construction, $S_{\boldsymbol{\ep}} \subset
 S_{\hat{\boldsymbol{\ep}}}$ and it can be easily checked that $\hat{\boldsymbol{\ep}}$ is the unique choice of signs that satisfies all the above properties.
 \end{proof}

\begin{figure}[h] \centering
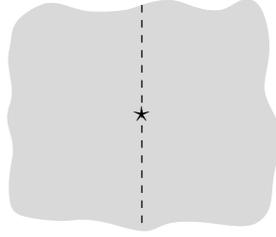
 \caption{The focus-focus value $c$ (displayed by $\star$) has multiplicity two and the two signs associated to $c$ are $+1$ and $-1$.} \label{LABEL8}
\end{figure}

\begin{rmk}\label{rmk:not_higher_mult} The above argument for the case of $S_{\boldsymbol{\ep}}$ not being path-connected only works because the system $\vat$ satisfies property~\ref{item:7}. If focus-focus values of higher multiplicity are allowed, then there may be no analogue of $\hat{\boldsymbol{\ep}}$ as in Lemma~\ref{lem:reduce}. In this case, the issue is that there is a choice of $\boldsymbol{\epsilon}$ such that there exists an~$x_i$ for which there is no interval on which the function $y \mapsto j_{\boldsymbol{\epsilon}}(x_i,y)$ is zero (see Fig.~\ref{LABEL8}).
\end{rmk}

\begin{rmk}\label{rmk:connected_suffices} The argument proving the case of Theorem~\ref{prop:rh} when $S_{\boldsymbol{\ep}}$ is not path-connected shows that, for a given faithful semitoric system, the set of all cartographic homeomorphisms that satisfy the hypotheses of Theorem~\ref{prop:rh} equals that of the cartographic homeomorphisms associated to those $\boldsymbol{\ep}$ for which $S_{\boldsymbol{\ep}}$ is path-connected (cf.\ \vungoc\ \cite[Proposition~4.1]{vu-ngoc}).
\end{rmk}

To finish this section, we prove Proposition~\ref{prop:im_equ} using Theorem~\ref{prop:rh}.

\begin{proof}[Proof of Proposition~\ref{prop:im_equ}] Let $(M_1,\omega_1,\Phi_1)$ be a faithful semitoric system with at least one focus-focus point. If $(M_2,\omega_2,\Phi_2)$ is a faithful semitoric system isomorphic to $(M_1,\omega_1,\Phi_1)$ as an integrable system, then the cardinality of the set of focus-focus values of the former equals that of the set of focus-focus values of the latter, as the isomorphism induces a~bijection between the focus-focus values. In order to prove the result, the non-trivial implication to check is that if $(M_1,\omega_1,\Phi_1 = (J_1,H_1))$ and $(M_2,\omega_2,\Phi_2 = (J_2,H_2))$ are isomorphic as integrable systems, then they are isomorphic as complexity one systems. In fact, we aim to achieve more: we show that any isomorphism $(\Psi,\psi)$ as integrable systems between $(M_1,\omega_1,\Phi_1 = (J_1,H_1))$ and $(M_2,\omega_2,\Phi_2 = (J_2,H_2))$ is necessarily an isomorphism as complexity one systems.

Fix an isomorphism $\big(\Psi,\psi = \big(\psi^{(1)},\psi^{(2)}\big)\big)$ as integrable systems between $(M_1,\omega_1,\Phi_1 = (J_1,H_1))$ and $(M_2,\omega_2,\Phi_2 = (J_2,H_2))$, let~$I$ denote the indexing set of the set of focus-focus values $B_{1, \mathrm{ff}}$ and set $B_{1,\mathrm{ff}} = \{c_i\}_{i \in I}$. Begin by observing that the function $\psi^{(1)} \circ \Phi_1 \colon M_1 \to \R$ is the moment map of an effective Hamiltonian $S^1$-action that is system-preserving. This is because $\big(\Psi,\psi = \big(\psi^{(1)},\psi^{(2)}\big)\big)$ is an isomorphism of integrable systems and because~$J_2$ is the moment map of an effective Hamiltonian $S^1$-action that preserves the integrable system $(M_2,\omega_2,\Phi_2)$. Fix a~choice of signs $\boldsymbol{\ep} \in \{+1,-1\}^I$ so that $S_{\boldsymbol{\ep}} \subset B_1$ is path-connected (Corollary~\ref{cor:there_exists_choice} ensures that such a choice exists), and fix a cartographic homeomorphism $f_{\boldsymbol{\ep}}\colon B_1 \to \R^2$ associated to $\boldsymbol{\ep}$ as in Theorem~\ref{prop:rh}. The above observation implies that the map $\psi^{(1)} \circ f^{-1}_{\boldsymbol{\ep}}|_{f_{\boldsymbol{\ep}}(S_{\boldsymbol{\ep}})}$ is smooth and $\Z$-affine. Since $f_{\boldsymbol{\ep}}(S_{\boldsymbol{\ep}})$ is connected and endowed with the restriction of the standard $\Z$-affine structure on $\R^2$, it follows that there exist $\zeta_1,\zeta_2 \in \Z$ and $a \in \R$ such that, for all $(x,y) \in f_{\boldsymbol{\ep}}(S_{\boldsymbol{\ep}})$, $\psi^{(1)} \circ f^{-1}_{\boldsymbol{\ep}}(x,y) = \zeta_1 x + \zeta_2 y + a$. As $f_{\boldsymbol{\ep}}(S_{\boldsymbol{\ep}}) \subset f_{\boldsymbol{\ep}}(B_1)$ is dense, it follows that, for all $(x,y) \in f_{\boldsymbol{\ep}}(B_1)$, $\psi^{(1)} \circ f^{-1}_{\boldsymbol{\ep}}(x,y) = \zeta_1 x + \zeta_2 y + a$. Therefore, for all $(x,y) \in B_1$,
 \begin{gather} \label{eq:18}
 \psi^{(1)}(x,y) = \zeta_1 x + \zeta_2 f^{(2)}_{\boldsymbol{\ep}}(x,y) + a,
 \end{gather}
since $f_{\boldsymbol{\ep}}(x,y) = \big(x,f^{(2)}_{\boldsymbol{\ep}}(x,y)\big)$. To prove the result, it suffices to show that $\zeta_1 \in \{+1,-1\}$ and that $\zeta_2 = 0$. To this end, fix $i \in I$ and let $U_i$ be an open neighborhood of $c_i$ in $B_1$ with the property that the subsystem of $(M_1,\omega_1,\Phi_1 = (J_1,H_1))$ relative to $U_i$ is a saturated regular neighborhood of the focus-focus fiber $\Phi^{-1}_1(c_i)$, i.e., all fibers in $\Phi^{-1}_1(U_i)$ are regular except for the focus-focus fiber $\Phi^{-1}_1(c_i)$ (see Section~\ref{sec:at_generalities}). Since $\psi^{(1)} \circ \Phi_1|_{\Phi_1^{-1}(U_i)}$ is the moment map of an effective Hamiltonian $S^1$-action that is system-preserving, Proposition~\ref{prop:S1-action} implies that, on~$\Phi^{-1}_1(U_i)$, $\pm d J_1 = d \big(\psi^{(1)} \circ \Phi_1\big)$. Since~$\Phi_1$ is a submersion away from the focus-focus point~$p_i$ lying on~$\Phi^{-1}_1(c_i)$ and $\Phi_1\big(\Phi^{-1}_1(U_i) \smallsetminus \{p_i\}\big) = U_i$, it follows that, on $U_i$,
 \begin{gather} \label{eq:19}
 d \psi^{(1)} = \pm d x.
 \end{gather}
 On the other hand, formula \eqref{eq:18} implies that, on $B_1$,
 \begin{gather} \label{eq:20}
 d \psi^{(1)} = \left(\zeta_1 + \zeta_2 \frac{\partial f^{(2)}}{\partial x}\right) dx + \zeta_2 \frac{\partial f^{(2)}}{\partial y}dy.
 \end{gather}
Comparing \eqref{eq:19} and \eqref{eq:20} and using the fact that $\zeta_1$, $\zeta_2$ are constant and that $\frac{\partial f^{(2)}}{\partial y} \neq 0$ on~$B_1$, we have that $\zeta_2 = 0$ and $\zeta_1 =\pm 1$, as desired.
 \end{proof}

\subsection{The set of cartographic pairs}\label{sec:set-cart-home}
Given a faithful semitoric system $\vat$, it is natural to ask for a description of the set of all cartographic pairs of~$\vat$ that satisfy the properties in Theorem~\ref{prop:rh}. These cartographic pairs and their images are henceforth referred to as {\em semitoric}. One reason why it is important to describe the set of semitoric cartographic images of a~faithful semitoric system with at least one focus-focus point is that it is an invariant of the isomorphism class of the system as a~strict complexity one system (see Definition~\ref{defn:im_faithful_st}).

\begin{Lemma}\label{lemma:im_invariant} Let $(M_1,\omega_1,\Phi_1)$ and $(M_2,\omega_2,\Phi_2)$ be faithful semitoric systems that are isomorphic as strict complexity one spaces and suppose that $(M_1,\omega_1,\Phi_1)$ has at least one focus-focus point. Then the sets of semitoric cartographic images of $(M_1,\omega_1,\Phi_1)$ and $(M_2,\omega_2,\Phi_2)$ are equal.
 \end{Lemma}
\begin{proof} Let $(\Psi,\psi)$ be an isomorphism as strict complexity one spaces between faithful semitoric systems $(M_1,\omega_1,\Phi_1)$ and $(M_2,\omega_2,\Phi_2)$, i.e., $\psi(x,y) = \big(x,\psi^{(2)}(x,y)\big)$ for some smooth function $\psi^{(2)}\colon B_1 \to \R$. First, observe that, since $(\Psi,\psi)$ is, in particular, an isomorphism as integrable systems, it induces a bijection between the sets of focus-focus values of $(M_1,\omega_1,\Phi_1)$ and $(M_2,\omega_2,\Phi_2)$. In what follows, this bijection is used tacitly to identify the indexing sets~$I_1$ and~$I_2$ of the sets of focus-focus values of $(M_1,\omega_1,\Phi_1)$ and $(M_2,\omega_2,\Phi_2)$. Suppose that $(f_{\boldsymbol{\ep}_2},S_{\boldsymbol{\ep}_2})$ is a~semitoric cartographic pair for $(M_2,\omega_2,\Phi_2)$. Corollary~\ref{cor:pull_back_iso_carto} implies that $\big(f_{\boldsymbol{\ep}_2} \circ \psi, \psi^{-1}(S_{\boldsymbol{\ep}_2})\big)$ is a~cartographic pair for $(M_1,\omega_1,\Phi_1)$. The aim is to show that this pair is semitoric, i.e., that
 \begin{enumerate}[label = (\arabic*), ref = (\arabic*), leftmargin=*]\itemsep=0pt
 \item \label{item:28} there exists a choice of signs $\boldsymbol{\ep}_1$ for $(M_1,\omega_1,\Phi_1)$ such that $S_{\boldsymbol{\ep}_1} =
 \psi^{-1}(S_{\boldsymbol{\ep}_2})$, and
 \item \label{item:29} the map $f_{\boldsymbol{\ep}_1}:= f_{\boldsymbol{\ep}_2} \circ\psi$ is associated to $\boldsymbol{\ep}_1$, i.e., as per the hypotheses of Theorem~\ref{prop:rh}, $f_{\boldsymbol{\ep}_1}$ is of the form $f_{\boldsymbol{\ep}_1}(x,y) = \big(x,f^{(2)}_{\boldsymbol{\ep}_1}(x,y)\big)$ for some continuous function $f^{(2)}_{\boldsymbol{\ep}_1} \colon B_1 \to \R$, and it satisfies properties \ref{item:8} and \ref{item:17}.
 \end{enumerate}
Begin by observing that, since $\psi(x,y) = \big(x,\psi^{(2)}(x,y)\big)$, the preimage of a vertical line under $\psi$ is a vertical line. This implies that the preimage of the vertical cuts determined by $\boldsymbol{\ep}_2$ are vertical half-lines. Moreover, since $(\Psi,\psi)$ is an isomorphism of integrable systems, the `starting point' of the above vertical lines are focus-focus values of $(M_1,\omega_1,\Phi_1)$. This determines uniquely a~choice of signs $\boldsymbol{\ep}_1$ for $(M_1,\omega_1,\Phi_1)$ such that $S_{\boldsymbol{\ep}_1} = \psi^{-1}(S_{\boldsymbol{\ep}_2})$, thus proving~\ref{item:28}. To prove~\ref{item:29}, observe that, by definition, $f_{\boldsymbol{\ep}_1} = f_{\boldsymbol{\ep}_2} \circ \psi$. Since $\psi(x,y) = \big(x,\psi^{(2)}(x,y)\big)$, the fact that $(f_{\boldsymbol{\ep}_2},S_{\boldsymbol{\ep}_2})$ satisfies the properties of Theorem~\ref{prop:rh} implies that
 \begin{gather*}
 f_{\boldsymbol{\ep}_1}(x,y) = \big(x, f^{(2)}_{\boldsymbol{\ep}_2}\big(x,\psi^{(2)}(x,y)\big)\big) =: \big(x,f^{(2)}_{\boldsymbol{\ep}_1}(x,y)\big).
 \end{gather*}
 Using the chain rule, for any $(x,y) \in S_{\boldsymbol{\ep}_1}$,
 \begin{gather} \label{eq:21}
 \frac{\partial f^{(2)}_{\boldsymbol{\ep}_1}}{\partial y}(x,y) = \frac{\partial f^{(2)}_{\boldsymbol{\ep}_2}}{\partial z}(x,z)\bigg\rvert_{z=\psi^{(2)}(x,y)}\frac{\partial \psi^{(2)}}{\partial y}(x,y).
 \end{gather}
Observe that, by construction, $(x,y) \in S_{\boldsymbol{\ep}_1}$ if and only if $(x,\psi^{(2)}(x,y)) \in S_{\boldsymbol{\ep}_2}$. Moreover, since $\psi(x,y) = \big(x,\psi^{(2)}(x,y)\big)$ is a diffeomorphism and $B_1$ is connected, the quantity $\operatorname{sgn}\big(\frac{\partial \psi^{(2)}}{\partial y}(x,y)\big)$ $=: \operatorname{sgn} (\psi )$ does not depend on $(x,y) \in B_1$. Using the fact that $f_{\boldsymbol{\ep}_2}$ satisfies property~\ref{item:8}, equation \eqref{eq:21} yields that $\operatorname{sgn}\big(\frac{\partial f^{(2)}_{\boldsymbol{\ep}_1}}{\partial y}(x,y)\big) =: \operatorname{sgn}(f_{\boldsymbol{\ep}_1})$ does not depend on $(x,y) \in S_{\boldsymbol{\ep}_1}$. In fact, $\operatorname{sgn}(f_{\boldsymbol{\ep}_1}) = \operatorname{sgn}(f_{\boldsymbol{\ep}_2})\operatorname{sgn}(\psi)$. This shows that $f_{\boldsymbol{\ep}_1}$ satisfies property~\ref{item:8}. To see that it satisfies property~\ref{item:17}, observe that, for any $(\mathsf{x},\mathsf{y}) \in B_{1,\mathrm{wt}}$,
 \begin{gather*}
 j_{\boldsymbol{\ep}_1}(\mathsf{x},\mathsf{y}) = j_{\boldsymbol{\ep}_2}(\psi(\mathsf{x},\mathsf{y})) \operatorname{sgn}(\psi).
 \end{gather*}
Using the above formula, the chain rule and the fact that $f_{\boldsymbol{\ep}_2}$ satisfies property~\ref{item:17}, it can be shown that $f_{\boldsymbol{\ep}_1}$ satisfies property~\ref{item:17} as well. Thus $\big(f_{\boldsymbol{\ep}_2} \circ \psi, \psi^{-1}(S_{\boldsymbol{\ep}_2})\big)$ is a~semitoric cartographic pair for $(M_1,\omega_1,\Phi_1)$, as desired. Hence, the set of semitoric cartographic images of $(M_2,\omega_2,\Phi_2)$ is contained in that of $(M_1,\omega_1,\Phi_1)$. Reversing the roles of $(M_1,\omega_1,\Phi_1)$ and $(M_2,\omega_2,\Phi_2)$ completes the proof.
\end{proof}

Providing a description of the set of semitoric cartographic images of a faithful semitoric system with at least one focus-focus point is the aim of this section, which generalizes, while being heavily inspired by, work of \vungoc\ \cite[Section~4]{vu-ngoc}. Henceforth, fix a faithful semitoric system~$\vat$ with at least one focus-focus point. The idea is to show that any cartographic homeomorphism of~$\vat$ as in Theorem~\ref{prop:rh} can be constructed from a fixed one by means of composing on the left with a suitable homeomorphism (see Corollary~\ref{cor:freedom_ff_sc} and Theorem~\ref{thm:different_carto} for precise statements). It is convenient to consider two separate cases:
\begin{itemize}[leftmargin=*]\itemsep=0pt
 \item Determine all cartographic homeomorphisms associated to a given choice of signs (Corollary~\ref{cor:freedom_ff_sc}).
 \item Determine how cartographic homeomorphisms associated to possibly distinct choices of signs are related (Theorem~\ref{thm:different_carto}).
\end{itemize}

First, we consider the set of all cartographic homeomorphisms associated to a given choice of signs; this is described in the following result, which is analogous to Remark~\ref{rmk:freedom_no_ff}. Recall that $\operatorname{Vert}(2;\Z)$ denotes the subgroup of~$\mathrm{AGL}(2;\Z)$ that preserves vertical lines (see Remark~\ref{rmk:toric_im}).

\begin{Corollary}\label{cor:freedom_ff_sc} Fix a faithful semitoric system $\vat$ and a choice of signs~$\boldsymbol{\epsilon}$. If $f_{\boldsymbol{\ep}}, \hat{f}_{\boldsymbol{\ep}}$ are cartographic homeomorphisms associated to $\boldsymbol{\ep}$, then there exists an element $h \in \operatorname{Vert}(2;\Z)$ with $\hat{f}_{\boldsymbol{\ep}} = h \circ f_{\boldsymbol{\ep}}$. Conversely, for any $h \in \operatorname{Vert}(2;\Z)$, the map $\hat{f}_{\boldsymbol{\ep}}:= h \circ f_{\boldsymbol{\ep}}$ is a~cartographic homeomorphism associated to $\boldsymbol{\ep}$.
\end{Corollary}

\begin{proof} It may be assumed without loss of generality that $S_{\boldsymbol{\ep}}$ is path-connected since the not path-connected case can be reduced to the path-connected one as in the proof of Theorem~\ref{prop:rh} (see Remark~\ref{rmk:connected_suffices}). Let $f_{\boldsymbol{\ep}}, \hat{f}_{\boldsymbol{\ep}} \colon B \to \R^2$ be cartographic homeomorphisms associated to~$\boldsymbol{\ep}$. Then their restrictions to $S_{\boldsymbol{\ep}}$ are developing maps for the induced $\Z$-affine structure on $S_{\boldsymbol{\ep}}$. Therefore, arguing as in Remark~\ref{rmk:freedom_no_ff}, there exists an element $h \in \operatorname{Vert}(2;\Z)$ with $\hat{f}_{\boldsymbol{\ep}}|_{S_{\boldsymbol{\ep}}} = h \circ f _{\boldsymbol{\ep}}|_{S_{\boldsymbol{\ep}}}$. Since $S_{\boldsymbol{\ep}} \subset B$ is dense, this implies that $\hat{f}_{\boldsymbol{\ep}} = h \circ f_{\boldsymbol{\ep}}$ as required.

Conversely, the proof of Theorem \ref{prop:rh} gives that composing a cartographic homeomorphism associated to $\boldsymbol{\ep}$ on the left with an element of $\operatorname{Vert}(2;\Z)$ yields another cartographic homeomorphism associated to $\boldsymbol{\ep}$.
\end{proof}

Having established Corollary \ref{cor:freedom_ff_sc}, we study the problem of relating cartographic homeomorphisms whose associated signs are not necessarily equal. Before stating the main result of this section we introduce some tools akin to those needed in \vungoc~\cite[Section~4]{vu-ngoc}, but slightly more involved as faithful semitoric systems allow for the presence of infinitely many focus-focus points (see Remarks~\ref{rmk:difference_san},~\ref{rmk:maps_r} and~\ref{rmk:compare_st} below).

As in Section \ref{sec:cuts-cart-home}, let $I$ be the indexing set of the set of focus-focus values $B_{\mathrm{ff}}$ defined in equation~\eqref{eq:2}. Also, fix the ordering on~$B_{\mathrm{ff}}$ as in the paragraph following equation \eqref{eq:2}, so elements of~$B_{\mathrm{ff}}$ are denoted by $c_i = (x_i,y_i)$ for $i \in I$. Henceforth, fix elements $\boldsymbol{\ep}, \hat{\boldsymbol{\ep}} \in \{+1,-1\}^I$ and associated cartographic homeomorphisms $f_{\boldsymbol{\ep}}, f_{\hat{\boldsymbol{\ep}}} \colon B \to \R^2$. Furthermore, fix a~basepoint $(\mathsf{x},\mathsf{y}) \in B_{\mathrm{wt}}$ with the property that $x_0 < \mathsf{x} < x_1$. (If $x_0$ or $ x_1$ is not defined, then only the other inequality is required.)

\begin{rmk}\label{rmk:difference_san} The above choice of basepoint agrees with that made in \vungoc\ \cite[proof of Proposition~4.1]{vu-ngoc}. However, there is an important difference that arises because of the possibility of having infinitely many focus-focus points for faithful semitoric systems. In what follows we must allow for the case in which there are focus-focus values `to the left' of the basepoint, i.e., with notation as above, for the case in which there exists $i \in I$ with $x_i < \mathsf{x}$. (If this is the case, then by the choices of indexing set $I$ of equation~\eqref{eq:2} and of basepoint, there are infinitely many such indices.)
\end{rmk}

Throughout this section, set $ T:= \left(\begin{smallmatrix} 1& 0 \\ 1 & 1 \end{smallmatrix}\right)$. For any $i \in I$, set
\begin{gather*} 
 k_i(\boldsymbol{\ep},\hat{\boldsymbol{\ep}}):= \operatorname{sgn}(f_{\hat{\boldsymbol{\ep}}}) \left(\frac{\ep_i - \hat{\ep}_i}{2}\right).
\end{gather*}
Moreover, for any $i \in I$ define $\mathfrak{l}_{i,\boldsymbol{\ep},\hat{\boldsymbol{\ep}}}\colon \R^2 \to \R^2$ as follows:
\begin{itemize}[leftmargin=*]\itemsep=0pt
\item If $i \leq 0$, let $\mathfrak{l}_{i,\boldsymbol{\ep},\hat{\boldsymbol{\ep}}}$ be the identity.
\item If $i > 0$, let $\mathfrak{l}_{i,\boldsymbol{\ep},\hat{\boldsymbol{\ep}}}$ be the piece-wise $\Z$-affine transformation that acts as the identity on the half-space $x < x_i$ and as the shear $T^{k_i(\boldsymbol{\ep},\hat{\boldsymbol{\ep}})}$ on the half-space $x \geq x_i$.
\end{itemize}
Analogously, for any $i \in I$ define $\mathfrak{r}_{i,\boldsymbol{\ep},\hat{\boldsymbol{\ep}}} \colon \R^2 \to \R^2$ as follows:
\begin{itemize}[leftmargin=*]\itemsep=0pt
\item If $i \leq 0$, let $\mathfrak{r}_{i,\boldsymbol{\ep},\hat{\boldsymbol{\ep}}}$ be the piece-wise $\Z$-affine transformation that acts as the shear $T^{-k_i(\boldsymbol{\ep},\hat{\boldsymbol{\ep}})}$ on the half-space $x < x_i$ and as the identity on the half-space $x \geq x_i$.
\item If $i > 0$ let $\mathfrak{r}_{i,\boldsymbol{\ep},\hat{\boldsymbol{\ep}}}$ be the identity.
\end{itemize}

\begin{rmk}\label{rmk:maps_r} While the maps $\mathfrak{l}_{i,\boldsymbol{\ep},\hat{\boldsymbol{\ep}}}$ are those used in \vungoc\ \cite[Section~4]{vu-ngoc}, the maps $\mathfrak{r}_{i,\boldsymbol{\ep},\hat{\boldsymbol{\ep}}}$ are needed in the following precisely because of the possibility of focus-focus values existing `to the left' of the basepoint (see Remark~\ref{rmk:difference_san}).
\end{rmk}

Explicitly, we have that if $i > 0$ then
\begin{gather*} 
 \mathfrak{l}_{i,\boldsymbol{\ep},\hat{\boldsymbol{\ep}}}:=
 \begin{cases}
 \mathrm{id} & \text{if } x < x_i, \\
 \left( \begin{pmatrix}
 1 & 0 \\
 k_i(\boldsymbol{\ep},\hat{\boldsymbol{\ep}}) & 1
 \end{pmatrix},
 \begin{pmatrix}
 0 \\
 - k_i(\boldsymbol{\ep},\hat{\boldsymbol{\ep}}) x_i
 \end{pmatrix}\right) & \text{if } x \geq x_i,
 \end{cases}
\end{gather*}
and if $i \leq 0$ then
\begin{gather*} 
 \mathfrak{r}_{i,\boldsymbol{\ep},\hat{\boldsymbol{\ep}}}=
 \begin{cases}
 \left( \begin{pmatrix}
 1 & 0 \\
 - k_i(\boldsymbol{\ep},\hat{\boldsymbol{\ep}})& 1
 \end{pmatrix},
 \begin{pmatrix}
 0 \\
 k_i(\boldsymbol{\ep},\hat{\boldsymbol{\ep}}) x_i
 \end{pmatrix}\right) & \text{if } x \leq x_i, \\
 \mathrm{id} & \text{if } x > x_i.
 \end{cases}
\end{gather*}

The maps $\mathfrak{l}_{i,\boldsymbol{\ep},\hat{\boldsymbol{\ep}}}$ and $\mathfrak{r}_{i,\boldsymbol{\ep},\hat{\boldsymbol{\ep}}}$ are well-defined and satisfy the following properties, the proofs of which are left to the reader.

\begin{Lemma}\label{cor:piece_z-aff_homeo} For each positive $($respectively non-positive$)$ $i\in I$, the map $ \mathfrak{l}_{i,\boldsymbol{\ep},\hat{\boldsymbol{\ep}}}$ $($respecti\-ve\-ly~$\mathfrak{r}_{i,\boldsymbol{\ep},\hat{\boldsymbol{\ep}}})$ is a~homeomorphism that fixes the vertical line $\{(x,y) \,|\, x = x_i\}$ pointwise and is a~$\Z$-affine isomorphism of $\R^2 \smallsetminus \{(x,y) \,|\, x = x_i\}$.
\end{Lemma}

\begin{Lemma}\label{cor:commute_piece} For any $i,j\in I$,
\begin{gather*}
 \mathfrak{l}_{i,\boldsymbol{\ep},\hat{\boldsymbol{\ep}}} \circ\mathfrak{l}_{j,\boldsymbol{\ep},\hat{\boldsymbol{\ep}}}
 = \mathfrak{l}_{j,\boldsymbol{\ep},\hat{\boldsymbol{\ep}}} \circ \mathfrak{l}_{i,\boldsymbol{\ep},\hat{\boldsymbol{\ep}}},\qquad
 \mathfrak{r}_{i,\boldsymbol{\ep},\hat{\boldsymbol{\ep}}} \circ\mathfrak{r}_{j,\boldsymbol{\ep},\hat{\boldsymbol{\ep}}} =
 \mathfrak{r}_{j,\boldsymbol{\ep},\hat{\boldsymbol{\ep}}} \circ \mathfrak{r}_{i,\boldsymbol{\ep},\hat{\boldsymbol{\ep}}},\qquad
 \mathfrak{l}_{i,\boldsymbol{\ep},\hat{\boldsymbol{\ep}}} \circ \mathfrak{r}_{j,\boldsymbol{\ep},\hat{\boldsymbol{\ep}}} = \mathfrak{r}_{j,\boldsymbol{\ep},\hat{\boldsymbol{\ep}}} \circ \mathfrak{l}_{i,\boldsymbol{\ep},\hat{\boldsymbol{\ep}}}.
\end{gather*}
\end{Lemma}

One ingredient in the construction of the homeomorphism that relates $f_{\boldsymbol{\ep}}$ and $f_{\hat{\boldsymbol{\ep}}}$ is the composition of the maps $\mathfrak{l}_{i,\boldsymbol{\ep},\hat{\boldsymbol{\ep}}}$ or $\mathfrak{r}_{i,\boldsymbol{\ep},\hat{\boldsymbol{\ep}}}$ as $i$ ranges over all indices in~$I$. Seeing as this may involve the composition of infinitely many maps different from the identity (as~$I$ may be infinite), some care is needed. To this end, we first introduce notation for the domains of these possibly infinite compositions. Set
\begin{gather*}
 D_{\sup} :=
 \begin{cases}
 \R^2 & \text{if } \mathsf{x}_{\sup} \in J(M), \\
 \big\{ (x,y) \in \R^2 \,|\, x < \mathsf{x}_{\sup} \big\} & \text{otherwise,}
 \end{cases} \\
 D_{\inf} :=
 \begin{cases}
 \R^2 & \text{if } \mathsf{x}_{\inf} \in J(M), \\
 \big\{ (x,y) \in \R^2 \,|\, x > \mathsf{x}_{\inf} \big\} & \text{otherwise,}
 \end{cases}
\end{gather*}
where $\mathsf{x}_{\sup}, \mathsf{x}_{\inf} \in \R$ are defined as in equation~\eqref{eq:9}. Observe that $B \subset D_{\sup} \cap D_{\inf}$.

Next, we define the desired compositions. If the cardinality of $I$ is finite, the situation is entirely analogous to the one considered in \vungoc\ \cite[Section~4]{vu-ngoc}. For, with the above conventions, $I$ being finite implies that $\mathsf{x}_{\sup}, \mathsf{x}_{\inf} \in J(M)$ and that~$I$ only contains positive elements. In this case, set $\mathfrak{l}_{\boldsymbol{\ep},\hat{\boldsymbol{\ep}}}\colon
D_{\sup} = \R^2 \to \R^2$ and $\mathfrak{r}_{\boldsymbol{\ep},\hat{\boldsymbol{\ep}}}\colon D_{\inf} = \R^2 \to \R^2$ to be equal to the finite compositions $\mathfrak{l}_{\lvert I \rvert, \boldsymbol{\ep},\boldsymbol{\hat{\ep}}} \circ \mathfrak{l}_{\lvert I \rvert -1, \boldsymbol{\ep},\boldsymbol{\hat{\ep}}} \circ \cdots \circ \mathfrak{l}_{1, \boldsymbol{\ep},\boldsymbol{\hat{\ep}}}$ and $\mathfrak{r}_{\lvert I \rvert, \boldsymbol{\ep},\boldsymbol{\hat{\ep}}} \circ \mathfrak{r}_{\lvert I \rvert -1, \boldsymbol{\ep},\boldsymbol{\hat{\ep}}} \circ \cdots \circ \mathfrak{r}_{1, \boldsymbol{\ep},\boldsymbol{\hat{\ep}}}$ respectively. Observe that the latter is, by definition, equal to the identity. It remains to consider the case in which $I$ is infinite. If $I$ has infinitely many positive elements, set, for any $(x,y)$ with $x \leq x_i$,
\begin{gather*}
 \mathfrak{l}_{\boldsymbol{\ep},\hat{\boldsymbol{\ep}}}(x,y):=\mathfrak{l}_{i, \boldsymbol{\ep},\boldsymbol{\hat{\ep}}} \circ\mathfrak{l}_{i -1, \boldsymbol{\ep},\boldsymbol{\hat{\ep}}} \circ \cdots \circ\mathfrak{l}_{1, \boldsymbol{\ep},\boldsymbol{\hat{\ep}}}(x,y).
\end{gather*}
 Analogously, if $I$ has infinitely many non-positive elements, set, for any $(x,y)$
with $x \geq x_i$,
\begin{gather*}
 \mathfrak{r}_{\boldsymbol{\ep},\hat{\boldsymbol{\ep}}}(x,y):=\mathfrak{r}_{i, \boldsymbol{\ep},\boldsymbol{\hat{\ep}}} \circ \mathfrak{r}_{i+1, \boldsymbol{\ep},\boldsymbol{\hat{\ep}}}\circ \cdots \circ \mathfrak{r}_{0, \boldsymbol{\ep},\boldsymbol{\hat{\ep}}}(x,y).
\end{gather*}
In the remaining cases, set $\mathfrak{l}_{\boldsymbol{\ep},\hat{\boldsymbol{\ep}}} = \mathrm{id} = \mathfrak{r}_{\boldsymbol{\ep},\hat{\boldsymbol{\ep}}}$. The following result is implied by the fact that any given point is only acted upon by finitely many elements of the composition (cf.\ Pelayo, Ratiu and \vungoc\ \cite[Lemma~2.1]{pvr_carto}).
\begin{Lemma}\label{lemma:wd} The maps $\mathfrak{l}_{\boldsymbol{\ep},\hat{\boldsymbol{\ep}}}\colon D_{\sup} \to D_{\sup}$ and $\mathfrak{r}_{\boldsymbol{\ep},\hat{\boldsymbol{\ep}}}\colon D_{\inf} \to D_{\inf}$ are well-defined.
\end{Lemma}

In fact, more is true.
\begin{Lemma}\label{lemma:well_defined} The maps $\mathfrak{l}_{\boldsymbol{\ep},\hat{\boldsymbol{\ep}}} \colon D_{\sup} \to D_{\sup}$ and $\mathfrak{r}_{\boldsymbol{\ep},\hat{\boldsymbol{\ep}}}\colon D_{\inf} \to D_{\inf}$ are homeomorphisms that are $\Z$-affine isomorphisms away from the set $\bigcup\limits_{i \in I} \{ (x,y) \,|\, x = x_i \}$.
\end{Lemma}

\begin{proof} The only non-trivial cases to consider are those of $\mathfrak{l}_{\boldsymbol{\ep},\hat{\boldsymbol{\ep}}} $ if~$I$ contains infinitely many positive elements and of $\mathfrak{r}_{\boldsymbol{\ep},\hat{\boldsymbol{\ep}}}$ if~$I$ contains infinitely many non-positive elements. Consider the former case (the latter is entirely analogous). To see that $\mathfrak{l}_{\boldsymbol{\ep},\hat{\boldsymbol{\ep}}} $ is a~homeomorphism, observe that the proof of Lemma~\ref{lemma:wd} implies that, for all $(x,y) \in D_{\sup}$, $\mathfrak{l}_{\boldsymbol{\ep},\hat{\boldsymbol{\ep}}}(x,y) = \big(x,\mathfrak{l}_{\boldsymbol{\ep},\hat{\boldsymbol{\ep}}}^{(2)}(x,y)\big)$, for some continuous function $\mathfrak{l}_{\boldsymbol{\ep},\hat{\boldsymbol{\ep}}}^{(2)} \colon D_{\sup} \to \R$. Moreover, it can be checked directly that, for any $x < \mathsf{x}_{\sup}$, the function $\mathfrak{l}_{\boldsymbol{\ep},\hat{\boldsymbol{\ep}}}^{(2)}(x,\cdot )$ is strictly increasing. Therefore, $\mathfrak{l}_{\boldsymbol{\ep},\hat{\boldsymbol{\ep}}}$ is a~homeomorphism onto its image. Since for any $i \geq 1$ the map $\mathfrak{l}_{i, \boldsymbol{\ep},\boldsymbol{\hat{\ep}}} \circ\mathfrak{l}_{i -1, \boldsymbol{\ep},\boldsymbol{\hat{\ep}}} \circ \cdots \circ\mathfrak{l}_{1, \boldsymbol{\ep},\boldsymbol{\hat{\ep}}}$ sends $\{(x,y) \,|\, x \leq x_i\}$ onto itself, it follows that $\mathfrak{l}_{\boldsymbol{\ep},\hat{\boldsymbol{\ep}}}(D_{\sup}) = D_{\sup}$. The fact that it is a $\Z$-affine isomorphism away from $\bigcup\limits_{i \in I} \{ (x,y) \,|\, x = x_i\}$ follows from the fact that, for any $i \geq 1$, $\mathfrak{l}_{i,\boldsymbol{\ep},\boldsymbol{\hat{\ep}}} \circ\mathfrak{l}_{ i -1, \boldsymbol{\ep},\boldsymbol{\hat{\ep}}} \circ \cdots \circ\mathfrak{l}_{1, \boldsymbol{\ep},\boldsymbol{\hat{\ep}}}$ also satisfies this property (see Lemma~\ref{cor:piece_z-aff_homeo}).
\end{proof}

With the maps $\mathfrak{l}_{\boldsymbol{\ep},\hat{\boldsymbol{\ep}}}$ and $\mathfrak{r}_{\boldsymbol{\ep},\hat{\boldsymbol{\ep}}} $ at hand, we can state the main result of this section. (Recall that $\operatorname{Vert}(2;\Z)$ denotes the subgroup of~$\mathrm{AGL}(2;\Z)$ preserving vertical lines, see Remark~\ref{rmk:toric_im}).

\begin{Theorem}\label{thm:different_carto} Let $(M,\omega,\Phi\!=\!(J,H))$ be a faithful semitoric system. Given any $\boldsymbol{\ep}, \hat{\boldsymbol{\ep}} \!\in\! \{{+}1{,}{-}1\}^I$ and any two cartographic homeomorphisms $f_{\boldsymbol{\ep}}, f_{\hat{\boldsymbol{\ep}}} \colon B \to \R^2$ associated to $\boldsymbol{\ep}$, $\hat{\boldsymbol{\ep}}$ respectively, there exists a transformation $h_{\boldsymbol{\ep}, \hat{\boldsymbol{\ep}}} \in \operatorname{Vert}(2;\Z)$ such that
 \begin{gather*} 
 f_{\hat{\boldsymbol{\ep}}} = \mathfrak{r}_{\boldsymbol{\ep},\hat{\boldsymbol{\ep}}} \circ
 \mathfrak{l}_{\boldsymbol{\ep},\hat{\boldsymbol{\ep}}} \circ h_{\boldsymbol{\ep}, \hat{\boldsymbol{\ep}}} \circ f_{\boldsymbol{\ep}}.
 \end{gather*}
\end{Theorem}

The main idea behind Theorem \ref{thm:different_carto} is not new; it first appears in \vungoc\ \cite[Proposition~4.1]{vu-ngoc} in the context of
semitoric systems. However, the more general context of faithful semitoric systems, where there may be infinitely many focus-focus points, deserves to be dealt with carefully. For instance, the transformation $\mathfrak{r}_{\boldsymbol{\ep},\hat{\boldsymbol{\ep}}}$, which is not needed in the study of semitoric systems, is necessary in this context (see Remarks~\ref{rmk:difference_san} and~\ref{rmk:maps_r}).

\begin{rmk}\label{rmk:compare_st} In fact, the statement and proof of Theorem~\ref{thm:different_carto} may be of use in the study of semitoric systems as well. The main references for these systems make the underlying (tacit) assumption that the signs of the cartographic homeomorphisms, as in Theorem~\ref{prop:rh}, are positive (cf.\ Pelayo and \vungoc\ \cite{pelayo-vu-ngoc-inventiones,pelayo-vu-ngoc-constr,vu-ngoc}).
\end{rmk}

\begin{proof}[Proof of Theorem \ref{thm:different_carto}] The proof is split into three steps:
 \begin{enumerate}[label=Step \arabic*:, ref = Step \arabic*, leftmargin=*]\itemsep=0pt
 \item \label{item:25} Construct the map $h_{\boldsymbol{\ep},\hat{\boldsymbol{\ep}}} \in \operatorname{Vert}(2;\Z)$.
 \item \label{item:26} Reduce to the simpler case in which $\boldsymbol{\ep}$ and $\hat{\boldsymbol{\ep}}$ differ in precisely one component.
 \item \label{item:27} Prove the simpler case.
 \end{enumerate}

{\bf Step 1: constructing the map $h_{\boldsymbol{\ep},\hat{\boldsymbol{\ep}}} \in \operatorname{Vert}(2;\Z)$.} As above, fix a basepoint $(\mathsf{x},\mathsf{y}) \in B_{\mathrm{wt}}$ with the property that $x_0 < \mathsf{x} < x_1$. (If $x_0$ or $ x_1$ is not defined, then only the other inequality is required.) Denote by $\mathfrak{S}$ the connected component of $S_{\boldsymbol{\ep}} \cap S_{\hat{\boldsymbol{\ep}}}$ containing $(\mathsf{x},\mathsf{y})$. Since both~$S_{\boldsymbol{\ep}}$ and~$S_{\hat{\boldsymbol{\ep}}}$ are open, so is $\mathfrak{S}$. Moreover it is path-connected by definition. Furthermore, it can be checked that $\mathfrak{S}$ intersects any vertical line either in an empty or in a~connected set. Therefore, by Proposition~\ref{prop:sub_mostly_vat}, the subsystem of $\vat$ relative to $\mathfrak{S}$ is faithful semitoric and contains $(\mathsf{x},\mathsf{y})$. Moreover, by construction, it contains no focus-focus value. By Lemma~\ref{cor:carto_subsys}, the maps $f_{\boldsymbol{\ep}}|_{\mathfrak{S}}$ and $f_{\hat{\boldsymbol{\ep}}}|_{\mathfrak{S}}$ are cartographic homeomorphisms for the subsystem of~$\vat$ relative to~$\mathfrak{S}$. Therefore, by Remark~\ref{rmk:freedom_no_ff}, there exists $h_{\boldsymbol{\ep},\hat{\boldsymbol{\ep}}} \in \operatorname{Vert}(2;\Z)$ such that $f_{\hat{\boldsymbol{\ep}}}|_{\mathfrak{S}} = h_{\boldsymbol{\ep},\hat{\boldsymbol{\ep}}}|_{f_{\boldsymbol{\ep}}(\mathfrak{S})} \circ f_{\boldsymbol{\ep}}|_{\mathfrak{S}}$. The map $h_{\boldsymbol{\ep},\hat{\boldsymbol{\ep}}}$ is the desired one.

{\bf Step 2: reducing to a simpler case.} Observe that, by Corollary~\ref{cor:freedom_ff_sc}, the map $h_{\boldsymbol{\ep},\hat{\boldsymbol{\ep}}} \circ f_{\boldsymbol{\ep}}$ is a cartographic homeomorphism associated to~$\boldsymbol{\ep}$. Moreover, the above argument shows that $f_{\hat{\boldsymbol{\ep}}}|_{\mathfrak{S}} = (h_{\boldsymbol{\ep},\hat{\boldsymbol{\ep}}} \circ
f_{\boldsymbol{\ep}}) |_{\mathfrak{S}}$. Thus, in order to prove the result in the statement of the theorem, it suffices to prove that, if $f_{\hat{\boldsymbol{\ep}}}|_{\mathfrak{S}} = f_{\boldsymbol{\ep}}|_{\mathfrak{S}}$, then $f_{\hat{\boldsymbol{\ep}}} = \mathfrak{r}_{\boldsymbol{\ep},\hat{\boldsymbol{\ep}}} \circ \mathfrak{l}_{\boldsymbol{\ep},\hat{\boldsymbol{\ep}}} \circ f_{\boldsymbol{\ep}}$. Henceforth, assume that $f_{\hat{\boldsymbol{\ep}}}|_{\mathfrak{S}} = f_{\boldsymbol{\ep}}|_{\mathfrak{S}}$, which implies $\operatorname{sgn}(f_{\hat{\boldsymbol{\ep}}}) = \operatorname{sgn}(f_{\boldsymbol{\ep}})$. In fact, we can simplify the argument further: it suffices to prove the claimed result under the assumption that all but one of the components of $\boldsymbol{\ep}$, $\hat{\boldsymbol{\ep}}$ are equal. For, if the latter holds, we can argue as follows. Consider a sequence of choices of signs~$\boldsymbol{\ep}_s$, for $s \in \{1,2,\ldots\} \cup \{\infty\}$, such that $\boldsymbol{\ep}_1 = \boldsymbol{\ep}$, $\boldsymbol{\ep}_{\infty} = \hat{\boldsymbol{\ep}}$, and, for any $s \geq 1$, all but one component of~$\boldsymbol{\ep}_s$ and~$\boldsymbol{\ep}_{s+1}$ are equal. Moreover, if $\boldsymbol{\ep}$, $\hat{\boldsymbol{\ep}}$ differ in finitely components, say in~$r$ components, choose the above sequence so that for all $s \geq r+1$, $\boldsymbol{\ep}_{s} = \hat{\boldsymbol{\ep}}$. For each $s \geq 1$, fix a~choice of cartographic homeomorphism $f_{\boldsymbol{\ep}_s}$ with the property that $f_{\boldsymbol{\ep}_s}|_{\mathfrak{S}} = f_{\boldsymbol{\ep}}|_{\mathfrak{S}}$. Moreover, require that $f_{\boldsymbol{\ep}_1} = f_{\boldsymbol{\ep}}$, that $f_{\boldsymbol{\ep}_{\infty}} = f_{\hat{\boldsymbol{\ep}}}$, and that, if $\boldsymbol{\ep}$, $\hat{\boldsymbol{\ep}}$ differ in precisely~$r$ components, then for all $s \geq r+1$, $f_{\boldsymbol{\ep}_{s}} = f_{\hat{\boldsymbol{\ep}}}$ .
Using the above sequence of signs and associated cartographic homeomorphisms and the fact that the claimed result holds when all but one component of the signs are equal, we obtain, for all $s \geq 1$, maps $\mathfrak{l}_{\boldsymbol{\ep}_s,\boldsymbol{\ep}_{s+1}}$, $\mathfrak{r}_{\boldsymbol{\ep}_s,\boldsymbol{\ep}_{s+1}}$ satisfying
\begin{gather} \label{eq:3}
 f_{\boldsymbol{\ep}_{s+1}} = \mathfrak{r}_{\boldsymbol{\ep}_s,\boldsymbol{\ep}_{s+1}} \circ \mathfrak{l}_{\boldsymbol{\ep}_s,\boldsymbol{\ep}_{s+1}} \circ f_{\boldsymbol{\ep}_s}.
\end{gather}
Therefore, iterating equation \eqref{eq:3}, for all $s \geq 1$,
\begin{gather} \label{eq:4}
 f_{\boldsymbol{\ep}_{s+1}} = \mathfrak{r}_{\boldsymbol{\ep}_s,\boldsymbol{\ep}_{s+1}} \circ
 \mathfrak{l}_{\boldsymbol{\ep}_s,\boldsymbol{\ep}_{s+1}} \circ \mathfrak{r}_{\boldsymbol{\ep}_{s-1},\boldsymbol{\ep}_{s}} \circ
 \mathfrak{l}_{\boldsymbol{\ep}_{s-1},\boldsymbol{\ep}_{s}} \circ
 \cdots \circ \mathfrak{r}_{\boldsymbol{\ep}_1,\boldsymbol{\ep}_{2}} \circ
 \mathfrak{l}_{\boldsymbol{\ep}_1,\boldsymbol{\ep}_{2}} \circ
 f_{\boldsymbol{\ep}},
\end{gather}
where we use the fact that $f_{\boldsymbol{\ep}} =f_{\boldsymbol{\ep}_1}$. If $\boldsymbol{\ep}$ and $\hat{\boldsymbol{\ep}}$ differ in precisely $r$ components, then, by construction, for all $s \geq r+1$,
\begin{gather*}
 \mathfrak{r}_{\boldsymbol{\ep}_s,\boldsymbol{\ep}_{s+1}} =\mathrm{id} = \mathfrak{l}_{\boldsymbol{\ep}_s,\boldsymbol{\ep}_{s+1}}.
\end{gather*}
Therefore, in this case, equation \eqref{eq:4} yields that
\begin{gather*}
 f_{\hat{\boldsymbol{\ep}}} = \mathfrak{r}_{\boldsymbol{\ep}_r,\boldsymbol{\ep}_{r+1}} \circ
 \mathfrak{l}_{\boldsymbol{\ep}_r,\boldsymbol{\ep}_{r+1}} \circ \mathfrak{r}_{\boldsymbol{\ep}_{r-1},\boldsymbol{\ep}_{r}} \circ
 \mathfrak{l}_{\boldsymbol{\ep}_{r-1},\boldsymbol{\ep}_{r}} \circ \cdots \circ \mathfrak{r}_{\boldsymbol{\ep}_1,\boldsymbol{\ep}_{2}} \circ \mathfrak{l}_{\boldsymbol{\ep}_1,\boldsymbol{\ep}_{2}} \circ f_{\boldsymbol{\ep}}.
\end{gather*}
Because the homeomorphisms in the above composition commute (see Lemma~\ref{cor:commute_piece}), by definition of $\mathfrak{l}_{\boldsymbol{\ep},\hat{\boldsymbol{\ep}}}$ and $\mathfrak{r}_{\boldsymbol{\ep},\hat{\boldsymbol{\ep}}}$,
\begin{gather*} 
 f_{\hat{\boldsymbol{\ep}}} = \mathfrak{r}_{\boldsymbol{\ep},\hat{\boldsymbol{\ep}}} \circ \mathfrak{l}_{\boldsymbol{\ep},\hat{\boldsymbol{\ep}}}.
\end{gather*}

Thus the result is proved if $\boldsymbol{\ep}$ and $\hat{\boldsymbol{\ep}}$ differ in finitely many components.

The case in which they differ by infinitely many components is entirely analogous, as we can consider the composite of infinitely many maps of the above form on, say, any compact subset of~$B$ (see the proof of Lemma~\ref{lemma:wd}) and use a compact exhaustion of~$B$. Thus, assuming that the result holds when the choices of signs differ in precisely one component, the result holds in general.

{\bf Step 3: proving the simple case.} Assume that $f_{\hat{\boldsymbol{\ep}}}|_{\mathfrak{S}} = f_{\boldsymbol{\ep}}|_{\mathfrak{S}}$ and that $\boldsymbol{\ep}$ and $\hat{\boldsymbol{\ep}}$ differ in precisely one component. Under these assumptions the result can be proved exactly as in
\vungoc\ \cite[Proposition~4.1]{vu-ngoc}, whose key ideas are explained below. Suppose that $\boldsymbol{\ep}$ and $\hat{\boldsymbol{\ep}}$ differ precisely in the $i$th component. By Corollary~\ref{cor:carto_iso}, $f_{\hat{\boldsymbol{\ep}}} \circ f_{\boldsymbol{\ep}}^{-1}$ is a homeomorphism that is piecewise $\Z$-affine. Using the proof of Theorem~\ref{prop:rh}, it may be assumed without loss of generality that $S_{\boldsymbol{\ep}}$ and $S_{\hat{\boldsymbol{\ep}}}$ are path-connected. As both sets are dense, it suffices to check the desired equality on $S_{\boldsymbol{\ep}} \cap S_{\hat{\boldsymbol{\ep}}}$. Since $\boldsymbol{\ep}$ and $\hat{\boldsymbol{\ep}}$ differ in precisely one component, it follows that $S_{\boldsymbol{\ep}}\cap S_{\hat{\boldsymbol{\ep}}}$ has two connected components, $\mathfrak{S}$ and $\mathfrak{S}'$, which are open and satisfy the assumptions of Proposition~\ref{prop:sub_mostly_vat}. Thus the subsystems of~$\vat$ relative to~$\mathfrak{S}$ and~$\mathfrak{S}'$ are faithful semitoric. By Lemma~\ref{cor:carto_subsys}, the restrictions of $f_{\boldsymbol{\ep}}, f_{\hat{\boldsymbol{\ep}}}$ to~$\mathfrak{S}$ and~$\mathfrak{S}'$ are cartographic homeomorphisms for the respective subsystems of~$\vat$. By construction, these subsystems contain no focus-focus points. Thus, Remark~\ref{rmk:freedom_no_ff} implies that there exist $h_{\mathfrak{S}},
h_{\mathfrak{S}'} \in \operatorname{Vert}(2;\Z)$ with $f_{\hat{\boldsymbol{\ep}}}|_{\mathfrak{S}} = h_{\mathfrak{S}} \circ f_{\boldsymbol{\ep}}|_{\mathfrak{S}}$ and $f_{\hat{\boldsymbol{\ep}}}|_{\mathfrak{S'}} = h_{\mathfrak{S'}} \circ f_{\boldsymbol{\ep}}|_{\mathfrak{S'}}$.

By assumption, $f_{\hat{\boldsymbol{\ep}}}|_{\mathfrak{S}} = f_{\boldsymbol{\ep}}|_{\mathfrak{S}}$, so $h_{\mathfrak{S}} = \mathrm{id}$; on the other hand, the above assumptions imply that $(\mathfrak{r}_{\boldsymbol{\ep},\hat{\boldsymbol{\ep}}} \circ \mathfrak{l}_{\boldsymbol{\ep},\hat{\boldsymbol{\ep}}})|_{\mathfrak{S}} = \mathrm{id}$. Thus the desired equality holds on $\mathfrak{S}$. It remains to check that it does on $\mathfrak{S}'$. Using property~\ref{item:17} for $f_{\boldsymbol{\ep}}$ and $f_{\hat{\boldsymbol{\ep}}}$ and the fact that $\operatorname{sgn}(f_{\boldsymbol{\ep}}) = \operatorname{sgn}( f_{\hat{\boldsymbol{\ep}}})$, it can be shown that the linear parts of~$h_{\mathfrak{S}'}$ and of~$(\mathfrak{r}_{\boldsymbol{\ep},\hat{\boldsymbol{\ep}}} \circ \mathfrak{l}_{\boldsymbol{\ep},\hat{\boldsymbol{\ep}}})|_{\mathfrak{S}'}$ are equal. To see that their translational components are equal, observe that the piecewise $\Z$-affine transformation given on~$f_{\boldsymbol{\ep}}(\mathfrak{S})$ and on $f_{\boldsymbol{\ep}}(\mathfrak{S}')$ by~$h_{\mathfrak{S}}$ and~$h_{\mathfrak{S}'}$, respectively, extends uniquely to a topological embedding of $f_{\boldsymbol{\ep}}(B)$ onto~$f_{\hat{\boldsymbol{\ep}}}(B)$ (which equals $f_{\hat{\boldsymbol{\ep}}} \circ f_{\boldsymbol{\ep}}^{-1}$). In particular, it acts as the identity on the vertical line containing the $i$th focus-focus value. This implies that the translational component of~$h_{\mathfrak{S}'}$ equals that of $(\mathfrak{r}_{\boldsymbol{\ep},\hat{\boldsymbol{\ep}}} \circ \mathfrak{l}_{\boldsymbol{\ep},\hat{\boldsymbol{\ep}}})|_{\mathfrak{S}'}$.
\end{proof}

To conclude this section, we enlarge the set of semitoric cartographic images of a faithful semitoric system with at least one focus-focus point to construct an invariant of the system up to isomorphisms of {\em integrable systems}. To this end, recall that $\mathsf{L} \subset \mathrm{AGL}(2;\Z)$ is the subset consisting of elements of the form
 \begin{gather*}
 \left( \begin{pmatrix}
 \eta_1 & 0 \\
 k & \eta_2
 \end{pmatrix}, \begin{pmatrix}
 a \\
 b
 \end{pmatrix}\right),
 \end{gather*}
where, for $i=1,2$, $\eta_i \in \{+1,-1\}$, $k \in \Z$, and $a,b \in \R$ (see Remark~\ref{rmk:toric_im}). Given a faithful semitoric system $\is$ with at least one focus-focus point, let $(f_{\boldsymbol{\ep}},S_{\boldsymbol{\ep}})$ be one of its semitoric cartographic pairs. By Lemma~\ref{lem:infinitely_many_carto}, given any element $\mathsf{l} \in \mathsf{L}$, $(\mathsf{l} \circ f_{\boldsymbol{\ep}},S_{\boldsymbol{\ep}})$ is a~cartographic pair for $\is$. Cartographic pairs of the form $(\mathsf{l} \circ f_{\boldsymbol{\ep}},S_{\boldsymbol{\ep}})$, where $\mathsf{l} \in \mathsf{L}$ and $(f_{\boldsymbol{\ep}},S_{\boldsymbol{\ep}})$ is a~semitoric cartographic pair, are said to be {\em complexity one}, and so are their images. The set of complexity one images of $\is$ is an invariant of the system up to isomorphisms of integrable systems.

\begin{Lemma}\label{lemma:im_inv_full} Let $(M_1,\omega_1,\Phi_1)$ and $(M_2,\omega_2,\Phi_2)$ be faithful semitoric systems that are isomorphic as integrable systems and suppose that $(M_1,\omega_1,\Phi_1)$ has at least one focus-focus point. Then the sets of complexity one cartographic images of $(M_1,\omega_1,\Phi_1)$ and $(M_2,\omega_2,\Phi_2)$ are equal.
 \end{Lemma}

\begin{proof} Let $(\Psi,\psi)$ be an isomorphism of $(M_1,\omega_1,\Phi_1)$ and $(M_2,\omega_2,\Phi_2)$ as integrable systems. Since $(M_1,\omega_1,\Phi_1)$ has at least one focus-focus point, Proposition~\ref{prop:im_equ} implies that $(\Psi,\psi)$ is an isomorphism as complexity one spaces, i.e., $\psi(x,y) = \big(\eta x + a,\psi^{(2)}(x,y)\big)$ for some $\eta \in \{\pm 1\}$ and some $a \in \R$ (see Definition~\ref{defn:im_faithful_st}). For an arbitrary $\mathsf{l}\in\mathsf{L}$ and a choice of signs~$\boldsymbol{\ep}_2$ that make $S_{\boldsymbol{\ep}_2} \subset B_2$ connected (a choice that exists by Corollary~\ref{cor:there_exists_choice}), consider the complexity one cartographic pair $(\mathsf{l} \circ f_{\boldsymbol{\ep}_2},S_{\boldsymbol{\ep}_2})$ for $(M_2,\omega_2,\Phi_2)$. Arguing as in the proof of Lemma~\ref{lemma:im_invariant}, $\big(\mathsf{l} \circ f_{\boldsymbol{\ep}_2}\circ \psi,\psi^{-1}(S_{\boldsymbol{\ep}_2})\big)$ is a~cartographic pair for $(M_1,\omega_1,\Phi_1)$. Moreover, since $\psi$ is of the above form, there exists a unique choice of signs $\boldsymbol{\ep}_1$ for $(M_1,\omega_1,\Phi_1)$ such that $\psi^{-1}(S_{\boldsymbol{\ep}_2}) = S_{\boldsymbol{\ep}_1}$. Since~$\psi$ is a diffeomorphism and $S_{\boldsymbol{\ep}_2}$ is connected, so is $S_{\boldsymbol{\ep}_1}$. Let $f_{\boldsymbol{\ep}_1} \colon B_1 \to \R^2$ be a cartographic homeomorphism associated to $\boldsymbol{\ep}_1$ as in Theorem~\ref{prop:rh}. Observe that, by construction $f_{\boldsymbol{\ep}_1}(S_{\boldsymbol{\ep}_1})$ is connected. By Corollary~\ref{cor:carto_iso}, the restriction of the map $\mathsf{l} \circ f_{\boldsymbol{\ep}_2} \circ \psi \circ f^{-1}_{\boldsymbol{\ep}_1}$ to $f_{\boldsymbol{\ep}_1}(S_{\boldsymbol{\ep}_1})$ is a~$\Z$-affine isomorphism onto $\mathsf{l} \circ f_{\boldsymbol{\ep}_2}(S_{\boldsymbol{\ep}_2})$. Since $f_{\boldsymbol{\ep}_1}(S_{\boldsymbol{\ep}_1})$ is connected, there exists an element $ h \in \mathrm{AGL}(2;\Z)$ such that $h = \mathsf{l} \circ f_{\boldsymbol{\ep}_2} \circ \psi \circ f^{-1}_{\boldsymbol{\ep}_1}$ on $f_{\boldsymbol{\ep}_1}(S_{\boldsymbol{\ep}_1})$. As $f_{\boldsymbol{\ep}_1}(S_{\boldsymbol{\ep}_1}) \subset f_{\boldsymbol{\ep}_1}(B_1)$ is dense, it follows that $h = \mathsf{l} \circ f_{\boldsymbol{\ep}_2} \circ \psi \circ f^{-1}_{\boldsymbol{\ep}_1}$ on $f_{\boldsymbol{\ep}_1}(B_1)$ or, equivalently, $h \circ f_{\boldsymbol{\ep}_1} = \mathsf{l} \circ f_{\boldsymbol{\ep}_2} \circ \psi$. Moreover, observe that $h \in \mathsf{L}$, for $\mathsf{l} \in \mathsf{L}$, $(f_{\boldsymbol{\ep}_i},S_{\boldsymbol{\ep}_i})$ are semitoric cartographic pairs and $\psi(x,y) = \big(\eta x + a,\psi^{(2)}(x,y)\big)$. This shows that the cartographic pair $\big(\mathsf{l} \circ f_{\boldsymbol{\ep}_2} \circ \psi, \psi^{-1}(S_{\boldsymbol{\ep}_2})\big) = (h \circ f_{\boldsymbol{\ep}_1},S_{\boldsymbol{\ep}_1})$ for $(M_1,\omega_1,\Phi_1)$ is complexity one. Since the set of images of complexity one cartographic pairs $(\mathsf{l} \circ f_{\boldsymbol{\ep}_2},S_{\boldsymbol{\ep}_2})$ such that $S_{\boldsymbol{\ep}_2}$ is connected equals the set of images of complexity one cartographic pairs of~$(M_2,\omega_2,\Phi_2)$ (see Remark~\ref{rmk:connected_suffices}), this shows that the set of complexity one cartographic images of $(M_2,\omega_2,\Phi_2)$ is contained in that of $(M_1,\omega_1,\Phi_1)$. Reversing the roles of $(M_1,\omega_1,\Phi_1)$ and $(M_2,\omega_2,\Phi_2)$ in the above argument completes the proof.
 \end{proof}

\subsection[$\eta$-cartographic faithful semitoric systems]{$\boldsymbol{\eta}$-cartographic faithful semitoric systems}\label{sec:choos-an-appr}
Let $\is$ be a faithful semitoric system. The presence of a focus-focus point implies that no faithful semitoric system isomorphic to~$\is$ has a cartographic moment map (see Section~\ref{sec:faithful}). On the other hand, Theorem~\ref{prop:rh} provides cartographic homeomorphisms associated to choices of vertical cuts. Fix any such cartographic homeomorphism~$f_{\boldsymbol{\ep}}$; while it is tempting to think of $(M,\omega,f_{\boldsymbol{\ep}}\circ\Phi)$ as an integrable system, the lack of smoothness of~$f_{\boldsymbol{\ep}}$ prevents it from being one. (If we were to adopt the non-standard convention of Harada and Kaveh \cite[Definition~2.1]{hk}, $(M,\omega,f_{\boldsymbol{\ep}}\circ\Phi)$ would be an integrable system.) The aim of this section is to show that, in some sense, the next best scenario holds: Given a choice of signs~$\boldsymbol{\ep}$ with~$S_{\boldsymbol{\ep}}$ connected, any cartographic homeomorphism associated to $\boldsymbol{\ep}$ can be modified in an arbitrarily small neighborhood of the cuts associated to $\boldsymbol{\ep}$ so that it becomes everywhere smooth (see Theorem~\ref{thm:rect-embedding} for a~precise statement). This smoothing of cartographic homeomorphisms generates representatives in the isomorphism class of a faithful semitoric system, which we call {\em $\boldsymbol{\eta}$-cartographic}, that are particularly useful when defining surgeries on (isomorphism classes of) faithful semitoric systems (cf.\ the forthcoming~\cite{HSSS-surgeries}). (The above sentence holds for the finest notion of isomorphism of faithful semitoric systems given by Definition~\ref{defn:im_faithful_st} and, therefore, for all.) Moreover, we show that cartographic homeomorphisms are, in some sense, limits of what we call {\em $\boldsymbol{\eta}$-cartographic embeddings} (see Proposition~\ref{prop:limit}).

Throughout the rest of this section, a faithful semitoric system $\vat$ containing at least one focus-focus singular point is fixed. (Recall that~$\vat$ is required to satisfy the generic assumption~\ref{item:7} throughout.) As above, set $B= \Phi(M)$, let $\{c_i\}_{i \in I} \subset \operatorname{Int}(B)$ denote the set of focus-focus values of $\is$, where~$I$ is ordered as in equation~\eqref{eq:2}, while~$l^{\boldsymbol{\ep}}$ denotes the union of the vertical cuts in $B$ associated to a~choice of signs $\boldsymbol{\ep} \in \{ +1,-1\}^I$. Also let $S_{\boldsymbol{\ep}} = B \smallsetminus l^{\boldsymbol{\ep}}$ denote the complement of those cuts. Finally, assume a~cartographic homeomorphism has positive sign (see Theorem~\ref{prop:rh}) unless otherwise stated.

\subsubsection{Admissible half-strips for faithful semitoric systems}\label{sec:admissible}

First, we define the (closed) neighborhoods of vertical half-lines that we use to construct the smoothing of a given cartographic homeomorphism.

\begin{figure}[h] \centering
 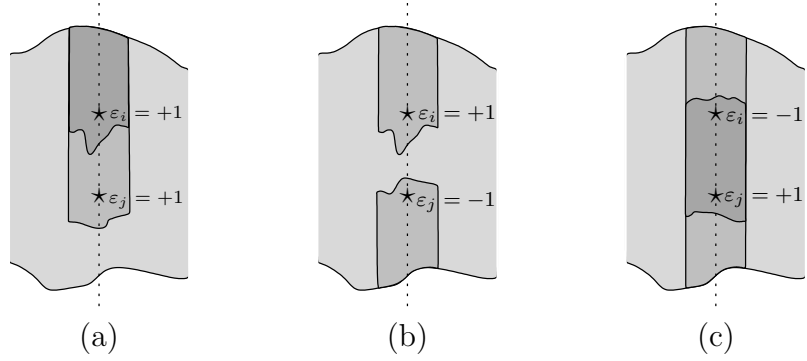
\caption{The symbol $\star$ indicates the points at which the half-strips are centered. Figures (a) and (b) show admissible half-strips with the same center line. Figure (c) shows half-strips that are not admissible.} \label{halfStrips}
\end{figure}

\begin{Definition}\label{defn:half-strip} Fix $\ep \in \{+1,-1\}$, $(x_0,y_0)\in\R^2$, $\eta>0$, and a continuous map $\gamma\colon \big[x_0 - \frac{\eta}{2},x_0 + \frac{\eta}{2} \big] \to \R$ satisfying $\ep y_0 > \ep \gamma(x)$ for all $x \in \big[x_0 - \frac{\eta}{2},x_0 + \frac{\eta}{2} \big]$. A {\em half-strip centered at $(x_0,y_0)$ of sign $\ep$ and width $\eta$ with bounding curve $\gamma$} is the following closed subset of $\R^2$:
 \begin{gather*}
 \sigma_{\eta,\gamma}^\ep(x_0,y_0): = \left\{(x,y) \,|\, x_0-\frac{\eta}{2} \leq x \leq x_0+\frac{\eta}{2} \text{ and } \ep y \geq \ep \gamma(x) \right\}
 \end{gather*}
(see Fig.~\ref{halfStrips}). The vertical line $ \{(x,y) \,|\, x =x_0 \}$ is called the {\em center line} of the half-strip. When the center point $(x_0,y_0)$ and the bounding curve $\gamma$ are not of particular concern, the half-strip is denoted by $\sigma_\eta^\ep$. The {\em base} of a half-strip $\sigma^{\ep}_{\eta,\gamma}(x_0,y_0)$ is the subset
 \begin{gather*}
 \sigma^{\ep}_{\eta,\gamma}(x_0,y_0) \cap \left\{ (x,y) \,|\, \ep y < \ep y_0 + \frac{\eta}{2}\right\}.
 \end{gather*}
\end{Definition}

Consider a choice of (countably many) points $\{(x_i,y_i)\}_{i \in I}$, of signs $\boldsymbol{\ep} \in \{+1,-1\}^I$, of positive numbers $\boldsymbol{\eta} \in \{\eta_i \}_{i \in I}$, and of continuous curves $\boldsymbol{\gamma} =\{\gamma_i\}_{i \in I}$. Let $\sigma^{\ep_i}_{\eta_i, \gamma_i}(x_i,y_i)$ be the half-strip centered at $(x_i,y_i)$ of sign $\ep_i$ and width $\eta_i>0$ with bounding curve $\gamma_i$. Moreover, set
\begin{gather*}
 \boldsymbol{\sigma}^{\boldsymbol{\ep}}_{\boldsymbol{\eta},\boldsymbol{\gamma}}:
= \bigcup_{i} \sigma^{\ep_i}_{\eta_i,\gamma_i} (x_i,y_i),
\end{gather*}
and denote the above choices of signs, widths and curves by the triple $(\boldsymbol{\ep},\boldsymbol{\eta},\boldsymbol{\gamma})$.

\begin{Definition}\label{defn:admissible_subset} Suppose $\mathsf{B} \subset \R^2$ has the property that its intersection with any vertical line is either empty or path-connected, and consider a countable set of points $\{(x_i,y_i)\}_{i \in I}$ therein. A triple $(\boldsymbol{\ep},\boldsymbol{\eta},\boldsymbol{\gamma})$ as above is {\em admissible} for the subset~$\mathsf{B}$ relative to the points $\{(x_i,y_i)\}_{i \in I}$ if it satisfies the following conditions:
 \begin{itemize}[leftmargin=*]\itemsep=0pt
 \item For all $i$, the base of the half-strip $\sigma^{\ep_i}_{\eta_i,\gamma_i}(x_i,y_i)$ is contained in $\operatorname{Int}(\mathsf{B})$.
 \item If $(x_i,y_i) \in \sigma^{\ep_j}_{\eta_j,\gamma_j}(x_j,y_j)$ for $i \neq j$, then $x_i = x_j$.
 \item Whenever the half-strips $\sigma^{\ep_i}_{\eta_i,\gamma_i}(x_i,y_i)$ and $\sigma^{\ep_j}_{\eta_j,\gamma_j}(x_j,y_j)$ share the same center line, $\eta_i = \eta_j$.
 \item The intersection of any two distinct half-strips is either empty or equal to one of the half-strips.
 \end{itemize}
 In this case, the corresponding half-strips are called {\em admissible} for $\mathsf{B}$ relative to the points $\{(x_i,y_i)\}_{i \in I}$.
\end{Definition}

Examples of admissible half-strips are sketched in Fig.~\ref{halfStrips}(a) and (b).

\begin{Definition}\label{defn:tabs} Let $\is$ be a faithful semitoric system whose set of focus-focus values is indexed by $I$ as in equation~\eqref{eq:2}. A~triple $(\boldsymbol{\ep},\boldsymbol{\eta},\boldsymbol{\gamma})$ as in Definition \ref{defn:admissible_subset} and their corresponding half-strips are {\em admissible} for $\is$ if they are admissible for~$B$ relative to the set of focus-focus values~$B_{\mathrm{ff}}= \{c_i\}_{i \in I}$.
\end{Definition}

Before establishing the existence of admissible half-strips for any faithful semitoric system (Proposition~\ref{prop:existence-admissible}), we derive the following necessary condition. Recall that the focus-focus values of a faithful semitoric system are ordered as in Section~\ref{sec:cuts-cart-home}.

\begin{Proposition}\label{prop:path-connected} If $(\boldsymbol{\ep},\boldsymbol{\eta},\boldsymbol{\gamma})$ is an admissible triple for $\is$, then $S_{\boldsymbol{\ep}} = B \smallsetminus l^{\boldsymbol{\ep}}$ is contractible.
\end{Proposition}

\begin{proof} Fix an admissible triple $(\boldsymbol{\ep},\boldsymbol{\eta},\boldsymbol{\gamma})$. By Lemma \ref{lemma:S_ep_path_conn} and Corollary \ref{cor:path_conn_iff_contr}, it suffices to check that, if $i > j$ and there are focus-focus points $c_i=(x_i,y_i)$, $c_j = (x_j,y_j)$ with $x_i = x_j$, then $\ep_i \geq \ep_j$. Suppose not, then the half-strips $\sigma^{\ep_i}_{\eta_i,\gamma_i}(x_i,y_i)$, $\sigma^{\ep_j}_{\eta_j,\gamma_j}(x_j,y_j)$ intersect, but neither is contained in the other (see Fig.~\ref{halfStrips}(c)), thus contradicting admissibility of the given triple.
\end{proof}

Recall that, by Corollary \ref{cor:path_conn_iff_contr}, $S_{\boldsymbol{\ep}}$ is path-connected if and only if it is contractible. The next result establishes the converse to Proposition~\ref{prop:path-connected}.

\begin{Proposition}\label{prop:existence-admissible} Given a faithful semitoric system $\is$ with at least one focus-focus value and any choice $\boldsymbol{\ep} \in \{+1,-1\}^I$ making $S_{\boldsymbol{\ep}}$ path-connected, there exist a choice of positive numbers $\boldsymbol{\eta} = \{\eta_i \}_{i \in I}$ and of continuous curves $\boldsymbol{\gamma} = \{ \gamma_i\}_{i \in I}$ such that the triple $(\boldsymbol{\ep},\boldsymbol{\eta}, \boldsymbol{\gamma})$ is admissible for~$\is$. Moreover, the widths $\boldsymbol{\eta}$ can be chosen so that if a half-strip $\sigma^{\ep_i}_{\eta_i}$ contains a corner of~$B$, then it contains precisely one, and that corner lies on the center line of the half-strip.
\end{Proposition}

\begin{proof} Fix a choice of $\boldsymbol{\ep} \in \{+1,-1\}^I$ as above and let $\{c_i = (x_i,y_i)\}_{i \in I}$ denote the set of focus-focus values of~$\vat$. By Proposition~\ref{lemma:iso_fixed_points}, the set of first coordinates of focus-focus values is a subset of the set of critical values of~$J$, which, by property~\ref{item:100}, does not contain any limit point in~$J(M)$. Moreover, focus-focus values are discrete in~$B$ and, by property~\ref{item:101}, there are finitely many of them on a given vertical line. The above facts imply that there exists a choice of positive numbers $\boldsymbol{\eta} = \{\eta_i\}_{i \in I}$ such that
 \begin{itemize}[leftmargin=*]\itemsep=0pt
 \item if $x_i = x_j$, $\eta_i = \eta_j$,
 \item if, for $i \neq j$, $x_j \in \big[x_i - \frac{\eta_i}{2}, x_i + \frac{\eta_i}{2}\big]$, then $x_i = x_j$, and
 \item for all $i$, $\big[x_i - \frac{\eta_i}{2}, x_i + \frac{\eta_i}{2}\big] \times \big]\ep_i y_i - \frac{\eta_i}{2}, \ep_iy_i + \frac{\eta_i}{2}\big[$ is contained in $\operatorname{Int}(B)$.
 \end{itemize}
For each $i \in I$, define $\gamma_i\colon \big[x_i - \frac{\eta_i}{2}, x_i + \frac{\eta_i}{2}\big] \to \R$ to be $\gamma_i(x):=y_i - \ep_i \frac{\eta}{2}$, and set $\boldsymbol{\gamma} = \{\gamma_i\}_{i \in I}$. It can be checked that the triple $(\boldsymbol{\ep},\boldsymbol{\eta},\boldsymbol{\gamma})$ is admissible for $\vat$. Moreover, if the $\eta_i$'s are chosen to be sufficiently small then any half-strip that contains a~corner contains precisely one, which lies on the center line of the strip.
\end{proof}

Henceforth, any admissible triple for a faithful semitoric system is assumed to satisfy the final property in Proposition~\ref{prop:existence-admissible} unless otherwise stated. Corollary \ref{cor:there_exists_choice} and Proposition~\ref{prop:existence-admissible} readily yield the following result, stated below without proof.

 \begin{Corollary}\label{cor:existence_admissible} Every faithful semitoric system has an admissible triple.
 \end{Corollary}

Fix an admissible triple $(\boldsymbol{\ep},\boldsymbol{\eta}, \boldsymbol{\gamma})$ for $\vat$. The next results can be interpreted as showing that the complement of the corresponding half-strips in $B$ behaves like~$S_{\boldsymbol{\ep}}$.

\begin{Lemma}\label{lemma:complement_adm_strip_open} If $(\boldsymbol{\ep},\boldsymbol{\eta}, \boldsymbol{\gamma})$ is admissible for a~faithful semitoric system $\is$, then $B \smallsetminus \boldsymbol{\sigma}^{\boldsymbol{\ep}}_{\boldsymbol{\eta},\boldsymbol{\gamma}}$ is open in~$B$.
\end{Lemma}

\begin{proof} As in the proof of Corollary \ref{cor:complement_open}, it suffices to show that $\boldsymbol{\sigma}^{\boldsymbol{\ep}}_{\boldsymbol{\eta},\boldsymbol{\gamma}}$ is closed in $B$. Let $\{(x_n,y_n)\} \subset \boldsymbol{\sigma}^{\boldsymbol{\ep}}_{\boldsymbol{\eta},\boldsymbol{\gamma}}$ be a sequence that converges to $(x_0,y_0) \in B$ and consider the sequence $\{x_n\} = \operatorname{pr}_1(\{(x_n,y_n)\})$ which converges to $x_0 \in J(M)$. Since $J(M)$ is locally compact, there exists a~compact neighborhood $K \subset J(M)$ of~$x_0$. Since $x_n \to x_0$, it follows that all but finitely many of the~$x_n$ are contained in~$K$. Since $K$ is compact and the critical values of $J$ are discrete in $J(M)$ by property~\ref{item:100}, $K$~contains at most finitely many critical values of~$J$. Therefore, by property~\ref{item:101} and Proposition~\ref{lemma:iso_fixed_points}, there are at most finitely many focus-focus values contained in $\operatorname{pr}^{-1}_1(K) \cap B$. Hence, all but finitely many of the $(x_n,y_n)$ are contained in the union of finitely many admissible half-strips, each of which is a closed subset of $\R^2$ and, hence, of~$B$. Thus $(x_0,y_0)$ belongs to this union of finitely many admissible half-strips and so to $\boldsymbol{\sigma}^{\boldsymbol{\ep}}_{\boldsymbol{\eta},\boldsymbol{\gamma}}$.
\end{proof}

\begin{Corollary}\label{cor:connected} If $(\boldsymbol{\ep},\boldsymbol{\eta}, \boldsymbol{\gamma})$ is admissible for a faithful semitoric system $\is$, then $B \smallsetminus \boldsymbol{\sigma}^{\boldsymbol{\ep}}_{\boldsymbol{\eta},\boldsymbol{\gamma}}$ is contractible.
\end{Corollary}

\begin{proof} By Lemma~\ref{lemma:complement_adm_strip_open} the subset $B \smallsetminus \boldsymbol{\sigma}^{\boldsymbol{\ep}}_{\boldsymbol{\eta},\boldsymbol{\gamma}}$ is open in~$B$. If the intersection of $B \smallsetminus \boldsymbol{\sigma}^{\boldsymbol{\ep}}_{\boldsymbol{\eta},\boldsymbol{\gamma}}$ with every vertical line were either empty or connected, then Proposition~\ref{prop:sub_mostly_vat} would imply that the subsystem relative to $B \smallsetminus \boldsymbol{\sigma}^{\boldsymbol{\ep}}_{\boldsymbol{\eta},\boldsymbol{\gamma}}$ would be faithful semitoric, after which Corollary~\ref{cor:contractible_mom_map} would ensure that $B \smallsetminus \boldsymbol{\sigma}^{\boldsymbol{\ep}}_{\boldsymbol{\eta},\boldsymbol{\gamma}}$, the moment map image of that faithful semitoric subsystem, would be contractible. Thus it suffices to show that the intersection of $B \smallsetminus \boldsymbol{\sigma}^{\boldsymbol{\ep}}_{\boldsymbol{\eta},\boldsymbol{\gamma}}$ with every vertical line is either empty or connected. Fix $x_0 \in \operatorname{pr}_1(B)$; if $x_0 \notin \operatorname{pr}_1(\boldsymbol{\sigma}^{\boldsymbol{\ep}}_{\boldsymbol{\eta},\boldsymbol{\gamma}})$,
 \begin{gather*}
 \big( B \smallsetminus \boldsymbol{\sigma}^{\boldsymbol{\ep}}_{\boldsymbol{\eta},\boldsymbol{\gamma}}\big)
 \cap \{(x,y) \,|\, x=x_0\} = B \cap \{(x,y) \,|\, x=x_0\},
 \end{gather*}
and the result follows from the fact that $\is$ is faithful semitoric. Suppose that $x_0 \in \operatorname{pr}_1(\boldsymbol{\sigma}^{\boldsymbol{\ep}}_{\boldsymbol{\eta},\boldsymbol{\gamma}})$ and call a half-strip $\sigma^{\ep_k}_{\eta_k,\gamma_k}(x_k,y_k)$ {\em maximal} if it is not a proper subset of $\sigma^{\ep_j}_{\eta_j,\gamma_j}(x_j,y_j)$ for any $j \neq k$. Because the triple $(\boldsymbol{\ep},\boldsymbol{\eta}, \boldsymbol{\gamma})$ is admissible for $\is$, there are at most two maximal half-strips $\sigma^{\ep_i}_{\eta_i,\gamma_i}(x_i,y_i)$, $\sigma^{\ep_j}_{\eta_j,\gamma_j}(x_j,y_j)$ with the property that for $s = i,j$, $x_0 \in \big[x_s - \frac{\eta_s}{2}, x_s + \frac{\eta_s}{2}\big]$. These half-strips are, by definition, disjoint and their bases are contained in $\operatorname{Int}(B)$. That property is sufficient to ensure that $\big( B \smallsetminus \boldsymbol{\sigma}^{\boldsymbol{\ep}}_{\boldsymbol{\eta},\boldsymbol{\gamma}}\big) \cap \{(x,y) \,|\, x=x_0\}$ is connected, as desired.
\end{proof}

Finally, we note that admissible triples behave well under isomorphisms of strict complexity one systems and taking saturated subsystems.

\begin{Corollary}\label{cor:adm_triple} Let $(\boldsymbol{\ep},\boldsymbol{\eta}, \boldsymbol{\gamma})$ be an admissible triple for the faithful semitoric system $\is$. Then any faithful semitoric system isomorphic to $\is$ as a~strict complexity one system inherits an admissible triple, as does any subsystem whose image contains every half-strip of $\boldsymbol{\sigma}^{\boldsymbol{\ep}}_{\boldsymbol{\eta},\boldsymbol{\gamma}}$ that it intersects. Moreover, for any cartographic homeomorphism $f_{\boldsymbol{\ep}}$ associated to $\boldsymbol{\ep}$, $(\boldsymbol{\ep},\boldsymbol{\eta}, \boldsymbol{\gamma})$ induces an admissible triple for $f_{\boldsymbol{\ep}}(B)$ relative to the image of the focus-focus values $\{f_{\boldsymbol{\ep}}(c_i)\}_{i\in I}$.
\end{Corollary}

\begin{proof} Let $(M',\omega',\Phi')$ be a faithful semitoric system isomorphic to $\is$ via the isomorphism $(\Psi,\psi)$ as strict complexity one systems (see Definition~\ref{defn:im_faithful_st}). The choices $\boldsymbol{\ep'} = \operatorname{sgn}(\det D \psi) \boldsymbol{\ep}$, $\boldsymbol{\eta'} = \boldsymbol{\eta}$ and $\boldsymbol{\gamma'} = \psi\circ \boldsymbol{\gamma} := \{ \psi \circ \gamma_i\}_{i \in I}$ define an admissible triple for $(M',\omega',\Phi')$ because of the special form of $\psi$ (see Definition~\ref{defn:vat}) and the connectedness of~$B$, ensuring that $\operatorname{sgn}(\det D \psi)$ is constant. Any faithful semitoric system whose image contains the half-strips of~$\boldsymbol{\sigma}^{\boldsymbol{\ep}}_{\boldsymbol{\eta},\boldsymbol{\gamma}}$ that it intersects inherits an admissible triple simply by restriction. Finally, given a~cartographic homeomorphism $f_{\boldsymbol{\ep}}$, the signs $\boldsymbol{\tilde{\ep}} = \operatorname{sgn}(f_{\boldsymbol{\ep}})\boldsymbol{\ep}$, widths $\boldsymbol{\tilde{\eta}} = \boldsymbol{\eta}$ and continuous curves $\boldsymbol{\tilde{\gamma}} = f_{\boldsymbol{\ep}} \circ \boldsymbol{\gamma} := \{ f_{\boldsymbol{\ep}} \circ \gamma_i\}_{i \in I}$ define an admissible triple for $f_{\boldsymbol{\ep}}(B)$ relative to $\{f_{\boldsymbol{\ep}}(c_i)\}_{i \in i}$.
\end{proof}

\subsubsection{Smoothing}\label{sec:smoothing}
With admissible half-strips at hand, we can state and prove the main result of this section.

\begin{Theorem}\label{thm:rect-embedding} Let $\vat$ be a faithful semitoric system and let $\boldsymbol{\ep}$ be a choice of signs such that $B \smallsetminus l^{\boldsymbol{\ep}}$ is connected. $($The existence of one such $\boldsymbol{\ep}$ is guaranteed by Corollary~{\rm \ref{cor:there_exists_choice}.)} Given any cartographic homeomorphism $f_{\boldsymbol{\ep}} \colon B \to \R^2$, there exists a smooth embedding $F_{\boldsymbol{\ep}} \colon B \to \R^2$ of the form
\begin{gather*}
F_{\boldsymbol{\ep}} (x,y) = \big(F^{(1)}_{\boldsymbol{\ep}},F^{(2)}_{\boldsymbol{\ep}}\big) (x,y)= \big(x, F^{(2)}_{\boldsymbol{\ep}}(x,y) \big)
\end{gather*}
agreeing with $f_{\boldsymbol{\ep}}$ on the complement of an arbitrarily small neighborhood of~$l^{\boldsymbol{\ep}}$.
 \end{Theorem}

\begin{figure}[h]\centering
 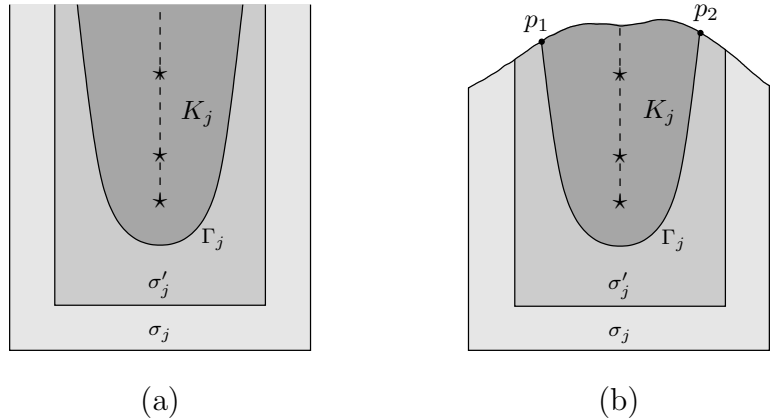
 \caption{In the above figures, $K_j$ is the darkest region, while $L_j$ is the region whose shading is `medium'. Consider the curve $\Ga_j$ as `seam' and glue the diffeomorphisms smoothly along $\Ga_j$. (a) sketches the case $\partial B_0 \cap \si_j' = \varnothing$ and (b) sketches $\partial B_0 \cap \si_j' \neq \varnothing$.} \label{smoothing}
\end{figure}

\begin{proof}Let $\boldsymbol{\ep}$ be as in the statement and fix an admissible triple $(\boldsymbol{\ep},\boldsymbol{\eta}, \boldsymbol{\gamma})$ for~$\is$. The map $f_{\boldsymbol{\ep}}$ is smooth on the complement of the cuts. As in the proof of Corollary~\ref{cor:connected}, say that a~half-strip is maximal if it is not a proper subset of any other half-strip. It is sufficient to modify $f_{\boldsymbol{\ep}}$ in the interior of maximal admissible half-strips, and since maximal half-strips are pairwise disjoint, it suffices to construct the modified map in the interior of each one separately.

Consider a maximal admissible half-strip, say $\si^{\ep_j}_{\eta_j,\gamma_j}(x_j,y_j)$. Without loss of generality assume that $\ep_j = +1$ so as to drop the notational dependence of the half-strip on $\ep_j$. Moreover, fix an admissible triple $(\boldsymbol{\ep},\boldsymbol{\eta}',\boldsymbol{\gamma}')$ for $\vat$, where, if $i \neq j$, $\eta'_i = \eta_i$ and $\gamma'_i = \gamma_i$, and if $i = j$ and $\eta'_j <\eta_j$ and $\gamma_j'(x) > \gamma_j(x)$ whenever both make sense.

There are two cases to consider, namely if $\eta'_j$ can be chosen so that $\partial B \cap \si_{\eta'_j, \gamma'_j} = \varnothing$ or not. Suppose that the former holds; then the set $W_j=B\cap \operatorname{Int}(\si_{\eta'_j,\gamma'_j})$ is open in $\R^2$. The situation is sketched in Fig.~\ref{smoothing}(a). Let $\Gamma_j$ be an embedded curve in $W_j$ of the form $\Gamma_j(x)=(x,h_j(x))$, where $h_j$ is a smooth function, that is disjoint from the cut $l^{\ep_j}$ and is such that $W_j\smallsetminus \Gamma_j$ has two components. Let $K_j$, $L_j$ be the closures in~$W_j$ of the two components of $W_j\smallsetminus \Gamma_j$, so $K_j\cap L_j=\Gamma_j$, and assume without loss of generality that the cut $l^{\ep_j}$ lies in $K_j$.

Recall that the cartographic homeomorphism $f_{\boldsymbol{\ep}}$ preserves orientation and is of the special form
\begin{gather*}
f_{\boldsymbol{\ep}}(x,y)= \big(f_{\boldsymbol{\ep}}^{(1)}(x,y), f_{\boldsymbol{\ep}}^{(2)}(x,y)\big) = \big(x,f_{\boldsymbol{\ep}}^{(2)}(x,y)\big).
\end{gather*}
Therefore, $f\circ\Gamma_j(x)=\big(x,f_{\boldsymbol{\ep}}^{(2)}(x,h_j(x))\big)$ and, since $\ep_j = +1$, if $(x,y)\in K_j$ then $y\ge h_j(x)$ and $f_{\boldsymbol{\ep}}^{(2)}(x,y)\ge f_{\boldsymbol{\ep}}^{(2)}(x,h_j(x))$, as $f_{\boldsymbol{\ep}}$ is orientation-preserving. Define $g_j\colon K_j \to \R^2$ by
\begin{gather*}
g_j(x,y):= \big(x, y + f^{(2)}_{\boldsymbol{\ep}}(x, h_j(x))-h_j(x)\big),
\end{gather*}
which is an orientation-preserving diffeomorphism of $K_j$ onto its image that satisfies
\begin{gather} \label{eq:17}
 g_j(x,h_j(x))= f_{\boldsymbol{\ep}}(x, h_j(x)).
\end{gather}
Now consider the map
\begin{gather*}
 F'_{\ep_j}\colon \ {\rm Int}(\si_{\eta'_j,\gamma'_j}) \cap B \to \R^2, \qquad (x,y) \mapsto \begin{cases}
 f_{\boldsymbol{\ep}}(x,y) & \mbox{if } (x,y) \in L_j \\
 g_j(x,y) & \mbox{if } (x,y) \in K_j,
 \end{cases}
\end{gather*}
which is a homeomorphism onto its image by equation \eqref{eq:17}. Furthermore, because $F'_{\ep_j}$ is a diffeomorphism on the complement of $\Gamma_j$, which is a closed submanifold of $W_j$, $F'_{\ep_j}$ can be isotoped to be a diffeomorphism onto the image $F'_{\ep_j}(W_j)$ via an isotopy that is supported in an arbitrarily small neighborhood of $\Gamma_j$ and is the identity in $L_j$ (cf.\ Hirsch \cite[Chapter~8]{hirsch}). By construction, $F'_{\ep_j}$ extends to all of $\sigma_{\eta_j,\gamma_j}$ as a~diffeomorphism, say $F_{\ep_j}$, on $\si_{\eta_j,\gamma_j}$ that agrees with~$f_{\boldsymbol{\ep}}$ on $\si_{\eta_j, \gamma_j}\smallsetminus K_j$. The map $F_{\ep_j}$ is the desired smoothing.

It remains to consider the case in which an admissible triple $(\boldsymbol{\ep},\boldsymbol{\eta}',\boldsymbol{\gamma}')$ as above does not exist, i.e., for any choices of~$\eta_j'$ and~$\gamma_j'$ as above, the corresponding half-strip also intersects $\partial B$ (see Fig.~\ref{smoothing}(b)). In this case, modify the argument as follows. Let $\Gamma_j\subset B\cap {\rm Int}(\si_{\eta_j,\gamma_j})$ be chosen as above, and so that all boundary points of $\Gamma_j$ also lie in path-connected components of $\partial B_0\cap{\rm Int}(\si_{\eta_j,\gamma_j})$. Because $f_{\boldsymbol{\ep}}$ is, by definition, smooth at a boundary point~$p$ of~$\Gamma_j$, the map $f_{\boldsymbol{\ep}}$ and the smooth curve $\Gamma_j$ can be extended to a neighborhood of~$p\in\R^2$. Make such an extension near the one or two boundary points of $\Gamma$, and let $W_j$ be an open tubular neighborhood of the extended curve~$\Gamma_j$. Let $K_j$ and $L_j$ be defined as in the first case (enlarged as per the extension just described), with the map~$F'_{\ep_j}$ defined as above. But to apply the smoothing argument, restrict attention to $K_j\cap W_j$ and $L_j\cap W_j$ so that $\Gamma_j$ is a closed submanifold of an open manifold, in this case the tubular neighborhood~$W_j$.
\end{proof}

\begin{Definition} Given an admissible triple $(\boldsymbol{\ep},\boldsymbol{\eta},\boldsymbol{\gamma})$, a~map $F_{\boldsymbol{\ep}}$ as constructed in Theorem~\ref{thm:rect-embedding} is referred to as an {\em $\boldsymbol{\eta}$-cartographic embedding}. Moreover, if the dependence on $\boldsymbol{\eta}$ is to be remembered, an $\boldsymbol{\eta}$-cartographic embedding is denoted by $F_{\boldsymbol{\ep},\boldsymbol{\eta}}$.
\end{Definition}

Theorem \ref{thm:rect-embedding} motivates introducing the following notion.

\begin{Definition}\label{defn:eta-carto} A faithful semitoric system is {\em$\boldsymbol{\eta}$-cartographic} if it admits an admissible triple $(\boldsymbol{\ep},\boldsymbol{\eta},\boldsymbol{\gamma})$ and a cartographic homeomorphism $f_{\boldsymbol{\ep}}$ whose restriction to the complement of the union of the corresponding admissible half-strips is the identity. If the choice of $(\boldsymbol{\ep},\boldsymbol{\eta},\boldsymbol{\gamma})$ is to be remembered, the system is said to be $\boldsymbol{\eta}$-cartographic {\em with respect to $(\boldsymbol{\ep},\boldsymbol{\eta},\boldsymbol{\gamma})$}.
\end{Definition}

The first application of Theorem \ref{thm:rect-embedding} is the following result.

\begin{Theorem}\label{thm:eta-carto} Any faithful semitoric system $(M,\omega,\Phi)$ is isomorphic, as a semitoric system, to an $\boldsymbol{\eta}$-cartographic one.
\end{Theorem}

\begin{proof} By Corollary \ref{cor:existence_admissible}, $\is$ has an admissible triple. Fix one such triple $(\boldsymbol{\ep},\boldsymbol{\eta},\boldsymbol{\gamma})$ for $\is$ and let $f_{\boldsymbol{\ep}} \colon B \to \R^2$ be the cartographic homeomorphism associated to $\boldsymbol{\ep}$. Let $F_{\boldsymbol{\ep}}\colon B \to \R^2$ be the associated smooth $\boldsymbol{\eta}$-cartographic embedding constructed in the proof of Theorem~\ref{thm:rect-embedding}. The form of $F_{\boldsymbol{\ep}}$ implies that, by construction, $(M,\omega,F_{\boldsymbol{\ep}}\circ \Phi)$ is a faithful semitoric system isomorphic to~$\is$. Moreover, since $F_{\boldsymbol{\ep}}$ is orientation-preserving, $(M,\omega,F_{\boldsymbol{\ep}}\circ \Phi)$ inherits an admissible triple $(\boldsymbol{\ep},\boldsymbol{\eta},F_{\boldsymbol{\ep}}\circ \boldsymbol{\gamma})$ by Corollary~\ref{cor:adm_triple}. The map $f_{\boldsymbol{\ep}} \circ F_{\boldsymbol{\ep}}^{-1}\colon F_{\boldsymbol{\ep}}(B) \to \R^2$ is a cartographic homeomorphism for $(M,\omega,F_{\boldsymbol{\ep}}\circ \Phi)$, which, by definition of $F_{\boldsymbol{\ep}}$, is the identity on the complement of the admissible half-strips for $(M,\omega,F_{\boldsymbol{\ep}}\circ \Phi)$ corresponding to $(\boldsymbol{\ep},\boldsymbol{\eta},F_{\boldsymbol{\ep}}\circ \boldsymbol{\gamma})$. Therefore, $(M,\omega,F_{\boldsymbol{\ep}}\circ \Phi)$ is $\boldsymbol{\eta}$-cartographic as required.
\end{proof}

To conclude this section, we show that the image of a cartographic homeomorphism can be seen as a `limit' of the moment map images of $\boldsymbol{\eta}$-cartographic systems. To make the above precise, let $\is$ be a faithful semitoric system, fix a choice of signs $\boldsymbol{\ep}$ for which $S_{\boldsymbol{\ep}}$ is path-connected (which exists by Corollary~\ref{cor:there_exists_choice}), and fix a~cartographic homeomorphism $f_{\boldsymbol{\ep}}$ associated to $\boldsymbol{\ep}$. Consider the set consisting of quadruples $(\boldsymbol{\ep}, \boldsymbol{\eta}, \boldsymbol{\gamma}, F_{\boldsymbol{\ep},\boldsymbol{\eta}})$, where $(\boldsymbol{\ep}, \boldsymbol{\eta}, \boldsymbol{\gamma})$ is an admissible triple for $\vat$ and $F_{\boldsymbol{\ep},\boldsymbol{\eta}}$ is an $\boldsymbol{\eta}$-cartographic embedding constructed starting from the cartographic homeomorphism $f_{\boldsymbol{\ep}}$. On this set, we define a~partial order $\succeq$ by setting
\begin{gather*}
(\boldsymbol{\ep}, \boldsymbol{\eta}, \boldsymbol{\ga},F_{\boldsymbol{\ep},\boldsymbol{\eta}} ) \succeq (\boldsymbol{\ep}, \boldsymbol{\tilde{\eta}}, \boldsymbol{\tilde{\ga}}, F_{\boldsymbol{\ep},\boldsymbol{\tilde{\eta}}}),
\end{gather*}
if and only if $\eta_j \geq \etti_j$ for all $j$, and $F_{\boldsymbol{\ep}, \boldsymbol{\eta}}(\boldsymbol{\si}_{\boldsymbol{\eta}, \boldsymbol{\ga}}^{\boldsymbol{\ep}}) \supseteq F_{\boldsymbol{\ep}, \boldsymbol{\tilde{\eta}}}(\boldsymbol{\si}_{\boldsymbol{\tilde{\eta}}, \boldsymbol{\tilde{\ga}}}^{\boldsymbol{\ep}})$. Note that $\succeq$ is reflexive and transitive. Moreover, for any two elements $(\boldsymbol{\ep}, \boldsymbol{\eta}, \boldsymbol{\ga},F_{\boldsymbol{\ep},\boldsymbol{\eta}}) $ and $\big(\boldsymbol{\ep}, \boldsymbol{\tilde{\eta}}, \boldsymbol{\tilde{\ga}}, F_{\boldsymbol{\ep},\boldsymbol{\tilde{\eta}}}\big)$, there exists an element $\big(\boldsymbol{\ep}, \boldsymbol{\tilde{\tilde{\eta}}}, \boldsymbol{\tilde{\tilde{\ga}}}, F_{\boldsymbol{\ep},\boldsymbol{\tilde{\tilde{\eta}}}}\big)$ such that $(\boldsymbol{\ep}, \boldsymbol{\eta}, \boldsymbol{\ga},F_{\boldsymbol{\ep},\boldsymbol{\eta}} ) \succeq \big(\boldsymbol{\ep}, \boldsymbol{\tilde{\tilde{\eta}}}, \boldsymbol{\tilde{\tilde{\ga}}}, F_{\boldsymbol{\ep},\boldsymbol{\tilde{\tilde{\eta}}}}\big) $ and $\big(\boldsymbol{\ep}, \boldsymbol{\tilde{\eta}}, \boldsymbol{\tilde{\ga}}, F_{\boldsymbol{\ep},\boldsymbol{\tilde{\eta}}}\big)\succeq \big(\boldsymbol{\ep}, \boldsymbol{\tilde{\tilde{\eta}}}, \boldsymbol{\tilde{\tilde{\ga}}}, F_{\boldsymbol{\ep},\boldsymbol{\tilde{\tilde{\eta}}}}\big)$,
as the construction of $F_{\boldsymbol{\ep}, \boldsymbol{\eta}}$ and $F_{\boldsymbol{\ep}, \boldsymbol{\tilde{\boldsymbol{\eta}}}}$ shows. Thus $\succeq$ turns the set of quadruples $(\boldsymbol{\ep}, \boldsymbol{\eta}, \boldsymbol{\gamma}, F_{\boldsymbol{\ep},\boldsymbol{\eta}})$ into a~directed set.

\begin{Proposition}\label{prop:limit} A cartographic moment map image of a faithful semitoric system associated to a choice of signs whose corresponding cuts do not disconnect the moment map image is the direct limit of $\boldsymbol{\eta}$-cartographic moment map images.
\end{Proposition}

\begin{proof} As above, fix a faithful semitoric system $\vat$, a choice of signs $\boldsymbol{\ep}$ ma\-king~$S_{\boldsymbol{\ep}}$ path-connected, and a cartographic homeomorphism $f_{\boldsymbol{\ep}}$ associated to $\boldsymbol{\ep}$. The set of quadruples $(\boldsymbol{\ep}, \boldsymbol{\eta}, \boldsymbol{\ga},F_{\boldsymbol{\ep},\boldsymbol{\eta}})$ is a~directed set with the above partial order $\succeq$. The condition
 \begin{gather*}
(\boldsymbol{\ep}, \boldsymbol{\eta}, \boldsymbol{\ga},F_{\boldsymbol{\ep},\boldsymbol{\eta}} )
 \succeq \big(\boldsymbol{\ep}, \boldsymbol{\tilde{\eta}}, \boldsymbol{\tilde{\ga}}, F_{\boldsymbol{\ep},\boldsymbol{\tilde{\eta}}}\big)
 \end{gather*}
implies that
\begin{gather} \label{inclusions}
 F_{\boldsymbol{\ep}, \boldsymbol{\eta}}\big(B \smallsetminus {\boldsymbol{\si}}_{\boldsymbol{\eta}, {\boldsymbol{\ga}}}^{\boldsymbol{\ep}}\big) \subseteq F_{\boldsymbol{\ep}, \tilde{\boldsymbol{\eta}}}\big(B \smallsetminus {\boldsymbol{\si}}_{\tilde{\boldsymbol{\eta}}, \boldsymbol{\tilde{\ga}}}^{\boldsymbol{\ep}}\big)
\end{gather}
such that $f_{\boldsymbol{\ep}}(B \smallsetminus l^{\boldsymbol{\ep}})$ coincides with the direct limit in the category of topological spaces given~by
\begin{gather*}
 \lim_{\longrightarrow} F_{\boldsymbol{\ep}, \boldsymbol{\eta}}\big(B \smallsetminus{\boldsymbol{\si}}_{\boldsymbol{\eta}, {\boldsymbol{\ga}}}^{\boldsymbol{\ep}} \big) = \left. \left( \bigsqcup_{\boldsymbol{\eta}} F_{\boldsymbol{\ep},
 \boldsymbol{\eta}}\big(B \smallsetminus {\boldsymbol{\si}}_{\boldsymbol{\eta},
 {\boldsymbol{\ga}}}^{\boldsymbol{\ep}}\big) \right) \right / \sim,
\end{gather*}
where $z \in F_{\boldsymbol{\ep}, \boldsymbol{\eta}}\big(B \smallsetminus{\boldsymbol{\si}}_{\boldsymbol{\eta}, {\boldsymbol{\ga}}}^{\boldsymbol{\ep}} \big) \sim \zti \in F_{\boldsymbol{\ep}, \boldsymbol{\tilde{\eta}}} \big(B \smallsetminus{\boldsymbol{\si}}_{\boldsymbol{\tilde{\eta}}, {\boldsymbol{\tilde{\ga}}}}^{\boldsymbol{\ep}} \big) $ if $z$ and $\zti$ get mapped under the corresponding inclusions in~\eqref{inclusions} to the same point.
\end{proof}

\subsection*{Acknowledgments}
The authors would like to thank the anonymous referees for their careful reading of an earlier draft of this paper, as well as for their suggestions that have significantly improved the quality of the paper. Furthermore, the authors would like to thank Eva Miranda for providing comments on an earlier draft of this paper. S.H.\ was partially supported by the Research Fund of the University of Antwerp and by SwissMAP. S.S.\ was partially supported by SFB-TRR 191 “Symplectic Structures in Geometry, Algebra and Dynamics”, funded by the Deutsche Forschungsgemeinschaft. D.S.\ was partially supported by the University of Cologne, SwissMAP, the NWO Veni grant 639.031.345 and by the CNPq Universal grant 409552/2016-0. M.S.\ was partially supported by Mercer University, the Institute of Pure and Applied Mathematics (IMPA) in Rio de Janeiro, the University of Cologne, and the Swiss Federal Institute of Technology (ETH) in Zurich.

\LastPageEnding

\end{document}

%% file: torusHeightFct.pdf_t
\begin{picture}(0,0)%
\includegraphics{torusHeightFct.pdf}%
\end{picture}%
\setlength{\unitlength}{4144sp}%
\begingroup\makeatletter\ifx\SetFigFont\undefined%
\gdef\SetFigFont#1#2#3#4#5{%
  \reset@font\fontsize{#1}{#2pt}%
  \fontfamily{#3}\fontseries{#4}\fontshape{#5}%
  \selectfont}%
\fi\endgroup%
\begin{picture}(1910,2544)(-551,-1693)
\put(-494,-421){\makebox(0,0)[lb]{\smash{{\SetFigFont{10}{12.0}{\rmdefault}{\mddefault}{\updefault}{\color[rgb]{0,0,0}$h$}%
}}}}
\end{picture}%

%% file: annulus-NEW.pdf_t
\begin{picture}(0,0)%
\includegraphics{annulus-NEW.pdf}%
\end{picture}%
\setlength{\unitlength}{4144sp}%
\begingroup\makeatletter\ifx\SetFigFont\undefined%
\gdef\SetFigFont#1#2#3#4#5{%
  \reset@font\fontsize{#1}{#2pt}%
  \fontfamily{#3}\fontseries{#4}\fontshape{#5}%
  \selectfont}%
\fi\endgroup%
\begin{picture}(1954,2096)(1253,-1327)
\end{picture}%

%% file: doubleCover.pdf_t
\begin{picture}(0,0)%
\includegraphics{doubleCover.pdf}%
\end{picture}%
\setlength{\unitlength}{4144sp}%
\begingroup\makeatletter\ifx\SetFigFont\undefined%
\gdef\SetFigFont#1#2#3#4#5{%
  \reset@font\fontsize{#1}{#2pt}%
  \fontfamily{#3}\fontseries{#4}\fontshape{#5}%
  \selectfont}%
\fi\endgroup%
\begin{picture}(3363,2544)(79,-1873)
\put(2746,-151){\makebox(0,0)[b]{\smash{{\SetFigFont{10}{12.0}{\rmdefault}{\mddefault}{\updefault}{\color[rgb]{0,0,0}$U$}%
}}}}
\end{picture}%

%% file: singLocNF.pdf_t
\begin{picture}(0,0)%
\includegraphics{singLocNF.pdf}%
\end{picture}%
\setlength{\unitlength}{4144sp}%
\begingroup\makeatletter\ifx\SetFigFont\undefined%
\gdef\SetFigFont#1#2#3#4#5{%
  \reset@font\fontsize{#1}{#2pt}%
  \fontfamily{#3}\fontseries{#4}\fontshape{#5}%
  \selectfont}%
\fi\endgroup%
\begin{picture}(3894,1824)(-1091,-973)
\end{picture}%

%% file: multipleFF.pdf_t
\begin{picture}(0,0)%
\includegraphics{multipleFF.pdf}%
\end{picture}%
\setlength{\unitlength}{4144sp}%
\begingroup\makeatletter\ifx\SetFigFont\undefined%
\gdef\SetFigFont#1#2#3#4#5{%
  \reset@font\fontsize{#1}{#2pt}%
  \fontfamily{#3}\fontseries{#4}\fontshape{#5}%
  \selectfont}%
\fi\endgroup%
\begin{picture}(2632,2301)(-1823,-1594)
\end{picture}%

%% file: spiralFlow1FF.pdf_t
\begin{picture}(0,0)%
\includegraphics{spiralFlow1FF.pdf}%
\end{picture}%
\setlength{\unitlength}{4144sp}%
\begingroup\makeatletter\ifx\SetFigFont\undefined%
\gdef\SetFigFont#1#2#3#4#5{%
  \reset@font\fontsize{#1}{#2pt}%
  \fontfamily{#3}\fontseries{#4}\fontshape{#5}%
  \selectfont}%
\fi\endgroup%
\begin{picture}(2389,2198)(200,-1587)
\end{picture}%

%% file: spiralFlow2FF.pdf_t
\begin{picture}(0,0)%
\includegraphics{spiralFlow2FF.pdf}%
\end{picture}%
\setlength{\unitlength}{4144sp}%
\begingroup\makeatletter\ifx\SetFigFont\undefined%
\gdef\SetFigFont#1#2#3#4#5{%
  \reset@font\fontsize{#1}{#2pt}%
  \fontfamily{#3}\fontseries{#4}\fontshape{#5}%
  \selectfont}%
\fi\endgroup%
\begin{picture}(1798,1958)(475,-1559)
\end{picture}%

%% file: image.pdf_t
\begin{picture}(0,0)%
\includegraphics{image.pdf}%
\end{picture}%
\setlength{\unitlength}{4144sp}%
\begingroup\makeatletter\ifx\SetFigFont\undefined%
\gdef\SetFigFont#1#2#3#4#5{%
  \reset@font\fontsize{#1}{#2pt}%
  \fontfamily{#3}\fontseries{#4}\fontshape{#5}%
  \selectfont}%
\fi\endgroup%
\begin{picture}(3769,2857)(429,-2636)
\put(721,-826){\makebox(0,0)[b]{\smash{{\SetFigFont{12}{14.4}{\rmdefault}{\mddefault}{\updefault}{\color[rgb]{0,0,0}$\star$}%
}}}}
\put(1081,-916){\makebox(0,0)[b]{\smash{{\SetFigFont{12}{14.4}{\rmdefault}{\mddefault}{\updefault}{\color[rgb]{0,0,0}$\star$}%
}}}}
\put(1981,-736){\makebox(0,0)[b]{\smash{{\SetFigFont{12}{14.4}{\rmdefault}{\mddefault}{\updefault}{\color[rgb]{0,0,0}$\star$}%
}}}}
\put(2341,-916){\makebox(0,0)[b]{\smash{{\SetFigFont{12}{14.4}{\rmdefault}{\mddefault}{\updefault}{\color[rgb]{0,0,0}$\star$}%
}}}}
\put(1981,-871){\makebox(0,0)[b]{\smash{{\SetFigFont{12}{14.4}{\rmdefault}{\mddefault}{\updefault}{\color[rgb]{0,0,0}$\star$}%
}}}}
\put(721,-1501){\makebox(0,0)[b]{\smash{{\SetFigFont{12}{14.4}{\rmdefault}{\mddefault}{\updefault}{\color[rgb]{0,0,0}$\star$}%
}}}}
\put(1441,-2446){\makebox(0,0)[b]{\smash{{\SetFigFont{12}{14.4}{\rmdefault}{\mddefault}{\updefault}{\color[rgb]{0,0,0}$\star$}%
}}}}
\put(1441,-1951){\makebox(0,0)[b]{\smash{{\SetFigFont{12}{14.4}{\rmdefault}{\mddefault}{\updefault}{\color[rgb]{0,0,0}$\star$}%
}}}}
\put(1441,-1771){\makebox(0,0)[b]{\smash{{\SetFigFont{12}{14.4}{\rmdefault}{\mddefault}{\updefault}{\color[rgb]{0,0,0}$\star$}%
}}}}
\put(3016,-1681){\makebox(0,0)[b]{\smash{{\SetFigFont{12}{14.4}{\rmdefault}{\mddefault}{\updefault}{\color[rgb]{0,0,0}$\star$}%
}}}}
\put(3331,-1276){\makebox(0,0)[b]{\smash{{\SetFigFont{12}{14.4}{\rmdefault}{\mddefault}{\updefault}{\color[rgb]{0,0,0}$\star$}%
}}}}
\put(3307,-2488){\makebox(0,0)[lb]{\smash{{\SetFigFont{10}{12.0}{\rmdefault}{\mddefault}{\updefault}{\color[rgb]{0,0,0}$\mathsf x_{\sup}$}%
}}}}
\put(3961,-2581){\makebox(0,0)[b]{\smash{{\SetFigFont{12}{14.4}{\rmdefault}{\mddefault}{\updefault}{\color[rgb]{0,0,0}$\R$}%
}}}}
\end{picture}%

%% file: cut.pdf_t
\begin{picture}(0,0)%
\includegraphics{cut.pdf}%
\end{picture}%
\setlength{\unitlength}{4144sp}%
\begingroup\makeatletter\ifx\SetFigFont\undefined%
\gdef\SetFigFont#1#2#3#4#5{%
  \reset@font\fontsize{#1}{#2pt}%
  \fontfamily{#3}\fontseries{#4}\fontshape{#5}%
  \selectfont}%
\fi\endgroup%
\begin{picture}(4829,1945)(154,-1925)
\put(1001,-553){\makebox(0,0)[b]{\smash{{\SetFigFont{7}{8.4}{\rmdefault}{\mddefault}{\updefault}{\color[rgb]{0,0,0}$\ep_i=+1$}%
}}}}
\put(997,-1051){\makebox(0,0)[b]{\smash{{\SetFigFont{7}{8.4}{\rmdefault}{\mddefault}{\updefault}{\color[rgb]{0,0,0}$\ep_j=+1$}%
}}}}
\put(2848,-1051){\makebox(0,0)[b]{\smash{{\SetFigFont{8}{9.6}{\rmdefault}{\mddefault}{\updefault}{\color[rgb]{0,0,0}$\ep_j=-1$}%
}}}}
\put(2849,-554){\makebox(0,0)[b]{\smash{{\SetFigFont{8}{9.6}{\rmdefault}{\mddefault}{\updefault}{\color[rgb]{0,0,0}$\ep_i=+1$}%
}}}}
\put(4696,-1049){\makebox(0,0)[b]{\smash{{\SetFigFont{8}{9.6}{\rmdefault}{\mddefault}{\updefault}{\color[rgb]{0,0,0}$\ep_j=+1$}%
}}}}
\put(4697,-556){\makebox(0,0)[b]{\smash{{\SetFigFont{8}{9.6}{\rmdefault}{\mddefault}{\updefault}{\color[rgb]{0,0,0}$\ep_i=-1$}%
}}}}
\put(721,-556){\makebox(0,0)[b]{\smash{{\SetFigFont{12}{14.4}{\rmdefault}{\mddefault}{\updefault}{\color[rgb]{0,0,0}$\star$}%
}}}}
\put(721,-1051){\makebox(0,0)[b]{\smash{{\SetFigFont{12}{14.4}{\rmdefault}{\mddefault}{\updefault}{\color[rgb]{0,0,0}$\star$}%
}}}}
\put(2566,-556){\makebox(0,0)[b]{\smash{{\SetFigFont{12}{14.4}{\rmdefault}{\mddefault}{\updefault}{\color[rgb]{0,0,0}$\star$}%
}}}}
\put(2566,-1051){\makebox(0,0)[b]{\smash{{\SetFigFont{12}{14.4}{\rmdefault}{\mddefault}{\updefault}{\color[rgb]{0,0,0}$\star$}%
}}}}
\put(4411,-556){\makebox(0,0)[b]{\smash{{\SetFigFont{12}{14.4}{\rmdefault}{\mddefault}{\updefault}{\color[rgb]{0,0,0}$\star$}%
}}}}
\put(4411,-1051){\makebox(0,0)[b]{\smash{{\SetFigFont{12}{14.4}{\rmdefault}{\mddefault}{\updefault}{\color[rgb]{0,0,0}$\star$}%
}}}}
\put(676,-1861){\makebox(0,0)[b]{\smash{{\SetFigFont{12}{14.4}{\rmdefault}{\mddefault}{\updefault}{\color[rgb]{0,0,0}(a)}%
}}}}
\put(2566,-1861){\makebox(0,0)[b]{\smash{{\SetFigFont{12}{14.4}{\rmdefault}{\mddefault}{\updefault}{\color[rgb]{0,0,0}(b)}%
}}}}
\put(4456,-1861){\makebox(0,0)[b]{\smash{{\SetFigFont{12}{14.4}{\rmdefault}{\mddefault}{\updefault}{\color[rgb]{0,0,0}(c)}%
}}}}
\end{picture}%

%% file: doubleCutFF.pdf_t
\begin{picture}(0,0)%
\includegraphics{doubleCutFF.pdf}%
\end{picture}%
\setlength{\unitlength}{4144sp}%
\begingroup\makeatletter\ifx\SetFigFont\undefined%
\gdef\SetFigFont#1#2#3#4#5{%
  \reset@font\fontsize{#1}{#2pt}%
  \fontfamily{#3}\fontseries{#4}\fontshape{#5}%
  \selectfont}%
\fi\endgroup%
\begin{picture}(1665,1421)(1519,-1116)
\put(2339,-441){\makebox(0,0)[b]{\smash{{\SetFigFont{12}{14.4}{\rmdefault}{\mddefault}{\updefault}{\color[rgb]{0,0,0}$\star$}%
}}}}
\end{picture}%

%% file: halfStrips.pdf_t
\begin{picture}(0,0)%
\includegraphics{halfStrips.pdf}%
\end{picture}%
\setlength{\unitlength}{4144sp}%
\begingroup\makeatletter\ifx\SetFigFont\undefined%
\gdef\SetFigFont#1#2#3#4#5{%
  \reset@font\fontsize{#1}{#2pt}%
  \fontfamily{#3}\fontseries{#4}\fontshape{#5}%
  \selectfont}%
\fi\endgroup%
\begin{picture}(4814,2146)(159,-1970)
\put(2863,-546){\makebox(0,0)[b]{\smash{{\SetFigFont{8}{9.6}{\rmdefault}{\mddefault}{\updefault}{\color[rgb]{0,0,0}$\ep_i=+1$}%
}}}}
\put(1003,-1049){\makebox(0,0)[b]{\smash{{\SetFigFont{7}{8.4}{\rmdefault}{\mddefault}{\updefault}{\color[rgb]{0,0,0}$\ep_j=+1$}%
}}}}
\put(1004,-543){\makebox(0,0)[b]{\smash{{\SetFigFont{7}{8.4}{\rmdefault}{\mddefault}{\updefault}{\color[rgb]{0,0,0}$\ep_i=+1$}%
}}}}
\put(2856,-1065){\makebox(0,0)[b]{\smash{{\SetFigFont{8}{9.6}{\rmdefault}{\mddefault}{\updefault}{\color[rgb]{0,0,0}$\ep_j=-1$}%
}}}}
\put(4699,-1036){\makebox(0,0)[b]{\smash{{\SetFigFont{8}{9.6}{\rmdefault}{\mddefault}{\updefault}{\color[rgb]{0,0,0}$\ep_j=+1$}%
}}}}
\put(4704,-558){\makebox(0,0)[b]{\smash{{\SetFigFont{8}{9.6}{\rmdefault}{\mddefault}{\updefault}{\color[rgb]{0,0,0}$\ep_i=-1$}%
}}}}
\put(2566,-556){\makebox(0,0)[b]{\smash{{\SetFigFont{12}{14.4}{\rmdefault}{\mddefault}{\updefault}{\color[rgb]{0,0,0}$\star$}%
}}}}
\put(2566,-1051){\makebox(0,0)[b]{\smash{{\SetFigFont{12}{14.4}{\rmdefault}{\mddefault}{\updefault}{\color[rgb]{0,0,0}$\star$}%
}}}}
\put(4411,-556){\makebox(0,0)[b]{\smash{{\SetFigFont{12}{14.4}{\rmdefault}{\mddefault}{\updefault}{\color[rgb]{0,0,0}$\star$}%
}}}}
\put(4411,-1051){\makebox(0,0)[b]{\smash{{\SetFigFont{12}{14.4}{\rmdefault}{\mddefault}{\updefault}{\color[rgb]{0,0,0}$\star$}%
}}}}
\put(2566,-1906){\makebox(0,0)[b]{\smash{{\SetFigFont{12}{14.4}{\rmdefault}{\mddefault}{\updefault}{\color[rgb]{0,0,0}(b)}%
}}}}
\put(4411,-1906){\makebox(0,0)[b]{\smash{{\SetFigFont{12}{14.4}{\rmdefault}{\mddefault}{\updefault}{\color[rgb]{0,0,0}(c)}%
}}}}
\put(721,-1906){\makebox(0,0)[b]{\smash{{\SetFigFont{12}{14.4}{\rmdefault}{\mddefault}{\updefault}{\color[rgb]{0,0,0}(a)}%
}}}}
\put(721,-1051){\makebox(0,0)[b]{\smash{{\SetFigFont{12}{14.4}{\rmdefault}{\mddefault}{\updefault}{\color[rgb]{0,0,0}$\star$}%
}}}}
\put(721,-556){\makebox(0,0)[b]{\smash{{\SetFigFont{12}{14.4}{\rmdefault}{\mddefault}{\updefault}{\color[rgb]{0,0,0}$\star$}%
}}}}
\end{picture}%

%% file: smoothing.pdf_t
\begin{picture}(0,0)%
\includegraphics{smoothing.pdf}%
\end{picture}%
\setlength{\unitlength}{4144sp}%
\begingroup\makeatletter\ifx\SetFigFont\undefined%
\gdef\SetFigFont#1#2#3#4#5{%
  \reset@font\fontsize{#1}{#2pt}%
  \fontfamily{#3}\fontseries{#4}\fontshape{#5}%
  \selectfont}%
\fi\endgroup%
\begin{picture}(4701,2665)(439,-2555)
\put(1351,-2491){\makebox(0,0)[b]{\smash{{\SetFigFont{12}{14.4}{\rmdefault}{\mddefault}{\updefault}{\color[rgb]{0,0,0}(a)}%
}}}}
\put(4096,-2491){\makebox(0,0)[b]{\smash{{\SetFigFont{12}{14.4}{\rmdefault}{\mddefault}{\updefault}{\color[rgb]{0,0,0}(b)}%
}}}}
\put(1351,-2041){\makebox(0,0)[b]{\smash{{\SetFigFont{8}{9.6}{\rmdefault}{\mddefault}{\updefault}{\color[rgb]{0,0,0}$\si_j$}%
}}}}
\put(1348,-1756){\makebox(0,0)[b]{\smash{{\SetFigFont{8}{9.6}{\rmdefault}{\mddefault}{\updefault}{\color[rgb]{0,0,0}$\si_j'$}%
}}}}
\put(4103,-2047){\makebox(0,0)[b]{\smash{{\SetFigFont{8}{9.6}{\rmdefault}{\mddefault}{\updefault}{\color[rgb]{0,0,0}$\si_j$}%
}}}}
\put(4100,-1762){\makebox(0,0)[b]{\smash{{\SetFigFont{8}{9.6}{\rmdefault}{\mddefault}{\updefault}{\color[rgb]{0,0,0}$\si_j'$}%
}}}}
\put(1599,-1490){\makebox(0,0)[lb]{\smash{{\SetFigFont{8}{9.6}{\rmdefault}{\mddefault}{\updefault}{\color[rgb]{0,0,0}$\Ga_j$}%
}}}}
\put(4351,-1496){\makebox(0,0)[lb]{\smash{{\SetFigFont{8}{9.6}{\rmdefault}{\mddefault}{\updefault}{\color[rgb]{0,0,0}$\Ga_j$}%
}}}}
\put(3606,-196){\makebox(0,0)[b]{\smash{{\SetFigFont{10}{12.0}{\rmdefault}{\mddefault}{\updefault}{\color[rgb]{0,0,0}$p_1$}%
}}}}
\put(4620,-132){\makebox(0,0)[b]{\smash{{\SetFigFont{10}{12.0}{\rmdefault}{\mddefault}{\updefault}{\color[rgb]{0,0,0}$p_2$}%
}}}}
\put(4103,-1282){\makebox(0,0)[b]{\smash{{\SetFigFont{11}{13.2}{\rmdefault}{\mddefault}{\updefault}{\color[rgb]{0,0,0}$\star$}%
}}}}
\put(4103,-1012){\makebox(0,0)[b]{\smash{{\SetFigFont{11}{13.2}{\rmdefault}{\mddefault}{\updefault}{\color[rgb]{0,0,0}$\star$}%
}}}}
\put(1351,-1276){\makebox(0,0)[b]{\smash{{\SetFigFont{11}{13.2}{\rmdefault}{\mddefault}{\updefault}{\color[rgb]{0,0,0}$\star$}%
}}}}
\put(1351,-1006){\makebox(0,0)[b]{\smash{{\SetFigFont{11}{13.2}{\rmdefault}{\mddefault}{\updefault}{\color[rgb]{0,0,0}$\star$}%
}}}}
\put(1351,-511){\makebox(0,0)[b]{\smash{{\SetFigFont{11}{13.2}{\rmdefault}{\mddefault}{\updefault}{\color[rgb]{0,0,0}$\star$}%
}}}}
\put(4103,-517){\makebox(0,0)[b]{\smash{{\SetFigFont{11}{13.2}{\rmdefault}{\mddefault}{\updefault}{\color[rgb]{0,0,0}$\star$}%
}}}}
\put(4339,-741){\makebox(0,0)[b]{\smash{{\SetFigFont{10}{12.0}{\rmdefault}{\mddefault}{\updefault}{\color[rgb]{0,0,0}$K_j$}%
}}}}
\put(1576,-745){\makebox(0,0)[b]{\smash{{\SetFigFont{10}{12.0}{\rmdefault}{\mddefault}{\updefault}{\color[rgb]{0,0,0}$K_j$}%
}}}}
\end{picture}%